%% file: TRUnderlyingProb.tex
\documentclass[10pt]{article}

\usepackage[dvips]{graphics,epsfig} 
\usepackage{amssymb,amsmath,amsthm,mathrsfs} 
\usepackage{graphicx}
\usepackage{lineno}

\usepackage[authoryear,sort&compress]{natbib}
\usepackage{url}
\usepackage{enumerate}
\usepackage{colordvi,color,pspicture}
\usepackage{psfrag}
\usepackage{rotating}
\usepackage{verbatim}

\graphicspath{{Figures/}}

\newcommand{\I}{\mathbf{I}}
\newcommand{\la}{\text{and}}
\newcommand{\lo}{\text{or}}

\newcommand{\E}{\mathbf{E}}

\newcommand{\U}{\mathcal{U}}

\newcommand{\g}{\gamma}
\newcommand{\G}{\Gamma}
\newcommand{\V}{\mathcal{V}}
\newcommand{\A}{\mathcal{A}}
\newcommand{\Y}{\mathcal{Y}}
\newcommand{\y}{\mathsf{y}}
\newcommand{\X}{\mathcal{X}}

\newcommand{\mE}{\mathcal{E}}
\newcommand{\mI}{\mathcal{I}}
\newcommand{\NY}{N_{\Y}}

\newcommand{\NPE}{N_{PE}}
\newcommand{\TY}{T(\Y_3)}
\newcommand{\N}{\mathcal{N}}
\newcommand{\R}{\mathbb{R}}

\newcommand{\msP}{\mathscr{P}}

\newcommand{\Var}{\mathbf{Var}}
\newcommand{\Cov}{\mathbf{Cov}}

\DeclareMathOperator{\argsup}{argsup}

\theoremstyle{plain}
\newtheorem{theorem}{Theorem}[section]
\newtheorem{lemma}[theorem]{Lemma}

\newtheorem{corollary}[theorem]{Corollary}

\theoremstyle{definition}

\theoremstyle{remark}

\newtheorem{remark}[theorem]{Remark}

\begin{document}

\title{Technical Report \# KU-EC-13-3:\\
Distribution of Relative Edge Density of the Graphs Based on
a Random Digraph Family}
\author{
Elvan Ceyhan\thanks{Address:
Department of Mathematics, Ko\c{c} University, 34450 Sar{\i}yer, Istanbul, Turkey.
e-mail: elceyhan@ku.edu.tr, tel:+90 (212) 338-1845, fax: +90 (212) 338-1559.
}
}
\date{\today}
\maketitle

{\bf short title:}
Relative Edge Density of Graphs Based on a Digraph Family

\begin{abstract}
\noindent
The vertex-random graphs called proximity catch digraphs (PCDs) have been introduced recently
and have applications in pattern recognition and spatial pattern analysis.
A PCD is a random directed graph (i.e., digraph) which is constructed from data using the
relative positions of the points from various classes.
Different PCDs result from different definitions of the
proximity region associated with each data point.
We consider the underlying and reflexivity graphs based on a family of PCDs
which is determined by a family of parameterized proximity maps called proportional-edge (PE) proximity map.
The graph invariant we investigate is the relative edge density of the underlying and reflexivity graphs.
We demonstrate that, properly scaled,
relative edge density of these graphs is a $U$-statistic,
and hence obtain the asymptotic normality of the relative edge density
for data from any distribution that satisfies mild regulatory conditions.
By detailed probabilistic and geometric calculations,
we compute the explicit form of the asymptotic normal distribution
for uniform data on a bounded region in the usual Euclidean plane.
We also compare the relative edge densities of the two types
of the graphs and the relative arc density of the PE-PCDs.
The approach presented here is also
valid for data in higher dimensions.
\end{abstract}



\vspace{.25 in}

\indent
$^\star$
This research was supported by the
European Commission under the Marie Curie International Outgoing Fellowship Programme
via Project \# 329370 titled PRinHDD.\\




\section{Introduction}
\label{sec:intro}
Classification and clustering have received considerable attention
in the probabilistic and statistical literature.
In this article, the probabilistic properties of a graph invariant
of a family of random graphs is investigated.
Vertex-random digraphs are directed graphs in which each vertex corresponds to a data point,
and directed edges (i.e., arcs) are defined in terms of a bivariate relation on the data points.
For example, nearest neighbor digraphs are defined by placing
an arc between each vertex and its nearest neighbor.
\cite{priebe:2001} introduced the class cover catch digraphs (CCCDs) in $\R$ and gave the exact and
the asymptotic distribution of the domination number of the CCCDs for uniform data on bounded intervals.
\cite{devinney:2002a}, \cite{marchette:2003}, \cite{priebe:2003b},
\cite{priebe:2003a}, and \cite{devinney:2006} applied the concept in higher dimensions and
demonstrated relatively good performance of CCCDs in classification.
Their methods involve \emph{data reduction} (i.e., \emph{condensing}) by
using approximate minimum dominating sets as \emph{prototype sets}
(since finding the exact minimum dominating set is an NP-hard
problem in general and for CCCD in multiple dimensions (see \cite{devinney:2006}).
Furthermore,
the exact and the asymptotic
distribution of the domination number of the CCCDs are not
analytically tractable in multiple dimensions.
For the domination number of CCCDs for one-dimensional data, a SLLN result is proved in \cite{devinney:2002b},
and this result is extended by \cite{wiermanSLLN:2008}; furthermore,
a generalized SLLN result is provided by \cite{wiermanSLLN:2008},
and a CLT is also proved by \cite{xiangCLT:2009}.
The asymptotic distribution of the domination number of CCCDs for non-uniform data
in $\mathbb{R}$ is also calculated in a rather general setting (\cite{ceyhan:dom-num-CCCD-NonUnif}).

\cite{ceyhan:Phd-thesis} generalized CCCDs to what is called \emph{proximity catch digraphs} (PCDs).
Let $(\Omega,\mathcal M)$ be a measurable space and $\X_n=\{X_1,X_2,\ldots,X_n\}$ and
$\Y_m=\{Y_1,Y_2,\ldots,Y_m\}$ be two sets of $\Omega$-valued random variables
from classes $\X$ and $\Y$, respectively, with joint probability distribution $F_{X,Y}$.
A PCD is comprised of a set of vertices and a set of arcs.
For example, in the two
class case, with classes $\X$ and $\Y$, the $\X$ points are the
vertices and there is an arc from $x_1 \in \X_n$ to $x_2\in \X_n$,
based on a binary relation which measures the relative
allocation of $x_1$ and $x_2$ with respect to $\Y$ points.
The PCDs are closely related to the class cover problem of \cite{cannon:2000}.
The \emph{class cover problem} for a target class, say $\X$, refers to finding a collection of neighborhoods,
$N_i$ around $X_i$ such that
(i) $\X_n \subseteq  \bigl(\bigcup_i N_i \bigr)$ and (ii) $\Y_m \cap \bigl(\bigcup_i N_i \bigr)=\emptyset$.
A collection of neighborhoods satisfying both conditions is called a {\em class cover}.
A cover satisfying (i) is a {\em proper cover} of class $\X$
while a cover satisfying (ii) is a {\em pure cover} relative to class $\Y$.
See \cite{priebe:2001} and \cite{cannon:2000} for more detail on the class cover problem.
The first PCD family is introduced by \cite{ceyhan:CS-JSM-2003};
the parameterized version of this PCD is developed by \cite{ceyhan:arc-density-CS}
where the relative arc density of the PCD is calculated and used for spatial pattern analysis.
\cite{ceyhan:dom-num-NPE-SPL} introduced another digraph family called \emph{proportional edge PCDs} (PE-PCDs)
and calculated the asymptotic distribution of its domination number and
used it for the same purpose (\cite{ceyhan:dom-num-NPE-MASA}; \cite{ceyhan:dom-num-NPE-Spat2011}).
The relative arc density of this PCD family is also computed
and used in spatial pattern analysis (\cite{ceyhan:arc-density-PE}).

The graphs based on digraphs are obtained by replacing arcs in the digraph
by edges based on bivariate relations.
If symmetric arcs are replaced by edges, then we obtain the reflexivity graph;
and if all arcs are replaced by edges disallowing multi-edges,
then we obtain the underlying graph (\cite{chartrand:2010}).
Properly scaled, we demonstrate that the relative
edge density of the underlying and reflexivity graphs of PE-PCDs is a $U$-statistic,
which has asymptotic normality by the general central limit theory of $U$-statistics.
Furthermore, we derive the explicit form of the parameters of the asymptotic normal distribution
of the relative edge density of the PCDs based on uniform data
in a bounded region in the Euclidean plane.

For the digraphs introduced by \cite{priebe:2001} (i.e., CCCDs),
whose relative arc density is also of the $U$-statistic form,
the asymptotic mean and variance of the relative density is not analytically tractable,
due to geometric difficulties encountered.
However, for the PCDs introduced in \cite{ceyhan:arc-density-PE}, and \cite{ceyhan:arc-density-CS},
the relative arc density has tractable asymptotic mean and variance.
The same holds for the underlying graphs as well.
We define the relative densities of graphs and digraphs and
derive their asymptotic distribution in general in Section \ref{sec:rel-dens-graph-digraph},
define the underlying and reflexivity graphs of PE-PCDs and
their relative edge densities in Section \ref{sec:underlying},
provide the asymptotic distribution of the relative edge density for uniform data
in Section \ref{sec:null-dist-edge-density}.
We treat the multiple triangle case in Section \ref{sec:multiple-triangle-case},
provide the discussion and conclusions in Section \ref{sec:discussion},
and the tedious calculations and long proofs are deferred to the Appendix.

\section{Relative Density of Graphs and Digraphs}
\label{sec:rel-dens-graph-digraph}
The main difference between a graph and a
digraph is that edges are directed in digraphs,
hence are called \emph{arcs}.
So the arcs are denoted as ordered pairs
while edges are not.

\subsection{Relative Edge Density of Graphs}
\label{sec:rel-dens-graph}
Let $G_n=(\V,\mE)$ be a graph with vertex set $\V=\{v_1,v_2,\ldots,v_n\}$
and edge set $\mE$.
The relative edge density of the graph $G_n$ which is of order $|\V| = n$,
denoted $\rho_e(G_n)$, is defined as
$$\rho_e(G_n) = \frac{2\,|\mE|}{n(n-1)}$$
where $|\cdot|$ denotes the set cardinality function (\cite{janson:2000}).
Thus $\rho_e(G_n)$ represents the ratio of the number of edges
in the graph $G_n$ to the number of edges in the complete graph of order $n$,
which is $n(n-1)/2$.
If $G_n$ is a random graph in which edges
result from a random process,
the \emph{edge probability} between vertices $v_i$ and $v_j$ is defined as
$p_e(i,j):=P(v_iv_j \in \mE)$ for all $i \not=j$, $i,j =1,2,\ldots,n$.

\begin{theorem}
\label{thm:edge-dens-result}
\textbf{(Main Result 1)}
Let $G_n=(\V,\mE)$ be a graph of order $n$
with $\V=\{v_1,v_2,\ldots,v_n\}$
and let $h^e_{ij}:=\I(v_iv_j \in \mE)$.
\begin{itemize}
\item[(a)]
If the set $\mE$ of edges result from a random process,
then $\rho_e(G_n)$ is a one-sample $U$-statistic of degree 2.
Moreover, if $p_e(i,j)=p_e$ for all $i \not=j$, $i,j =1,2,\ldots,n$
(i.e., the edge probability is constant for each pair of vertices $v_i,v_j$),
then $\rho_e(G_n)$ is an unbiased estimator of $p_e$.

\item[(b)]
If the set $\mE$ of edges result from a random process,
such that
$h^e_{ij}$ are identically distributed with $p_e(i,j)=p_e$ for all $i \not=j$, $i,j =1,2,\ldots,n$,
and $h^e_{ij}$ and $h^e_{kl}$ are independent for distinct $i,j,k,l$,
and $\nu_e:=\Cov(h^e_{ij},h^e_{ik}) >0 $ for all $i \not=j\not=k$, $i,j,k =1,2,\ldots,n$,
then $\sqrt{n}\,\bigl[\rho_e(G_n)-p_e]\stackrel{\mathcal L}{\longrightarrow} \N(0,4\,\nu_e)$
as $n \rightarrow \infty$,
where $\stackrel{\mathcal L}{\longrightarrow}$ stands for convergence in law or distribution
and $\N(\mu,\sigma^2)$ stands for the normal distribution with mean $\mu$
and variance $\sigma^2$.
\end{itemize}
\end{theorem}

\noindent \textbf{Proof:}
(a)
Assume the edges $\mE$ result from a random process
and let $G_n$ be the corresponding graph.
Let $h^e_{ij}=\I(v_iv_j \in \mE)$.
Since the edge $v_iv_j \in \mE$ can equivalently be expressed as
$v_jv_i \in \mE$ for all $i,j$,
we have $h^e_{ij}=h^e_{ji}$ and so $h^e_{ij}$ is symmetric in $i,j$.
Additionally, $|\mE|=\sum_{i < j} \,h^e_{ij}$.
So
$$\rho_e(G_n)=\frac{1}{({n \atop 2})}\sum_{i < j}\,h^e_{ij}.$$
Thus, $\rho_e(G_n)$ is a one-sample $U$-statistic of degree 2 with symmetric kernel $h^e_{ij}$
(\cite{lehmann:2004}).
Assume, moreover,
$P(v_iv_j \in \mE)=p_e$ for all $i \not=j$, $i,j =1,2,\ldots,n$.
Then for $i \not=j$,
we have $\E[h^e_{ij}]=\E[h^e_{12}]=\E[\I(v_1v_2 \in \mE)]=P(v_1v_2 \in \mE)=p_e$.
Hence $p_e$ is an estimable parameter of degree 2.
Furthermore,
$$\E[\rho_e(G_n)]=
\frac{2}{n(n-1)}\,\E[|\mE|]=
\sum_{i < j}\,\E[h^e_{ij}]=
\sum_{i < j}\,p_e=
\frac{2}{n(n-1)}\left(\frac{n(n-1)}{2} p_e\right)=p_e.$$
Then, $\rho_e(G_n)$ is actually an unbiased estimator of $p_e$.

(b) Assume the conditions for $G_n=(\V,\mE)$ stated in the hypothesis.
In part (a) we have shown that $p_e$ is an estimable parameter of degree 2,
and $\rho_e(G_n)$ is a one-sample $U$-statistic of degree 2 with symmetric kernel $h^e_{ij}$.
Furthermore, $(h^e_{ij})^2=h^e_{ij}$, since $(\I(v_iv_j \in \mE))^2=\I(v_iv_j \in \mE)$.
So $\E[(h^e_{ij})^2]= \E[h^e_{ij}]=p_e<\infty$ and
$\Cov(h^e_{ij},h^e_{ik})=\E[h^e_{ij}h^e_{ik}]-p_e^2 < \infty$,
since $\E[h^e_{ij}h^e_{ik}]=P((h^e_{ij},h^e_{ik})=(1,1))$.
By the hypothesis,
$\nu_e=\Cov(h^e_{ij},h^e_{ik})>0$.
Then by Theorem 3.3.13 in \cite{randles:1979},
we have $\sqrt{n}\,\bigl[\rho_e(G_n)-p_e]\stackrel{\mathcal L}{\longrightarrow} \N(0,4\,\nu_e)$
as $n \rightarrow \infty$. $\blacksquare$

In part (b) of Theorem \ref{thm:edge-dens-result},
we have
$\Cov(h^e_{ij},h^e_{ik})=\E[h^e_{ij}h^e_{ik}]-\E[h^e_{ij}]\E[h^e_{ik}]=\E[h^e_{ij}h^e_{ik}]-p_e^2$
where $h^e_{ij}h^e_{ik}=\I(v_iv_j \in \mE)\I(v_iv_k \in \mE)=\I(\{v_iv_j,v_iv_k\} \subset \mE)$,
so $\E[h^e_{ij}h^e_{ik}]=P(\{v_iv_j,v_iv_k\} \subset \mE)$.
Hence $\nu_e>0$ iff $P(\{v_iv_j,v_iv_k\} \subset \mE)>p_e^2$.
Notice that
$\E[|h^e_{ij}|^3]=\E[h^e_{ij}]=p_e<\infty$
and assuming  $\nu_e >0$,
then the sharpest rate of convergence in the asymptotic normality of
$\rho_e(G_n)$ is (\cite{callaert:1978})
as follows:
$$\sup_{t\in \R} \left| P \left( \frac{\sqrt{n}(\rho_e(G_n)-p_e)}{2\,\sqrt{\nu_e}}\le t \right)-\Phi(t) \right|
\le C\,\cdot p_e\, \cdot ( 4\,\nu_e)^{-3/2}\cdot n^{-1/2}=
C_e\frac{p_e}{\sqrt{n\,\nu_e^3}}$$
where $C_e$ is a constant and $\Phi(\cdot)$ is the standard normal
distribution function.
Furthermore, we have
$$\Var[h^e_{ij}]=
\E[(h^e_{ij})^2]-(\E[h^e_{ij}])^2=
\E[h^e_{ij}]-p_e^2=
p_e-p_e^2=p_e(1-p_e).$$

The graph $G_n$ in Theorem \ref{thm:edge-dens-result} is not a deterministic graph,
but a random one.
In general a random graph
is obtained by starting with a set of $n$ vertices and adding edges between them at random.
Most commonly studied is the Erd\H{o}s--–R\'{e}nyi model,
denoted $G(n,p)$, in which every possible edge occurs independently with probability $p$ (\cite{erdos:1959}).
Notice that the random graph $G(n,p)$ satisfies part (a) of Theorem \ref{thm:edge-dens-result},
so the relative edge density of $G(n,p)$ is a $U$-statistic;
however, the asymptotic distribution of its relative edge density is degenerate
(with $\rho(G(n,p))\stackrel{\mathcal L}{\longrightarrow} p$ as $n \rightarrow \infty$)
since the covariance term is zero due to the independence between the edges.

\subsection{Relative Arc Density of Digraphs}
\label{sec:rel-dens-digraph}
Let $D_n=(\V,\A)$ be a digraph with vertex set $\V=\{v_1,v_2,\ldots,v_n\}$
and arc set $\A$.
The relative arc density of the digraph $D_n$ which is of order $|\V| = n$,
denoted $\rho_a(D_n)$, is defined as
$$\rho_a(D_n) = \frac{|\A|}{n(n-1)}.$$
Thus $\rho_a(D_n)$ represents the ratio of the number of arcs
in the digraph $D_n$ to the number of arcs in the complete digraph of order $n$,
which is $n(n-1)$.
If $D_n=(\V,\A)$ is a random digraph in which arcs 
result from a random process,
the \emph{arc probability} between vertices $v_i,v_j$ is defined as
$p_a(i,j):=P((v_i,v_j) \in \A)$ for all $i \not=j$, $i,j =1,2,\ldots,n$.

\begin{theorem}
\label{thm:arc-dens-result}
\textbf{(Main Result 2)}
Let $D_n=(\V,\A)$ be a digraph of order $n$
with $\V=\{v_1,v_2,\ldots,v_n\}$
and let $g_{ij}:=\I((v_i,v_j) \in \A)$.
\begin{itemize}
\item[(a)]
If the set $\A$ of arcs result from a random process,
then $\rho_a(D_n)$ is a one-sample $U$-statistic of degree 2.
Moreover, if $p_a(i,j)=p_a$ for all $i \not=j$, $i,j =1,2,\ldots,n$,
(i.e., the arc probability is constant for each pair of vertices $v_i,v_j$),
then $\rho_a(D_n)$ is an unbiased estimator of $p_a$.

\item[(b)]
If the set $\A$ of arcs result from a random process
such that
$g_{ij}$ are identically distributed with $p_a(i,j)=p_a$ for all $i \not=j$, $i,j =1,2,\ldots,n$,
$g_{ij}$ and $g_{kl}$ are independent for distinct $i,j,k,l$,
and $\Cov(g_{ij},g_{kl})>0$ for all $i \not=j$ and $k\not=l$
and exactly one of $i,j$ is equal to exactly one of $k,l$ for $i,j,k,l =1,2,\ldots,n$,
then $\sqrt{n}\,\bigl[\rho_a(D_n)-p_a]\stackrel{\mathcal L}{\longrightarrow} \N(0,\nu_a)$
as $n \rightarrow \infty$,
where $\nu_a=\lim_{n \rightarrow \infty}n\,\Var[\rho_a(D_n)]$.
\end{itemize}
\end{theorem}

\noindent \textbf{Proof:}
(a)
Assume that the arcs $\A$ result from a random process
and let $D_n$ be the corresponding digraph.
Let $g_{ij}=\I((v_i,v_j) \in \A)$.
The arcs $(v_i,v_j) \in \A$ and $(v_j,v_i) \in \A$ are distinct for $i \not=j$,
so $g_{ij}$ is not symmetric in $i,j$.
But we can define a symmetric kernel as
$h^a_{ij}=(g_{ij}+g_{ji})/2$.
Then we have, $|\A|=\sum_{i < j} \,h^a_{ij}$.
So
$$\rho_a(D_n)=\frac{1}{({n \atop 2})}\sum_{i < j}\,h^a_{ij}.$$
Thus, $\rho_a(D_n)$ is a one-sample $U$-statistic of degree 2 with symmetric kernel $h^a_{ij}$.
Assume, moreover,
$P((v_i,v_j) \in \A)=p_a$ for all $i \not=j$, $i,j =1,2,\ldots,n$.
Then for $i \not=j$,
we have
$$\E[h^a_{ij}]=\E[(g_{ij}+g_{ji})/2]=(\E[g_{ij}]+\E[g_{ji}])/2=\E[g_{ij}]=\E[g_{12}]=P((v_1,v_2) \in \A)=p_a.$$
Hence $p_a$ is an estimable parameter of degree 2.
Furthermore,
$$\E[\rho_a(D_n)]=
\frac{2}{n(n-1)}\,\E[|\A|]=
\frac{2}{n(n-1)}\,\sum_{i < j}\,\E[h^a_{ij}]=
\frac{2}{n(n-1)}\,\sum_{i < j}\,\E[g_{ij}]=
\frac{2}{n(n-1)}\,\sum_{i < j}\,p_a=
p_a.$$
Then, $\rho_a(D_n)$ is actually an unbiased estimator of $p_a$.

(b) Assume the conditions for $D_n=(\V,\A)$ stated in the hypothesis.
In part (a) we have shown that $p_a$ is an estimable parameter of degree 2,
and $\rho_a(D_n)$ is a one-sample $U$-statistic of degree 2 with symmetric kernel $h^a_{ij}$.
Furthermore, $(h^a_{ij})^2=(g^2_{ij}+2\,g_{ij}g_{ji}+g^2_{ji})/4=(g_{ij}+2\,g_{ij}g_{ji}+g_{ji})/4$,
since $(\I((v_i,v_j) \in \A))^2=\I((v_i,v_j) \in \A)$.
So $\E[(h^a_{ij})^2]= (\E[g_{ij}]+2\,\E[g_{ij}g_{ji}]+\E[g_{ji}])/4=(2\,p_a+2\,\E[g_{ij}g_{ji}])/4$.
Since $\E[g_{ij}g_{ji}]=\E[\I((v_i,v_j) \in \A)\I((v_j,v_i) \in \A)]=
\E[\I(\{(v_i,v_j),(v_j,v_i)\} \subset \A)]=P(\{(v_i,v_j),(v_j,v_i)\} \subset \A)$,
we have $\E[(h^a_{ij})^2]<\infty$.
By the hypothesis,
$\nu_a=\Cov(h^a_{ij},h^a_{ik})=
\Cov((g_{ij}+g_{ji})/2,(g_{ik}+g_{ki})/2)=
(\Cov[g_{ij},g_{ik}]+\Cov[g_{ij},g_{ki}]+\Cov[g_{ji},g_{ik}]+\Cov[g_{ji},g_{ki}])/4>0$.
Additionally,
$\Cov(g_{ij},g_{kl})=\E[g_{ij}g_{kl}]-p_a^2 < \infty$,
since $\E[g_{ij}g_{kl}]=P((g_{ij},g_{kl})=(1,1))$.
Hence $\nu_a< \infty$ as well.
Then we have $\sqrt{n}\,\bigl[\rho_a(D_n)-p_a]\stackrel{\mathcal L}{\longrightarrow} \N(0,\nu_a)$
as $n \rightarrow \infty$.
$\blacksquare$

In part (b) of Theorem \ref{thm:arc-dens-result},
we have
$$\nu_a=\Cov(h^a_{ij},h^a_{ik})=\E[h^a_{ij}h^a_{ik}]-\E[h^a_{ij}]\E[h^a_{ik}]=\E[h^a_{ij}h^a_{ik}]-p_a^2$$
where
\begin{multline*}
h^a_{ij}h^a_{ik}=(g_{ij}+g_{ji})(g_{ik}+g_{ki})/4=(g_{ij}g_{ik}+g_{ij}g_{ki}+g_{ji}g_{ik}+g_{ji}g_{ki})/4=
(\I((v_i,v_j) \in \A)\I((v_i,v_k) \in \A)+\\
\I((v_i,v_j) \in \A)\I((v_k,v_i) \in \A)+
\I((v_j,v_i) \in \A)\I((v_i,v_k) \in \A)+
\I((v_j,v_i) \in \A)\I((v_k,v_i) \in \A))/4=\\
(\I(\{(v_i,v_j),(v_i,v_k)\} \subset \A)+
\I(\{(v_i,v_j),(v_k,v_i)\} \subset \A)+
\I(\{(v_j,v_i),(v_i,v_k)\} \subset \A)+
\I(\{(v_j,v_i),(v_k,v_i)\} \subset \A))/4.
\end{multline*}
So, $\nu_a>0$ iff
{\small
$$(P(\{(v_i,v_j),(v_i,v_k)\} \subset \A)+P(\{(v_i,v_j),(v_k,v_i)\} \subset \A)+
P(\{(v_j,v_i),(v_i,v_k)\} \subset \A)+P(\{(v_j,v_i),(v_k,v_i)\} \subset \A))/4 > p_a^2.$$
}
Notice that
\begin{multline*}
\E[|h^a_{ij}|^3]=
\E[(g_{ij}+g_{ji})^3/8]=
\E[g^3_{ij}+3\,g^2_{ij}g_{ji}+3\,g_{ij}g^2_{ji}+g^3_{ji}]/8=
\E[g_{ij}+3\,g_{ij}g_{ji}+3\,g_{ij}g_{ji}+g_{ji}]/8=\\
(2\,\E[g_{ij}]+6\,\E[g_{ij}g_{ji}])/8=
(p_a+3\,p_{sa})/4 < \infty,
\end{multline*}
where $p_{sa}=P(g_{ij}g_{ji}=1)=P(\{(v_i,v_j),(v_j,v_i)\} \subset \A)$
is the symmetric arc probability in $D_n$.
Assuming  $\nu_a >0$,
then the sharpest rate of convergence in the asymptotic normality of
$\rho_a(D_n)$ is (\cite{callaert:1978})
as follows:
$$\sup_{t\in \R} \left| P \left( \frac{\sqrt{n}(\rho_a(D_n)-p_a)}{\sqrt{\nu_a}}\le t \right)-\Phi(t) \right|
\le C\,\cdot p_a\, \cdot (\nu_a)^{-3/2}\cdot n^{-1/2}=
C_a\,\frac{p_a}{\sqrt{n\,\nu_a^3}}$$
where $C_a$ is a constant.
Furthermore, we have
\begin{multline*}
\Var[h^a_{ij}]=
\E[(h^a_{ij})^2]-(\E[h^a_{ij}])^2=
\E[(g_{ij}+g_{ji})^2/4]-p_a^2=
\E[g^2_{ij}+2\,g_{ij}g_{ji}+g^2_{ji}]/4-p_a^2=\\
\E[g_{ij}+2\,g_{ij}g_{ji}+g_{ji}]/4-p_a^2=
(\E[g_{ij}]+2\,\E[g_{ij}g_{ji}]+\E[g_{ji}])/4-p_a^2=
(p_a+2\,p_{sa}+p_a)/4-p_a^2=\\
(p_a+p_{sa})/2-p_a^2.
\end{multline*}

The digraph $D_n$ in Theorem \ref{thm:arc-dens-result} is not a deterministic digraph,
but a random one.
In general a random digraph, just like a random graph,
can be obtained by starting with a set of $n$ vertices and adding arcs between them at random.
We can consider the counterpart of the Erd\H{o}s–--R\'{e}nyi model for digraphs,
denoted $D(n,p)$, in which every possible arc occurs independently with probability $p$.
Notice that the random digraph $D(n,p)$ satisfies part (a) of Theorem \ref{thm:arc-dens-result},
so the relative arc density of $D(n,p)$ is a $U$-statistic,
however, the asymptotic distribution of its relative arc density is degenerate
(with $\rho(D(n,p))\stackrel{\mathcal L}{\longrightarrow} p$ as $n \rightarrow \infty$)
since the covariance term is zero due to the independence between the arcs.

\section{Relative Edge Density of the Graphs Based on PCDs}
\label{sec:underlying}

\subsection{Proximity Catch Digraphs and the Corresponding Graphs}
\label{sec:under-PCDs}
Let $(\Omega,\mathcal{M})$ be a measurable space and
$d(\cdot,\cdot):\Omega\times \Omega \rightarrow [0,\infty)$ be any distance function.
Consider $N:\Omega \rightarrow \mathcal P(\Omega)$,
where $\mathcal P(\cdot)$ represents the power set functional.
Then given $\Y_m \subset \Omega$,
the {\em proximity map}
$N(\cdot)$ associates with each point $x \in \Omega$
a {\em proximity region} $N(x) \subseteq \Omega$.
The region $N(x)$ is defined in terms of the distance
between $x$ and $\Y_m$.
we define the vertex-random PCD, $D_n$, with vertex set
$\V=\{X_1,X_2,\ldots,X_n\}$ and arc set $\A$ by
$(X_i,X_j) \in \A \iff X_j \in N(X_i)$ where point $X_i$ ``catches" point $X_j$.
The random digraph $D_n$ depends on the (joint) distribution of the $X_i$
and on the map $N(\cdot)$.
The adjective {\em proximity} --- for the catch digraph $D_n$ and for the map $N(\cdot)$ ---
comes from thinking of the region $N(x)$ as representing those
points in $\Omega$ ``close'' to $x$ (\cite{toussaint:1980} and \cite{jaromczyk:1992}).
The $\G_1$-region $\G_1(\cdot,N):\Omega \rightarrow \mathcal P(\Omega)$
associates the region $\G_1(x,N):=\{z \in \Omega: x \in N(z)\}$ with each point $x \in \Omega$.
A $\G_1$-region is sort of a ``dual" of the corresponding proximity region
and is closely associated with domination number being equal to one.
If $X_1,X_2,\ldots,X_n$ are $\Omega$-valued random variables,
then the $N(X_i)$ (and $\G_1(X_i,N)$), $i=1,2,\ldots,n$ are random sets.
If the $X_i$ are independent and identically distributed,
then so are the random sets $N(X_i)$ (and $\G_1(X_i,N)$).

%
If $X_1,X_2,\ldots,X_n \stackrel{iid}{\sim} F$,
then, by Theorem \ref{thm:arc-dens-result},
the relative arc density
of the associated vertex-random proximity catch digraph, $D_n$,
denoted $\rho_a(D_n)$, is a $U$-statistic.
See \cite{ceyhan:arc-density-CS, ceyhan:arc-density-PE}
for its derivation and other details.

The \emph{reflexivity graph} for digraph $D_n=(\V, \A)$ is the graph
$G_{\la}(D_n)=(\V, \mE_{\la})$ where $\mE_{\la}$ is the set of edges
such that  $uv \in \mE_{\la}$ iff $(u,v) \in \A$ \emph{\textbf{and}} $(v,u) \in \A$.
The \emph{underlying graph} of a
digraph is the graph obtained by replacing each arc $(u,v) \in \A$ or
each symmetric arc, $\{(u,v),\,(v,u)\} \subset \A$ by the edge $uv$.
Then, the underlying graph for $D_n=(\V, \A)$ is the graph
$G_{\lo}(D_n)=(\V, \mE_{\lo})$ where $\mE_{\lo}$ is the set of edges such that
$uv \in \mE_{\lo}$ iff $(u,v) \in \A$ \emph{\textbf{or}} $(v,u) \in \A$.


Consider the vertex-random PCD, $D_n$,
with vertex set $\V=\{X_1,X_2,\ldots,X_n\}$
and arc set $\A$ defined by
$(X_i,X_j) \in \A \iff X_j \in N(X_i)$.
The {\em reflexivity graph}, $G_{\la}$, of $D_n$ with the vertex set $\V$ and
the edge set $\mE_{\la}$ is defined by $X_iX_j \in \mE_{\la}$ iff
$(X_i,X_j) \in \A \text{ and } (X_j,X_i) \in \A$.
Likewise, the {\em underlying graph}, $G_{\lo}$, of $D_n$ with the vertex set $\V$
and the edge set $\mE_{\lo}$ is defined by
$X_iX_j \in \mE_{\lo} \iff (X_i,X_j) \in \A \text{ or } (X_j,X_i) \in \A$.
Then
$X_iX_j \in \mE_{\la}$ iff $X_j \in N(X_i) \text{ and } X_i \in N(X_j)$
iff $X_j \in N(X_i) \text{ and } X_j \in \G_1(X_i,N)$ iff $X_j \in N(X_i) \cap  \G_1(X_i,N)$.
Similarly, $X_iX_j \in \mE_{\lo}$ iff $X_j \in N(X_i) \cup  \G_1(X_i,N)$.

\subsection{Relative Arc Density of the PCDs}
\label{sec:arc-density}
The relative arc density of the PCD, $D_n$, is denoted as $\rho_a(D_n)$.
Let $g_{ij}=\I(\left( X_i,X_j \right)\in \A)=\I(X_j \in N(X_i))$.
Then for $X_i \stackrel{iid}{\sim}F$, $i=1,2,\ldots,n$,
$\rho_a(D_n)$ can be written as
$$\rho_a(D_n)=\frac{2}{n\,(n-1)}\sum\hspace*{-0.1 in}\sum_{i < j \hspace*{0.25 in}}   \hspace*{-0.1 in} \,h_{ij}$$
where
$2\,h_{ij} = (g_{ij}+g_{ji})=\I(X_j \in N(X_i))+\I(X_i \in N(X_j))$
is the number of arcs between $X_i$ and $X_j$ in $D_n$.
Note that $h_{ij}$ is a symmetric kernel with finite variance since $0 \le h_{ij} \le 1$.
Moreover, $\rho_a(D_n)$ is a random
variable that depends on $n$, $F$, and $N(\cdot)$ (i.e., $\mathcal Y$).
But $\E\left[ \rho_a(D_n) \right]$ only depends on $F$ and $N(\cdot)$.
That is,
\begin{equation}
\label{eqn:E[rho-arc-D]}
0\le \E\left[ \rho_a(D_n) \right]=\frac{2}{n\,(n-1)}\sum\hspace*{-0.1 in}\sum_{i < j \hspace*{0.25 in}}
\hspace*{-0.1 in} \,\E[h_{ij}]=\E\left[ h_{12} \right]
\end{equation}
where $2\,\E\left[ h_{12} \right]=
\E[\I((X_1,X_2) \in \A) +\I((X_2,X_1) \in \A)]=
\E[\I((X_1,X_2) \in \A)]+ \E[\I((X_2,X_1) \in \A)]=
P((X_1,X_2) \in \A)+ P((X_2,X_1) \in \A)=
P(X_2 \in N(X_1))+P(X_1 \in N(X_2))=2\,p_a(N)$.
Hence $\E\left[ h_{12} \right]=p_a(N)$,
which is the arc probability for the PCD, $D_n$.
Notice also that $p_a(N)=P(X_j \in N(X_i))$ for $i \not= j$.
Furthermore,
\begin{equation}
\label{eqn:Var[rho-and-D]}
0 \le \Var\left[ \rho_a(D_n) \right]=
\frac{4}{n^2\,(n-1)^2}\Var\left[\sum\hspace*{-0.1 in}\sum_{i < j \hspace*{0.25 in}}   \hspace*{-0.1 in} h_{ij}\right].
\end{equation}
Expanding this expression,
we have
$$\Var\left[ \rho_a(D_n) \right]=\frac{2}{n\,(n-1)}\Var\left[ h_{12} \right]+
\frac{4(n-2)}{n\,(n-1)} \, \Cov\left[ h_{12},h_{13} \right].$$
As in Section \ref{sec:rel-dens-digraph}, we have
$$\Var\left[ h_{12} \right]=(p_a(N)-\left[p_a(N)\right]^2)=
p_a(N)\left(1-p_a(N)\right).$$
Moreover, the covariance is as follows
$$\Cov\left[ h_{12},h_{13} \right]=
\E\left[ h_{12}.h_{13} \right]-\E\left[ h_{12} \right]\E\left[ h_{13} \right],$$
where $\E\left[ h_{12} \right]=\E\left[ h_{13} \right]=p_a(N)$ and,
\begin{eqnarray*}
4\,\E\left[ h_{12}.h_{13} \right]& = & \E[(g_{12}+g_{21})(g_{13}+g_{31})]=\E[g_{12}g_{13}+g_{12}g_{31}+g_{21}g_{13}+g_{21}g_{31}]\\
& = & \E[\I(X_2 \in N(X_1)\I( X_3 \in N(X_1))+\I(X_2 \in N(X_1)\I( X_1 \in N(X_3))+\\
&  & \I(X_1 \in N(X_2)\I( X_3 \in N(X_1))]+\I(X_1 \in N(X_2)\I( X_1 \in N(X_3))]\\
& = & \E[\I(\{X_2,X_3\} \subset N(X_1))+\I(X_2 \in N(X_1)\I(X_3 \in \G_1(X_3,N))+\\
& & \I(X_2 \in \G_1(X_1)\I( X_3 \in N(X_1))]+\I(X_2 \in \G_1(X_1,N)\I( X_3 \in \G_1(X_1,N))]\\
& = & P(\{X_2,X_3\} \subset N(X_1))+2\,P(X_2 \in N(X_1),\I(X_3 \in \G_1(X_3,N))+P(\{X_2,X_3\} \subset \G_1(X_1,N)).
\end{eqnarray*}

The digraph $D_n$ is a random digraph
where the arc probability is $P((X_i,X_j) \in \A)=p_a(N)$ for $i\not=j$
and is an estimable parameter of degree 2.
Using Equation \eqref{eqn:E[rho-arc-D]},
we have that $\rho_a(D_n)$ is an unbiased estimator of $p_a(N)$.
Notice that for PCDs,
the set of vertices $\V=\X_n$ is a random sample from a distribution $F$
(i.e., the vertices directly result from a random process),
and the arcs are defined based on the random sets (i.e., proximity regions) $N(X_i)$
as described before.
Hence the set of arcs $\A$ (indirectly) result from a random process
such that
$g_{ij}$ are identically distributed
and $g_{ij}$ and $g_{kl}$ are independent for distinct $i,j,k,l$.
Furthermore, we have $\nu_a(N):=\Cov\left[ h_{ij},h_{ik} \right]<\infty$ as before.
Then we have the following corollary to the Main Result 2.
\begin{corollary}
\label{cor:arc-dens-PCD}
The relative arc density, $\rho_a(D_n)$, of the PCD,
$D_n$, is a one-sample $U$-statistic of degree 2
and is an unbiased estimator of $p_a(N)$.
If, additionally, $\nu_a(N)=\Cov\left[ h_{ij},h_{ik} \right]>0$ for all $i \not=j\not=k$, $i,j,k =1,2,\ldots,n$,
then $\sqrt{n}\,\bigl[\rho_a(D_n)-\mu(N)]\stackrel{\mathcal L}{\longrightarrow} \N(0,4\,\nu_a(N))$
as $n \rightarrow \infty$.
\end{corollary}

In the above corollary,
$\nu_a(N) > 0$ iff
$P(\{X_2,X_3\} \subset N(X_1))+2\,P(X_2 \in N(X_1),\I(X_3 \in \G_1(X_3,N))+P(\{X_2,X_3\} \subset \G_1(X_1,N))>4\,p^2_a(N)$.

\subsubsection{The Joint Distribution of $\left( h_{12},h_{13} \right)$}
The pair $\left( h_{12},h_{13} \right)$ is a bivariate discrete random
variable with nine possible values:
$$\left( h_{12},h_{13} \right)\in \{(0,0),(0,1),(0,2),(1,0),(1,1),(1,2),(2,0),(2,1),(2,2)\}.$$
Then finding the joint distribution of $\left( h_{12},h_{13} \right)$
is equivalent to finding the joint probability mass function of
$\left( h_{12},h_{13} \right)$.

First, note that
$$\left( h_{12},h_{13} \right)=(0,0)\text{ iff } g_{12}=g_{21}=g_{13}=g_{31}=0 \text{ iff }$$
$$\I(X_2 \in N(X_1))=\I(X_1 \in N(X_2))=\I(X_3 \in N(X_1))=\I(X_1 \in N(X_3))=0 \text{ iff }$$
$$\I(X_2 \in N(X_1))=\I(X_2 \in \G_1(X_1,N))=\I(X_3 \in N(X_1))=\I(X_3 \in \G_1(X_1,N))=0 \text{ iff }$$
$$\I(X_2 \in \TY\setminus N(X_1))=\I(X_2 \in \TY\setminus \G_1(X_1,N))=\I(X_3 \in \TY\setminus N(X_1))=\I(X_3 \in \TY\setminus \G_1(X_1,N))=1 \text{ iff }$$
$$\I(\{X_2,X_3 \} \subset \TY\setminus N(X_1))=\I(\{X_2,X_3 \} \subset \TY\setminus\G_1(X_1,N))=1 \text{ iff }$$
$$\I(\{X_2,X_3 \} \subset \TY\setminus (N(X_1)\cup \G_1(X_1,N)))=1.$$
Hence $P(\left( h_{12},h_{13} \right)=(0,0))=P(\{X_2,X_3 \} \subset \TY\setminus (N(X_1)\cup \G_1(X_1,N)))$.

Furthermore, by symmetry,
$P(\left( h_{12},h_{13} \right)=(0,1))=P(\left( h_{12},h_{13} \right)=(1,0))$,
$P(\left( h_{12},h_{13} \right)=(0,2))=P(\left( h_{12},h_{13} \right)=(2,0))$,
and
$P(\left( h_{12},h_{13} \right)=(1,2))=P(\left( h_{12},h_{13} \right)=(2,1))$.
So it suffices to calculate one of each pair of the probabilities in the above cases.

Finally,
$$\left( h_{12},h_{13} \right)=(2,2)\text{ iff } g_{12}=g_{21}=g_{13}=g_{31}=1 \text{ iff }$$
$$\I(X_2 \in N(X_1))=\I(X_1 \in N(X_2))=\I(X_3 \in N(X_1))=\I(X_1 \in N(X_3))=1 \text{ iff }$$
$$\I(X_2 \in N(X_1))=\I(X_2 \in \G_1(X_1,N))=\I(X_3 \in N(X_1))=\I(X_3 \in \G_1(X_1,N))=1 \text{ iff }$$
$$\I(\{X_2,X_3 \} \subset N(X_1))=\I(\{X_2,X_3 \} \subset \G_1(X_1,N))=1 \text{ iff }$$
$$\I(\{X_2,X_3 \} \subset (N(X_1)\cap \G_1(X_1,N)))=1.$$
Hence $P(\left( h_{12},h_{13} \right)=(2,2))=P(\{X_2,X_3 \} \subset (N(X_1)\cap \G_1(X_1,N)))$.
Finally,
$P(\left( h_{12},h_{13} \right)=(1,1))$ can be found by subtracting
the sum of the probabilities in the other cases from 1.

\subsection{Relative Edge Density of the Reflexivity Graphs Based on PCDs}
\label{sec:AND-edge-density}
The relative edge density of the reflexivity graph, $G_{\la}(D_n)$, based on the PCD, $D_n$,
is denoted as $\rho_{\la}(D_n)$.
For $X_i \stackrel{iid}{\sim}F$, $i=1,2,\ldots,n$,
one can write down the relative edge density as
$$\rho_{\la}(D_n)=\frac{2}{n\,(n-1)}\sum\hspace*{-0.1 in}\sum_{i < j \hspace*{0.25 in}}   \hspace*{-0.1 in} \,h^{\la}_{ij}$$
where
\begin{eqnarray*}
h^{\la}_{ij} & = & \I(\left( X_i,X_j \right)\in \mE_{\la})= \I((X_i,X_j) \in \A,\;(X_j,X_i) \in \A)\\
&=&\I(X_j \in N(X_i),\;X_i \in N(X_j))=\I(X_j \in N(X_i),\;X_j \in \G_1\left( X_i,N \right))\\
&=&\I(X_j \in N(X_i)\cap \G_1\left( X_i,N \right))
\end{eqnarray*}
is the number of edges between $X_i$ and $X_j$ in $G_{\la}(D_n)$
or
number of symmetric arcs between $X_i$ and $X_j$ in $D_n$.
Note that $h^{\la}_{ij}$ is a symmetric kernel with finite variance since $0 \le h^{\la}_{ij} \le 1$.
Moreover, $\rho_{\la}(D_n)$ is a random
variable that depends on $n$, $F$, and $N(\cdot)$ (i.e., $\mathcal Y$).
But $\E\left[ \rho_{\la}(D_n) \right]$ only depends on $F$ and $N(\cdot)$.
That is,
\begin{equation}
\label{eqn:E[rho-and-D]}
0\le \E\left[ \rho_{\la}(D_n) \right]=\frac{2}{n\,(n-1)}\sum\hspace*{-0.1 in}\sum_{i < j \hspace*{0.25 in}}
\hspace*{-0.1 in} \,\E[h^{\la}_{ij}]=\E\left[ h^{\la}_{12} \right]
\end{equation}
where $\E\left[ h^{\la}_{12} \right]=
\E[\I((X_1,X_2) \in \A\, \, , \,\, (X_2,X_1) \in \A)]=
P((X_1,X_2) \in \A\, \, , \,\, (X_2,X_1) \in \A)=
P(X_2 \in N(X_1)\cap \G_1(X_1,N))=p_{\la}(N)$.
Notice that $p_{\la}(N)$ is the edge probability for
the underlying graph $G_{\la}(D_n)$,
but it is symmetric arc probability $p_{sa}(D_n)$ for the PCD, $D_n$.
Notice also that $p_{\la}(N)=P(X_j \in N(X_i)\cap \G_1\left( X_i,N \right))$ for $i \not= j$.
Furthermore,
\begin{equation}
\label{eqn:Var[rho-and-D]}
0 \le \Var\left[ \rho_{\la}(D_n) \right]=
\frac{4}{n^2\,(n-1)^2}\Var\left[\sum\hspace*{-0.1 in}\sum_{i < j \hspace*{0.25 in}}   \hspace*{-0.1 in} h^{\la}_{ij}\right].
\end{equation}
Expanding this expression,
we have
$$\Var\left[ \rho_{\la}(D_n) \right]=\frac{2}{n\,(n-1)}\Var\left[ h^{\la}_{12} \right]+
\frac{4\,(n-2)}{n\,(n-1)} \, \Cov\left[ h^{\la}_{12},h^{\la}_{13} \right].$$
Here, as in Section \ref{sec:rel-dens-graph}, we have
$$\Var\left[ h^{\la}_{12} \right]=p_{\la}(N)-\left[p_{\la}(N)\right]^2=
p_{\la}(N)\left(1-p_{\la}(N)\right).$$
Moreover, the covariance is as follows
$$\Cov\left[ h^{\la}_{12},h^{\la}_{13} \right]=
\E\left[ h^{\la}_{12}.h^{\la}_{13} \right]-\E\left[ h^{\la}_{12} \right]\E\left[ h^{\la}_{13} \right].$$
Since $\E\left[ h^{\la}_{12} \right]=\E\left[ h^{\la}_{13} \right]=p_{\la}(N)$ and,
\begin{eqnarray*}
\E\left[ h^{\la}_{12}.h^{\la}_{13} \right]
& = & \E[\I(X_2 \in N(X_1)\cap \G_1(X_1,N))\, . \,\I( X_3 \in N(X_1)\cap \G_1(X_1,N))]\\
& = & \E[\I(X_2 \in N(X_1)\cap \G_1(X_1,N)\, , \, X_3 \in N(X_1)\cap \G_1(X_1,N))]\\
& = & P(X_2 \in N(X_1)\cap \G_1(X_1,N)\, , \, X_3 \in N(X_1)\cap \G_1(X_1,N))\\
& = & P(\{X_2,X_3\} \subset N(X_1)\cap \G_1(X_1,N)),
\end{eqnarray*}
it follows that
$$\Cov\left[ h^{\la}_{12},h^{\la}_{13} \right]=P(\{X_2,X_3\} \subset N(X_1)\cap \G_1(X_1,N))-\left[p_{\la}(N)\right]^2.$$

The underlying graph $G_{\la}(D_n)$ is a random graph
where the edge probability is $P(X_iX_j \in \mE_{\la})=p_{\la}(N)$ for $i\not=j$
and is an estimable parameter of degree 2.
Using Equation \eqref{eqn:E[rho-and-D]},
we have that $G_{\la}(D_n)$ is an unbiased estimator of $p_{\la}(N)$.
Notice that for the reflexivity graphs based on the PCDs,
the set of vertices $\V=\X_n$ is a random sample from a distribution $F$
(i.e., the vertices directly result from a random process),
and the edges are defined based on the random sets $N(X_i)\cap \G_1(X_1,N)$
as described before.
Hence the set of edges $\mE_{\la}$ (indirectly) result from a random process
such that
$h^{\la}_{ij}$ are identically distributed
and $h^{\la}_{ij}$ and $h^{\la}_{kl}$ are independent for distinct $i,j,k,l$.
Furthermore, we have $\Cov\left[ h^{\la}_{ij},h^{\la}_{kl} \right]<\infty$ as before.
Then we have the following corollary to the Main Result 1.

\begin{corollary}
\label{cor:edge-dens-and}
The relative edge density, $\rho_{\la}(D_n)$, of the reflexivity graph,
$G_{\la}(D_n)$, is a one-sample $U$-statistic of degree 2
and is an unbiased estimator of $p_{\la}(N)$.
If, additionally, $\nu_{\la}(N):=\Cov\left[ h^{\la}_{ij},h^{\la}_{ik} \right]>0$ for all $i \not=j\not=k$, $i,j,k =1,2,\ldots,n$,
then $\sqrt{n}\,\bigl[\rho_{\la}(D_n)-p_{\la}(N)]\stackrel{\mathcal L}{\longrightarrow} \N(0,4\,\nu_{\la}(N))$
as $n \rightarrow \infty$.
\end{corollary}
In the above corollary,
$\nu_{\la}(N) > 0$ iff
$P(\{X_2,X_3\} \subset N(X_1)\cap \G_1(X_1,N))>\left[p_{\la}(N)\right]^2$.

\subsubsection{The Joint Distribution of $\left( h^{\la}_{12},h^{\la}_{13} \right)$}
By definition $\left( h^{\la}_{12},h^{\la}_{13} \right)$ is a discrete random
variable with four possible values:
$$\left( h^{\la}_{12},h^{\la}_{13} \right)\in \{(0,0),(1,0),(0,1),(1,1)\}.$$
Then finding the joint distribution of $\left( h^{\la}_{12},h^{\la}_{13} \right)$
is equivalent to finding the joint probability mass function of
$\left( h^{\la}_{12},h^{\la}_{13} \right)$.

First, note that
$$\left( h^{\la}_{12},h^{\la}_{13} \right)=(0,0)\text{ iff }$$
$$\I(X_2 \in N(X_1)\cap \G_1(X_1,N))=\I(X_3 \in N(X_1)\cap \G_1(X_1,N))=0 \text{ iff }$$
$$\I(X_2 \in \TY\setminus N(X_1)\cap \G_1(X_1,N))=\I(X_3 \in \TY\setminus N(X_1)\cap \G_1(X_1,N))=1 \text{ iff }$$
$$\I(\{X_2,X_3 \} \subset \TY\setminus [N(X_1)\cap \G_1(X_1,N)])=1.$$

Hence $P(\left( h^{\la}_{12},h^{\la}_{13} \right)=(0,0))=P(\{X_2,X_3 \} \subset
\TY\setminus [N(X_1)\cap \G_1(X_1,N)])$.

Next, the pair $\left( h^{\la}_{12},h^{\la}_{13} \right)=(1,1)$ iff
$h^{\la}_{12}=h^{\la}_{13}=1$. So
$P\left(\left( h^{\la}_{12},h^{\la}_{13} \right)=(1,1) \right)=\E\left[ h^{\la}_{12}.h^{\la}_{13} \right]$.

Furthermore, by symmetry
$P\left(\left( h^{\la}_{12},h^{\la}_{13} \right)=(0,1)\right)=P(\left( h^{\la}_{12},h^{\la}_{13} \right)=(1,0)).$
So it follows that
\begin{eqnarray*}
P\left(\left( h^{\la}_{12},h^{\la}_{13} \right)=(0,1)\right)
&=&P\left(\left( h^{\la}_{12},h^{\la}_{13} \right)=(1,0)\right)\\
&=&\frac{1}{2}\,\left[1-\left(P\left(\left( h^{\la}_{12},h^{\la}_{13} \right)=(0,0)\right)+
    P\left(\left( h^{\la}_{12},h^{\la}_{13} \right)=(1,1)\right)\right)\right].
\end{eqnarray*}

\subsection{Relative Edge Density of the Graphs Based on PCDs}
\label{sec:OR-edge-density}
The relative edge density of the underlying graph, $G_{\lo}(D_n)$, based on the PCD, $D_n$,
is denoted as $\rho_{\lo}(D_n)$.f
For $X_i \stackrel{iid}{\sim}F$, $i=1,2,\ldots,n$,
one can write down the relative edge density as
$$\rho_{\lo}(D_n)=\frac{2}{n\,(n-1)}\sum\hspace*{-0.1 in}\sum_{i < j \hspace*{0.25 in}}   \hspace*{-0.1 in} \,h^{\lo}_{ij}$$
where
\begin{eqnarray*}
h^{\lo}_{ij} & = & \I(\left( X_i,X_j \right)\in \mE_{\lo})= \I((X_i,X_j) \in \A\text{ or }(X_j,X_i) \in \A)\\
&=&\I(X_j \in N(X_i)\text{ or }X_i \in N(X_j))=\I(X_j \in N(X_i)\text{ or }X_j \in \G_1\left( X_i,N \right))\\
&=&\I(X_j \in N(X_i)\cup \G_1\left( X_i,N \right))
\end{eqnarray*}
is the number of edges between $X_i$ and $X_j$ in $G_{\lo}(D_n)$.
Note that $h^{\lo}_{ij}$ is a symmetric kernel with finite variance since
$0 \le h^{\lo}_{ij} \le 1$.
Moreover, $\rho_{\lo}(D_n)$ is a
random variable that depends on $n$, $F$, and $N(\cdot)$ (i.e.,  $\mathcal Y$).
But $\E[\rho_{\lo}(D_n)]$ does only depend on $F$ and $N(\cdot)$.
That is,
\begin{equation}
\label{eqn:E[rho-or-D]}
0\le \E\left[ \rho_{\lo}(D_n) \right]=
\frac{2}{n\,(n-1)}\sum\hspace*{-0.1 in}\sum_{i < j \hspace*{0.25 in}}
\hspace*{-0.1 in} \,\E[h^{\lo}_{ij}]=\E\left[ h^{\lo}_{12} \right]
\end{equation}
where $\E\left[ h^{\lo}_{12} \right]=
\E[\I(X_2 \in N(X_1)\cup \G_1\left( X_1,N \right))]=
P(X_2 \in N(X_1)\cup \G_1(X_1,N))=p_{\lo}(N)$.
Notice that $p_{\lo}(N)$ is the edge probability for
the underlying graph $G_{\lo}(D_n)$
and that $p_{\lo}(N)=P(X_j \in N(X_i)\cup \G_1\left( X_i,N \right))$ for $i \not= j$.

Similar to the reflexivity graph case,
we have,
\begin{equation}
\label{eqn:Var[rho-or-D]}
0 \le \Var\left[ \rho_{\lo}(D_n) \right]=
\frac{2}{n\,(n-1)}\Var\left[ h^{\lo}_{12} \right]+
\frac{4\,(n-2)}{n\,(n-1)} \, \Cov[h^{\lo}_{12},h^{\lo}_{13}].
\end{equation}
As before, it follows that
$$\Var\left[ h^{\lo}_{12} \right]=p_{\lo}(N)-\left[p_{\lo}(N)\right]^2=
p_{\lo}(N)\left(1-p_{\lo}(N)\right).$$
and
$$\Cov\left[ h^{\lo}_{12},h^{\lo}_{13} \right]=
\E\left[ h^{\lo}_{12}.h^{\lo}_{13} \right]-\E\left[ h^{\lo}_{12} \right]\E\left[ h^{\lo}_{13} \right].$$
Since $\E\left[ h^{\lo}_{12} \right]=\E\left[ h^{\lo}_{13} \right]=p_{\lo}(N)$ and,
\begin{eqnarray*}
\E\left[ h^{\lo}_{12}.h^{\lo}_{13} \right]
& = & \E[\I(X_2 \in N(X_1)\cup \G_1(X_1,N))\, . \,\I( X_3 \in N(X_1)\cup \G_1(X_1,N))]\\
& = & \E[\I(X_2 \in N(X_1)\cup \G_1(X_1,N)\, , \, X_3 \in N(X_1)\cup \G_1(X_1,N))]\\
& = & P(X_2 \in N(X_1)\cup \G_1(X_1,N)\, , \, X_3 \in N(X_1)\cup \G_1(X_1,N))\\
& = & P(\{X_2,X_3\} \subset N(X_1)\cup \G_1(X_1,N)),
\end{eqnarray*}
it follows that
$$\Cov\left[ h^{\lo}_{12},h^{\lo}_{13} \right]=P(\{X_2,X_3\} \subset N(X_1)\cup \G_1(X_1,N))-\left[p_{\lo}(N)\right]^2.$$

The underlying graph $G_{\lo}(D_n)$ is a random graph
where the edge probability is $P(X_iX_j \in \mE_{\lo})=p_{\lo}(N)$ for $i\not=j$
and is an estimable parameter of degree 2.
Using Equation \eqref{eqn:E[rho-or-D]},
we obtain that $G_{\lo}(D_n)$ is an unbiased estimator of $p_{\lo}(N)$.
Notice that for the underlying and reflexivity graphs based on the PCDs,
the set of vertices $\V=\X_n$ is a random sample from a distribution $F$
(i.e., the vertices directly result from a random process),
and the edges are defined based on the random sets $N(X_i) \cup \G_1(X_1,N)$
as described before.
Hence the set of edges $\mE_{\lo}$ (indirectly) result from a random process
such that
$h^{\lo}_{ij}$ are identically distributed
and $h^{\lo}_{ij}$ and $h^{\lo}_{kl}$ are independent for distinct $i,j,k,l$.
Furthermore, we have $\Cov\left[ h^{\lo}_{ij},h^{\lo}_{kl} \right]<\infty$ as before.
Then we have the following corollary to the Main Result 1.
\begin{corollary}
\label{cor:edge-dens-or}
The relative edge density, $\rho_{\lo}(D_n)$, of the underlying graph,
$G_{\lo}(D_n)$, is a one-sample $U$-statistic of degree 2
and is an unbiased estimator of $p_{\lo}(N)$.
If, additionally, $\nu_{\lo}(N):=\Cov\left[ h^{\lo}_{ij},h^{\lo}_{ik} \right]>0$ for all $i \not=j\not=k$, $i,j,k =1,2,\ldots,n$,
then $\sqrt{n}\,\bigl[\rho_{\lo}(D_n)-p_{\lo}(N)]\stackrel{\mathcal L}{\longrightarrow} \N(0,4\,\nu_{\lo}(N))$
as $n \rightarrow \infty$.
\end{corollary}


In the above corollary,
$\nu_{\lo}(N) > 0$ iff
$P(\{X_2,X_3\} \subset N(X_1)\cup \G_1(X_1,N))>\left[p_{\lo}(N)\right]^2$.

\subsubsection{The Joint Distribution of $\left( h^{\lo}_{12},h^{\lo}_{13} \right)$}
Finding the joint distribution of $\left( h^{\lo}_{12},h^{\lo}_{13} \right)$ is
equivalent to finding the joint probability mass function of
$\left( h^{\lo}_{12},h^{\lo}_{13} \right)$, i.e.,  finding
$P(\left( h^{\lo}_{12},h^{\lo}_{13} \right)=(i,j)) \text{ for each } (i,j)\in
\{(0,0),(1,0),(0,1),(1,1)\}.$

First, note that
$$\left( h^{\lo}_{12},h^{\lo}_{13} \right)=(0,0)\text{ iff }$$
$$\I(X_2 \in N(X_1)\cup \G_1(X_1,N))=\I(X_3 \in N(X_1)\cup \G_1(X_1,N))=0 \text{ iff }$$
$$\I(X_2 \in \TY\setminus N(X_1)\cup \G_1(X_1,N))=\I(X_3 \in \TY\setminus N(X_1)\cup \G_1(X_1,N))=1 \text{ iff }$$
$$\I(\{X_2,X_3 \} \subset \TY\setminus [N(X_1)\cup \G_1(X_1,N)])=1.$$

Hence $P(\left( h^{\lo}_{12},h^{\lo}_{13} \right)=(0,0))=P(\{X_2,X_3 \} \subset
\TY\setminus [N(X_1)\cup \G_1(X_1,N)])$.

Next, note that $\left( h^{\lo}_{12},h^{\lo}_{13} \right)=(1,1)$ iff
$h^{\lo}_{12}=h^{\lo}_{13}=1$.
$P(\left( h^{\lo}_{12},h^{\lo}_{13} \right)=(1,1))=\E\left[ h^{\lo}_{12}.h^{\lo}_{13} \right]$.

By symmetry
$P(\left( h^{\lo}_{12},h^{\lo}_{13} \right)=(0,1))=P(\left( h^{\lo}_{12},h^{\lo}_{13} \right)=(1,0))$.
Hence
\begin{eqnarray*}
P(\left( h^{\lo}_{12},h^{\lo}_{13} \right)=(0,1)) & = & P(\left( h^{\lo}_{12},h^{\lo}_{13} \right)=(1,0)) \\
& = & \frac{1}{2}\,\left(1-\left[P(\left( h^{\lo}_{12},h^{\lo}_{13} \right)=(0,0))+
P(\left( h^{\lo}_{12},h^{\lo}_{13} \right)=(1,1))\right]\right).
\end{eqnarray*}

\begin{remark}
Note that $2\,h_{ij}=h^{\la}_{ij}+h^{\lo}_{ij}$,
since
if $g_{ij}=g_{ji}=0$, then $2\,h_{ij}=0$, and $h^{\la}_{ij}=h^{\lo}_{ij}=0$;
if $g_{ij}=g_{ji}=1$, then $2\,h_{ij}=2$, and $h^{\la}_{ij}=h^{\lo}_{ij}=1$;
and
if $g_{ij}=0$ and $g_{ji}=1$, then $2\,h_{ij}=1$, and $h^{\la}_{ij}=0$ and $h^{\lo}_{ij}=1$;
by symmetry, the same holds when $g_{ij}=1$ and $g_{ji}=0$.
$\square$
\end{remark}

\subsection{Proportional-Edge Proximity Maps and the Associated Regions}
\label{sec:prop-edge}
\begin{figure} [ht]
    \centering
    \scalebox{.4}{\input{Gammaofnu2.pstex_t}}
    \caption{Construction of the $\G_1$-region, $\G_1^{r=2}(x)$ (shaded region) for an $x \in R(\y_1)$. }
\label{fig:ProxMapDef2}
    \end{figure}
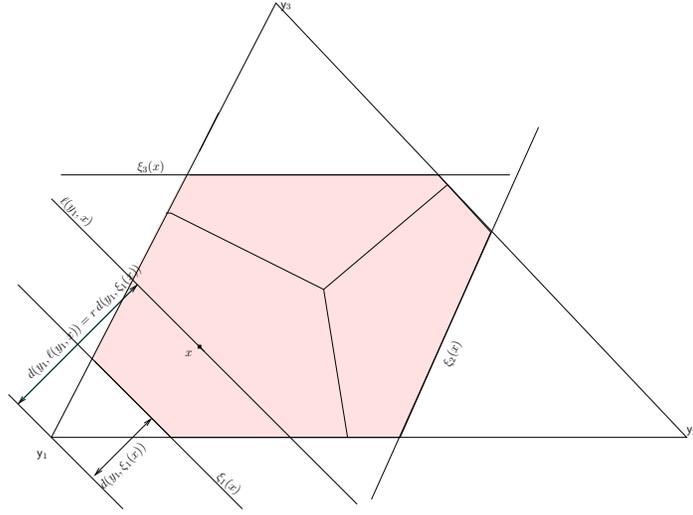

Let $\Omega = \mathbb{R}^2$ and $\Y_3 = \{\y_1,\y_2,\y_3\} \subset \mathbb{R}^2$ be three non-collinear points.
Denote by $\TY$ the triangle (including the interior) formed by these three points.
For $r \in [1,\infty]$
define $\NPE^r(x)$ to be the {\em proportional-edge}  proximity map with parameter $r$
and $\G_1^r(x):=\G_1\left( x,\NPE^r \right)$ to be the corresponding $\G_1$-region as follows;
see also Figures  \ref{fig:ProxMapDef1} and \ref{fig:ProxMapDef2}.
Let ``vertex regions'' $R(\y_1)$, $R(\y_2)$, $R(\y_3)$
partition $\TY$ using segments from the
center of mass of $\TY$ to the edge midpoints.
For $x \in \TY \setminus \Y_3$, let $v(x) \in \Y_3$ be the
vertex whose region contains $x$; $x \in R(v(x))$.
If $x$ falls on the boundary of two vertex regions,
or at the center of mass, we assign $v(x)$ arbitrarily.
Let $e(x)$ be the edge of $\TY$ opposite $v(x)$.
Let $\ell(v(x),x)$ be the line parallel to $e(x)$ through $x$.
Let $d(v(x),\ell(v(x),x))$ be the Euclidean (perpendicular) distance from $v(x)$ to $\ell(v(x),x)$.
For $r \in [1,\infty)$ let $\ell_r(v(x),x)$ be the line parallel to $e(x)$
such that
$$d(v(x),\ell_r(v(x),x)) = rd(v(x),\ell(v(x),x))
\text{ and }
d(\ell(v(x),x),\ell_r(v(x),x)) < d(v(x),\ell_r(v(x),x)).$$
Let $T_r(x)$ be the triangle similar to
and with the same orientation as $\TY$
having $v(x)$ as a vertex and $\ell_r(v(x),x)$ as the opposite edge.
Then the proportional-edge proximity region
$\NPE^r(x)$ is defined to be $T_r(x) \cap \TY$.

Furthermore, let $\xi_i(x)$ be the line such that $\xi_i(x)\cap \TY \not=\emptyset$ and
$r\,d(\y_i,\xi_i(x))=d(\y_i,\ell(\y_i,x))$ for $i=1,2,3$.
Then $\G_1^r(x)\cap R(\y_i)=\{z \in R(\y_i): d(\y_i,\ell(\y_i,z)) \ge d(\y_i,\xi_i(x)\}$, for $i=1,2,3$.
Hence  $\G_1^r(x)=\bigcup_{i=1}^3 (\G_1^r(x)\cap R(\y_i))$.
Notice that $r \ge 1$ implies $x \in \NPE^r(x)$ and $x \in \G_1^r(x)$.
Furthermore,
$\lim_{r \rightarrow \infty} \NPE^r(x) = \TY$ for all $x \in \TY \setminus \Y_3$,
and so we define $\NY^{\infty}(x) = \TY$ for all such $x$.
For $x \in \Y_3$, we define $\NPE^r(x) = \{x\}$ for all $r \in [1,\infty]$.
Then,
for $x \in R(\y_i)$,
$\lim_{r \rightarrow \infty} \G_1^r(x) = \TY \setminus \{\y_j,\y_k\}$
for distinct $i,j$, and $k$.

Notice that $X_i \stackrel{iid}{\sim} F$,
with the additional assumption
that the non-degenerate two-dimensional
probability density function $f$ exists
with support in $\TY$,
implies that the special cases in the construction of $\NPE^r$
---$X$ falls on the boundary of two vertex regions,
or at the center of mass, or $X \in \Y_3$ ---
occur with probability zero.
Note that for such an $F$,
$\NPE^r(x)$ is a triangle a.s.
and $\G_1^r(x)$ is a convex or nonconvex polygon.

\begin{figure} [ht]
    \centering
   \scalebox{.4}{\input{Nofnu2.pstex_t}}
\caption{Construction of proportional-edge proximity region,
$\NPE^{r=2}(x)$ (shaded region) for an $x \in R(\y_1)$. }
\label{fig:ProxMapDef1}
    \end{figure}
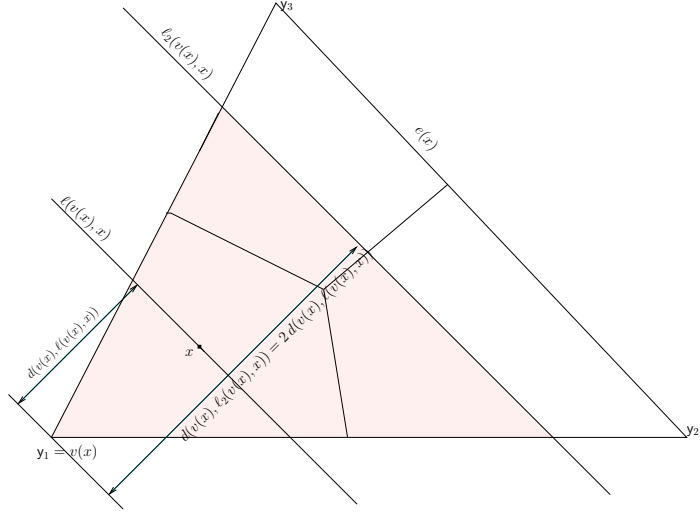

\subsection{Relative Edge Density of the Underlying and Reflexivity Graphs of Proportional-Edge PCDs}
\label{sec:relative-density-PEPCD}
Let $\X_n=\{X_1,X_2,\ldots,X_n\}$ be a sample from a distribution $F$
with support in $\TY$.
Let $D_n(r)$ be the proportional edge PCD with
vertex set $\V=\X_n$ and arc set $\A$
defined by $\left( X_i,X_j \right) \in \A \iff X_j \in \NPE^r(X_i)$.
Consider the underlying and reflexivity graphs of the vertex-random PCD, $D_n(r)$.
Recall that $X_iX_j \in \mE_{\la}$ iff $X_j \in \NPE^r(X_i) \cap
\G_1\left( X_i,\NPE^r \right)$ and $X_iX_j \in \mE_{\lo}$ iff $X_j \in \NPE^r(X_i)
\cup  \G_1\left( X_i,\NPE^r \right)$.

Let $h^{\la}_{ij}(r):=\I(X_iX_j \in \mE_{\la})=\I(X_j \in \NPE^r(X_i)\cap \G^r_1\left( X_i \right))$
and $h^{\lo}_{ij}(r):=\I(X_iX_j \in \mE_{\lo})=\I(X_j \in \NPE^r(X_i)\cup \G^r_1\left( X_i \right))$
for $i \not= j$.
The relative edge density $\rho^{\la}_n(r) := \rho_{\la}(D_n(r))$
depends on $n$ explicitly, and on $F$ and $\NPE^r$ implicitly.
The expectation $\E\left[\rho^{\la}_n(r)\right]$, however, is independent of $n$
and depends on only $F$ and $\NPE^r$.
Let $p_{\la}(F,r):=\E\left[ h^{\la}_{12}(r) \right]$
and
$\nu_{\la}(F,r):=\Cov\left[ h^{\la}_{12}(r),h^{\la}_{13}(r) \right]$.
Then
\begin{eqnarray}
0 \le \E\left[ \rho^{\la}_n(r) \right] = \E\left[ h^{\la}_{12}(r) \right] \le 1.
\end{eqnarray}
The variance $\Var\left[ \rho^{\la}_n(r) \right]$ simplifies to
\begin{eqnarray}
0 \leq
  \Var\left[ \rho^{\la}_n(r) \right] =
     \frac{2}{n(n-1)} \Var\left[ h^{\la}_{12}(r) \right] +
     \frac{4\,(n-2)}{n(n-1)} \Cov\left[ h^{\la}_{12}(r),h^{\la}_{13}(r) \right]
  \leq 1/4.
\end{eqnarray}
By Theorem \ref{thm:edge-dens-result}, it follows that
\begin{eqnarray}
\sqrt{n}\left(\rho^{\la}_n(r)-p_{\la}(F,r) \right)
\stackrel{\mathcal{L}}{\longrightarrow}
\N\left(0,4\,\nu_{\la}(F,r) \right)
\end{eqnarray}
provided $\nu_{\la}(F,r) > 0$.
The asymptotic variance of $\rho^{\la}_n(r)$ is $4\,\nu_{\la}(F,r)$ and
depends on only $F$ and $\NPE^r$.
Thus we need determine only $p_{\la}(F,r)$ and
$\nu_{\la}(F,r)$ in order to obtain the normal
approximation
\begin{eqnarray}
\label{eqn:asy-norm-and}
\rho^{\la}_n(r) \stackrel{\text{approx}}{\sim}
\N\left(p_{\la}(F,r),\frac{4\,\nu_{\la}(F,r)}{n}\right).
\end{eqnarray}

The above paragraph holds for
$\rho^{\lo}_n(r)=\rho_{\lo}(D_n(r))$ also with $\rho^{\la}_n(r)$ is
replaced by $\rho^{\lo}_n(r)$, $h^{\la}_{12}(r)$ and $h^{\la}_{13}(r)$ are
replaced by $h^{\lo}_{12}(r)$ and $h^{\lo}_{13}(r)$, respectively.


For $r=1$, $\NPE^{r=1}(x) \cap \G_1^{r=1}(x)=\ell(v(x),x)$
which has zero $\R^2$-Lebesgue measure.
Then we have
$$\E\left[ \rho^{\la}_n(r=1) \right] = \E\left[ h^{\la}_{12}(r=1) \right]=
\mu_{\la}(r=1)=
P(X_2 \in \NPE^{r=1}(X_1) \cap \G_1^{r=1}(X_1))=0.$$
Similarly,
$P(\{X_2,X_3\} \subset \NPE^{r=1}(X_1) \cap \G_1^{r=1}(X_1))=0$.
Thus, $\nu_{\la}(r=1)=0$.
Furthermore, for $r=\infty$,
$\NPE^{r=\infty}(x) \cap \G_1^{r=\infty}(x)=\TY$ for all $x \in \TY \setminus \Y_3$.
Then
$$\E\left[ \rho^{\la}_n(r=\infty) \right] =
\E\left[ h^{\la}_{12}(r=\infty) \right]=
\mu_{\la}(r=\infty)=
P(X_2 \in \NPE^{r=\infty}(X_1) \cap \G_1^{r=\infty}(X_1)=
P(X_2 \in \TY)=1.$$
Similarly, $P(\{X_2,X_3\} \subset \NPE^{r=\infty}(X_1) \cap \G_1^{r=\infty}(X_1))=1$.
Hence $\nu_{\la}(r=\infty)=0$.
Therefore, the CLT result in Equation \eqref{eqn:asy-norm-and} does not hold
for $r \in \{1,\infty\}$.
Furthermore, $\rho^{\la}_n(r=1)=0$ a.s. and $\rho^{\la}_n(r=\infty)=1$ a.s.
For $r \in (1,\infty)$,
since $h^{\la}_{12}(r) = \I(X_2 \in \NPE^r(X_1)\cap \G_1^r(X_1)$
is the number of edges in the reflexivity graph,
$h^{\la}_{12}(r)$ tends to be high if the intersection region is large.
In such a case, $h^{\la}_{13}(r)$ tends to be high also.
That is, $h^{\la}_{12}(r)$ and $h^{\la}_{13}(r)$
tend to be high and low together.
So, for $r \in (1,\infty)$,
we have $\nu_{\la}(F,r)>0$.
See also Figure \ref{fig:asymptotics} (right)
and Appendix 1.

For $r=1$, $\NPE^{r=1}(x) \cup \G_1^{r=1}(x)$
has positive $\R^2$-Lebesgue measure.
Then
$P(\{X_2,X_3\} \subset \NPE^{r=1}(X_1) \cup \G_1^{r=1}(X_1))>0$.
Thus, $\nu_{\lo}(r=1) \not= 0$.
On the other hand, for $r=\infty$,
$\NPE^{r=\infty}(X_1) \cup \G_1^{r=\infty}(X_1))=\TY$ for all $X_1 \in \TY$.
Then
$$\E\left[ \rho^{\lo}_n(r=\infty) \right] = \E\left[ h^{\lo}_{12}(r=\infty) \right]=
P(X_2 \in \NPE^{r=\infty}(X_1) \cup \G_1^{r=\infty}(X_1))=
\mu_{\lo}(r=\infty)= P(X_2 \in \TY)=1.$$
Similarly, $P(\{X_2,X_3\} \subset \NPE^{r=\infty}(X_1) \cup \G_1^{r=\infty}(X_1))=1$.
Hence $\nu_{\lo}(r=\infty)=0$.
Therefore, the CLT result for the underlying graph case does not hold for $r =\infty$.
Moreover $\rho^{\lo}_n(r=\infty)=1$ a.s.
For $r \in [1,\infty)$,
since $h^{\lo}_{12}(r) = \I(X_2 \in \NPE^r(X_1)\cup \G_1^r(X_1)$
is the number of edges in the underlying graph,
$h^{\lo}_{12}(r)$ tends to be high if the union region is large.
In such a case, $h^{\lo}_{13}(r)$ tends to be high also.
That is, $h^{\lo}_{12}(r)$ and $h^{\lo}_{13}(r)$
tend to be high and low together.
So, for $r \in [1,\infty)$,
we have $\nu_{\lo}(F,r)>0$.
See also Figure \ref{fig:asymptotics} (right)
and Appendix 2.

\begin{remark}
\label{rem:arc-density}
\textbf{Relative Arc Density of Proportional-Edge PCDs:}

Let $h_{ij}(r):=\I( (X_i,X_j) \in \A)=\I(X_j \in \NPE^r(X_i))$ for $i \not= j$
and the relative arc density $\rho_n(r) := \rho(D_n(r))$.
Let $p(r):=\E\left[ \rho_n(r) \right]$
and
$4\,\nu(r):=\Cov\left[ h_{12}(r),h_{13}(r) \right]$.
By Theorem \ref{thm:arc-dens-result}, we have
\begin{eqnarray}
\sqrt{n}\left(\rho_n(r)-p(r)\right)
\stackrel{\mathcal{L}}{\longrightarrow}
\N\left(0,4\,\nu(r)\right)
\end{eqnarray}
provided $\nu(r) > 0$.
The explicit forms of asymptotic mean $p(r)$
and variance $4\,\nu(r)$ for uniform data are provided in \cite{ceyhan:arc-density-PE}.
$\square$
\end{remark}

\section{Asymptotic Distribution of Relative Edge Density for Uniform Data}
\label{sec:null-dist-edge-density}
Let $X_i \stackrel{iid}{\sim} \mathcal{U}(\TY)$ for $i=1,2,\ldots,n$,
where $\mathcal{U}(\TY)$ is the the uniform distribution on the triangle $\TY$.

We first present a ``geometry invariance" result which will simplify our subsequent analysis
by allowing us to consider the special case of the equilateral triangle.
Let $\rho^{\la}_n(r) := \rho_{\la}(\mathcal{U}(\TY),r)$ and
$\rho^{\lo}_n(r) := \rho_{\lo}(\mathcal{U}(\TY),r)$.
\begin{theorem}
\label{thm:geo-inv-NYr-under}
\textbf{Geometry Invariance:}
Let $\Y_3 = \{\y_1,\y_2,\y_3\} \subset \mathbb{R}^2$
be three non-collinear points.
For $i=1,2,\ldots,n$,
let $X_i \stackrel{iid}{\sim} \mathcal{U}(\TY)$.
Then for any $r \in [1,\infty]$
the distribution of $\rho^{\la}_n(r)$ and $\rho^{\lo}_n(r)$
is independent of $\Y_3$,
and hence the geometry of $\TY$.
\end{theorem}

\noindent \textbf{Proof:}
A composition of translation, rotation, reflections, and scaling
will take any given triangle $T_o = T(\y_1,\y_2,\y_3)$
to the ``basic'' triangle $T_b = T((0,0),(1,0),(c_1,c_2))$
with $0 < c_1 \le 1/2$, $c_2 > 0$ and $(1-c_1)^2+c_2^2 \le 1$,
preserving uniformity.
The transformation $\phi: \mathbb{R}^2 \rightarrow \mathbb{R}^2$
given by $\phi(u,v) = \left( u+\frac{1-2\,c_1}{2\,c_2}\,v,\frac{\sqrt{3}}{2\,c_2}\,v \right)$
takes $T_b$ to
the equilateral triangle
$T_e = T\left((0,0),(1,0),\left(1/2,\sqrt{3}/2\right)\right)$.
Investigation of the Jacobian shows that $\phi$
also preserves uniformity.
Furthermore, the composition of $\phi$ with the rigid motion transformations and scaling
maps
     the boundary of the original triangle $T_o$
  to the boundary of the equilateral triangle $T_e$,
     the median lines of $T_o$
  to the median lines of $T_e$,
and  lines parallel to the edges of $T_o$
  to lines parallel to the edges of $T_e$.
(A median line in a triangle is the line joining a vertex with the center of mass.)
Since the joint distribution of any collection of the $h^{\la}_{ij}(r)$ and $h^{\lo}_{ij}(r)$
involves only probability content of unions and intersections
of regions bounded by precisely such lines,
and the probability content of such regions is preserved since uniformity is preserved,
the desired result follows.
$\blacksquare$

Based on Theorem \ref{thm:geo-inv-NYr-under},
for our proportional-edge proximity map and the uniform data,
we may assume that
$\TY$ is a standard equilateral triangle, $T_e$,
with vertices $\Y_3 = \{(0,0),(1,0),(1/2,\sqrt{3}/2)\}$,
henceforth.

In the case of the (proportional-edge proximity map, uniform data) pair,
the asymptotic distribution of
$\rho^{\la}_n(r)$ and
$\rho^{\lo}_n(r)$
as a function of $r$ can be derived.
Recall that
$p_{\la}(r)=\E\left[ h^{\la}_{12}(r) \right]=P(X_2 \in \NPE^r(X_1)\cap \G_1^r(X_1))$ and
$p_{\lo}(r)=\E\left[ h^{\lo}_{12} \right]=P(X_2 \in \NPE^r(X_1)\cup \G_1^r(X_1))$
are the edge probabilities in the reflexivity and underlying graphs, respectively.

\begin{theorem}
\label{thm:asy-norm-under}
\textbf{(Main Result 3)}
For $r \in (1,\infty)$,
$$\sqrt{n}\left(\rho^{\la}_n(r)-p_{\la}(r)\right)\Big/\sqrt{4\,\nu_{\la}(r)}
\stackrel{\mathcal{L}}{\longrightarrow}  \mathcal{N}(0,1)$$
and
for $r \in [1,\infty)$,
$$\sqrt{n}\left(\rho^{\lo}_n(r)-p_{\lo}(r)\right)\Big/\sqrt{4\,\nu_{\lo}(r)}
\stackrel{\mathcal{L}}{\longrightarrow}  \mathcal{N}(0,1).$$
\end{theorem}
where the asymptotic means are
\begin{eqnarray}
\label{eqn:Asymean_and}
p_{\la}(r) =
 \begin{cases}
{\frac{(1-r)(5\,r^5-148\,r^4+245\,r^3-178\,r^2-232\,r+128)}{54\,r^2(r+2)(r+1)}} &\text{for} \quad r \in [1,4/3), \\
  -{\frac{101\,r^5-801\,r^4+1302\,r^3-732\,r^2-536\,r+672}{216\,r(r+2)(r+1)}}                                 &\text{for} \quad r \in [4/3,3/2), \\
  {\frac{r^8-13\,r^7+30\,r^6+148\,r^5-448\,r^4+264\,r^3+288\,r^2-368\,r+96}{8\,r^4(r+2)(r+1)}}  &\text{for} \quad r \in [3/2,2), \\
  \frac{(r^3+3\,r^2-2+2\,r)(-1+r)^2}{r^4(r+1)} &\text{for} \quad r \in [2,\infty),
 \end{cases}
\end{eqnarray}
\begin{eqnarray}
\label{eqn:Asymean_or}
p_{\lo}(r) =
\begin{cases}
           \frac{47\,r^6-195\,r^5+860\,r^4-846\,r^3-108\,r^2+720\,r-256}{108\,r^2(r+2)(r+1)}&\text{for} \quad r \in [1,4/3),\\
           \frac{175\,r^5-579\,r^4+1450\,r^3-732\,r^2-536\,r+672}{216\,r\,(r+2)(r+1)} &\text{for} \quad r \in [4/3,3/2),\\
           -\frac{3\,r^8-7\,r^7-30\,r^6+84\,r^5-264\,r^4+304\,r^3+144\,r^2-368\,r+96}{8\,r^4(r+2)(r+1)} &\text{for} \quad r \in [3/2,2),\\
           \frac{r^5+r^4-6\,r+2}{r^4(r+1)} &\text{for} \quad r \in [2,\infty),\\
         \end{cases}
\end{eqnarray}
and the asymptotic variances are
\begin{equation}
\label{eqn:Asyvar_and}
\nu_{\la}(r)=\sum_{i=1}^{11}\vartheta^{\la}_i(r)\,\I(\mI_i),
\end{equation}
\begin{equation}
\label{eqn:Asyvar_or}
\nu_{\lo}(r)=\sum_{i=1}^{11}\vartheta^{\lo}_i(r)\,\I(\mI_i).
\end{equation}
The explicit forms of $\vartheta^{\la}_i(r)$ and $\vartheta^{\lo}_i(r)$ are provided in Appendix Sections 1 and 2,
and the derivations of $p_{\la}(r)$, $\nu_{\la}(r)$,
$p_{\lo}(r)$, and $\nu_{\lo}(r)$ are provided in Appendix Section 3.

\begin{figure} [ht]
\centering
\scalebox{.3}{\input{means3.pstex_t}}
\scalebox{.3}{\input{var3.pstex_t}}
\caption{
\label{fig:asymptotics}
Result of Theorem \ref{thm:asy-norm-under}:
asymptotic null means $p(r)$, $p_{\la}(r)$, and $p_{\lo}(r)$ (left) and
variances $4\,\nu(r)$, $4\,\nu_{\la}(r)$, and $4\,\nu_{\lo}(r)$ (right)
for $r \in [1,5]$.
Some values of note:
$\mu(1) = 37/216$, $\mu_{\la}(1) =0$, and $\mu_{\lo}(1) = 37/108$,
$\lim_{r \rightarrow \infty} p(r)=
\lim_{r \rightarrow \infty} p_{\la}(r)=\lim_{r \rightarrow \infty} p_{\lo}(r) = 1$,
$4\,\nu_{\la}(r=1)=0$ and $\lim_{r \rightarrow \infty}4\,\nu_{\la}(r)=0$,
$4\,\nu_{\lo}(r=1)=1/3240$ and $\lim_{r \rightarrow \infty}4\,\nu_{\lo}(r)=0$,
and $\argsup_{r \in [1,\infty]} 4\,\nu(r) \approx 2.045$ with
$\sup_{r \in [1,\infty]} 4\,\nu(r) \approx .1305$,
$\argsup_{r \in [1,\infty]}4\,\nu_{\la}(r) \approx 2.69$ with
$\sup_{r \in [1,\infty]} 4\,\nu_{\la}(r) \approx .0537$,
$\argsup_{r \in [1,\infty]} 4\,\nu_{\lo}(r) \approx 1.765$ with
$\sup_{r \in [1,\infty]} 4\,\nu_{\lo}(r) \approx .0318$.}
\end{figure}
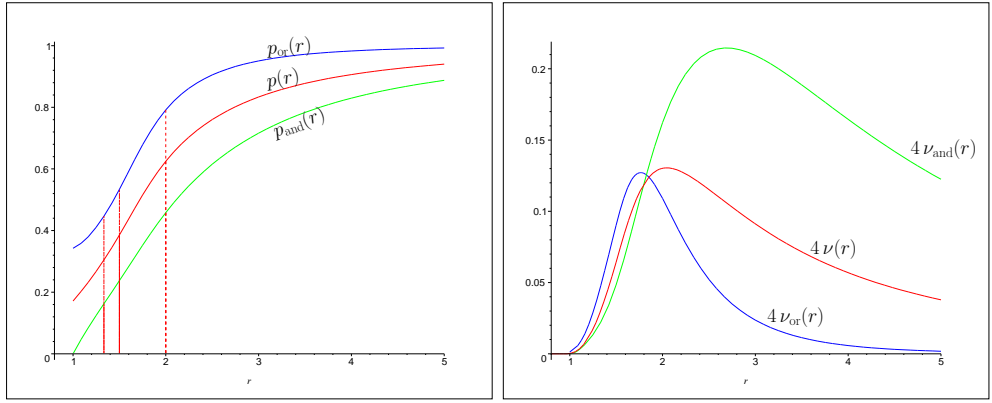

\begin{figure}
\centering
\psfrag{r}{\normalsize{$r$}}
\epsfig{figure=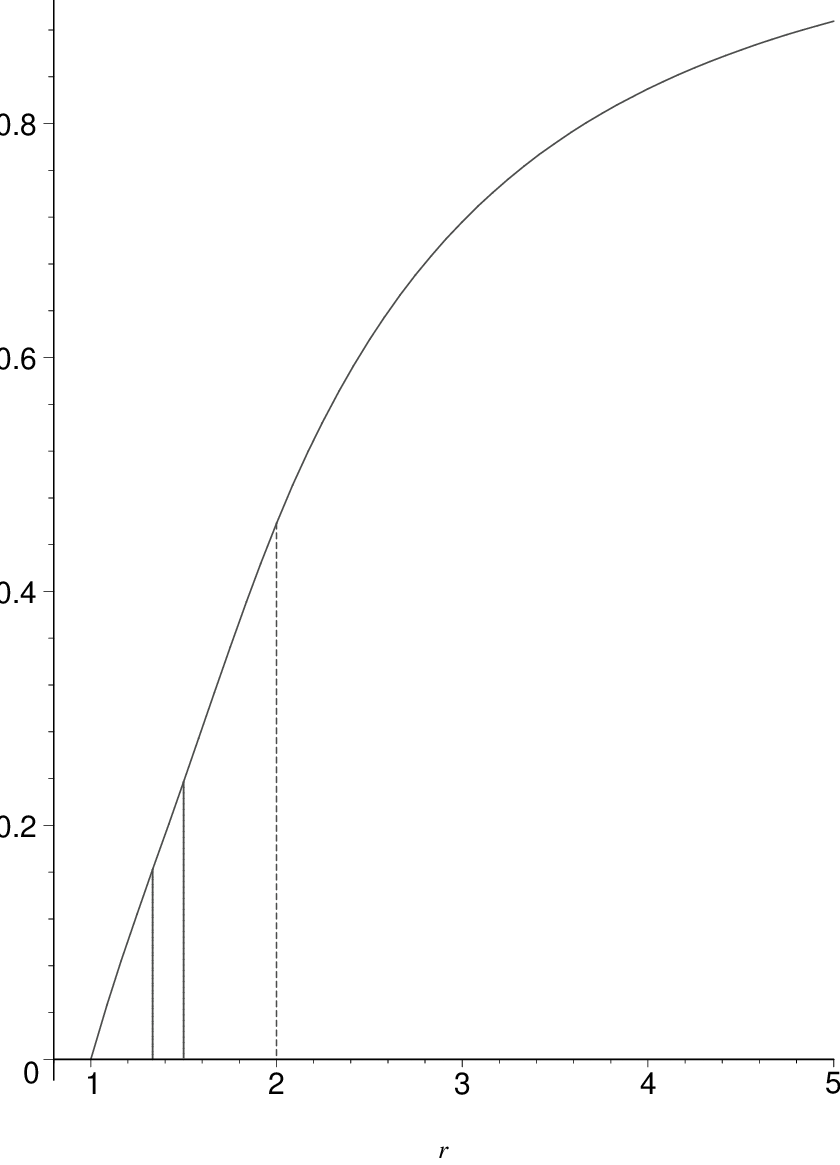, height=150pt, width=200pt}
\epsfig{figure=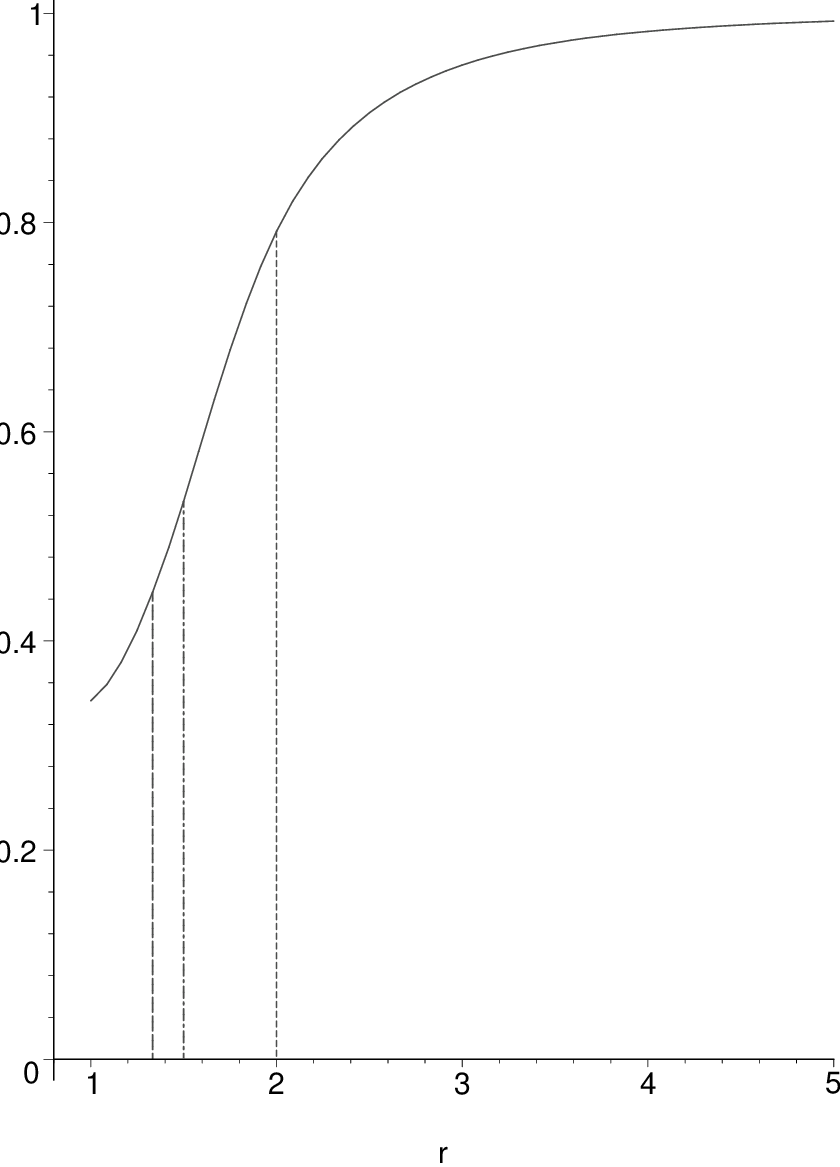, height=150pt, width=200pt}\\
\caption{
\label{fig:p_sa and p_or of nu}
The edge probabilities $p_{\la}(r)$ (left) and $p_{\lo}(r)$ (right) for $r \in [1,5]$.}
\end{figure}

The expectation $\E\left[ h^{\la}_{12}(r) \right]=p_{\la}(r)$
is as in Equation \eqref{eqn:Asymean_and};
and
$\E\left[ h^{\lo}_{12}(r) \right]=p_{\lo}(r)$
is as in Equation \eqref{eqn:Asymean_or} (see Figure \ref{fig:p_sa and p_or of nu}.
Notice that $\mu_{\la}(r=1)=0$ and $\lim_{r\rightarrow \infty}p_{\la}(r)=1$ (at rate $O(r^{-1})$);
and $\mu_{\lo}(r=1)=37/108$ and $\lim_{r\rightarrow \infty}p_{\lo}(r)=1$ (at rate $O(r^{-1})$).

To illustrate the limiting distribution, for example, $r=2$ yields
$$\frac{\sqrt{n}(\rho^{\la}_n(2)-\mu_{\la}(2))}{\sqrt{4\,\nu_{\la}(2)}}=
\sqrt{\frac{362880\,n}{58901}} \left(\rho^{\la}_n(2) -
\frac{11}{24}\right)\stackrel{\mathcal{L}}{\longrightarrow}
\N(0,1)$$ and
$$\frac{\sqrt{n}(\rho^{\lo}_n(2)-\mu_{\lo}(2))}{\sqrt{4\,\nu_{\lo}(2)}}=
\sqrt{\frac{120960\,n}{13189}} \left(\rho^{\lo}_n(2) -
\frac{19}{24}\right)\stackrel{\mathcal{L}}{\longrightarrow} \N(0,1);$$
or equivalently,
$$\rho^{\la}_n(2) \stackrel{\text{approx}}{\sim} \N\left(\frac{11}{24},\frac{58901}{362880\,n}\right)
\text{  and  } \rho^{\lo}_n(2) \stackrel{\text{approx}}{\sim}
\N\left(\frac{19}{24},\frac{13189}{120960\,n}\right).$$

By construction of the underlying and reflexivity graphs,
there is a natural ordering of the means of relative arc and edge densities.
\begin{lemma}
The means of the relative edge densities and arc density (i.e., the edge and arc probabilities)
have the following ordering:
$p_{\la}(r) < p(r) < p_{\lo}(r)$ for all $r \in [1,\infty)$.
Furthermore, for $r=\infty$,
we have $p_{\la}(r) = p(r) = p_{\lo}(r)=1$.
\end{lemma}

\noindent \textbf{Proof:}
Recall that $p_{\la}(r) = \E[\rho_n^{\la}(r)]=P(X_2 \in \NPE^r(X_1) \cap \G_1^r(X_1))$,
$p(r) = \E[\rho_n(r)]=P(X_2 \in \NPE^r(X_1))$, and
$p_{\lo}(r) = \E[\rho_n^{\lo}(r)]=P(X_2 \in \NPE^r(X_1) \cup \G_1^r(X_1))$.
And $\left[ \NPE^r(X_1) \cap \G_1^r(X_1) \right] \subseteq \NPE^r(X_1) \subseteq \left[ \NPE^r(X_1) \cup \G_1^r(X_1) \right]$
with probability 1 for all $r \ge 1$ with equality holding for $r=\infty$ only.
Then the desired result follows
(See also Figure \ref{fig:asymptotics}).
$\blacksquare$

Note that the above lemma holds for all $X_i$ that has a continuous distribution on $\TY$.
There is also a stochastic ordering for the relative edge and
arc densities as follows.
\begin{theorem}
For sufficiently small $r$, $\rho_n^{\la}(r) <^{ST} \rho_n(r) <^{ST} \rho_n^{\lo}(r)$ as $n \rightarrow \infty$.
\end{theorem}

\noindent \textbf{Proof:}
Above we have proved that $p_{\la}(r) < p(r) < p_{\lo}(r)$ for all $r \in [1,\infty)$.
For small $r$ ($r \lesssim 1.8$) the asymptotic variances have the same ordering,
$4\,\nu_{\la}(r) < 4\,\nu(r) < 4\,\nu_{\lo}(r)$.
Since $\rho_n^{\la}(r), \, \rho_n(r),\, \rho_n^{\lo}(r)$ are asymptotically normal,
then the desired result follows
(See also Figure \ref{fig:asymptotics}).
$\blacksquare$

\begin{figure}[ht]
\centering
\psfrag{Density}{\Huge{density}}
\rotatebox{-90}{ \resizebox{2.0 in}{!}{ \includegraphics{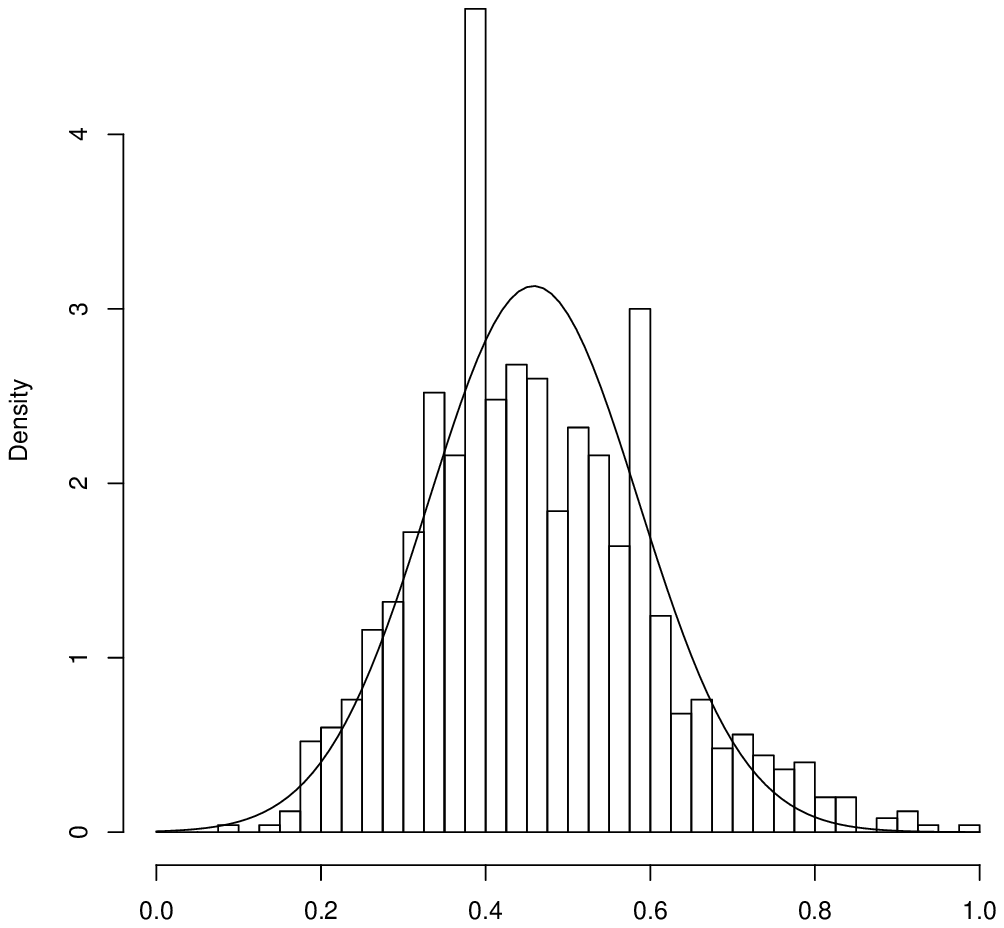} } }
\rotatebox{-90}{ \resizebox{2.0 in}{!}{ \includegraphics{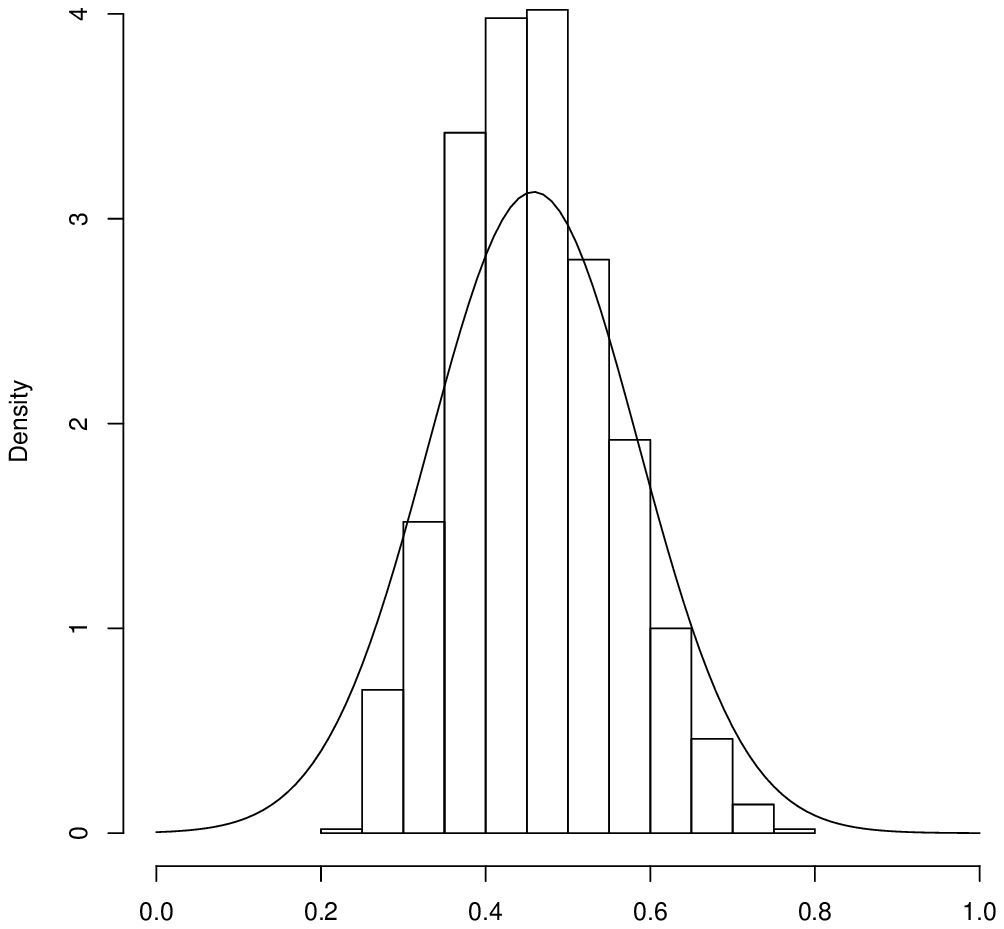} } }
\rotatebox{-90}{ \resizebox{2.0 in}{!}{ \includegraphics{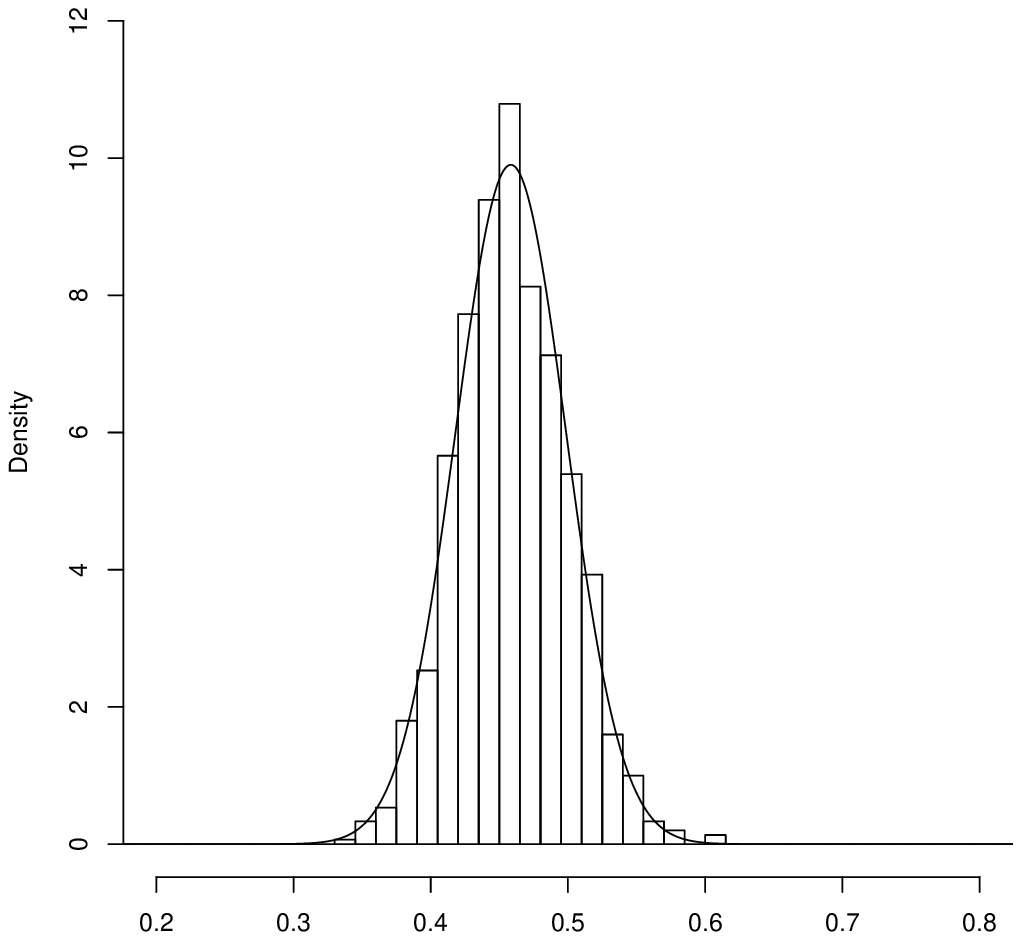} } }
\caption{
\label{fig:NormApprox_1}
Depicted are $\rho^{\la}_n(2) \stackrel{\text{approx}}{\sim} \mathcal{N}\left(\frac{11}{24},\frac{58901}{362880\,n}\right)$
for $n=10,\,20,\,100$ (left to right).
Histograms are based on 1000 Monte Carlo replicates.
Solid lines are the corresponding normal densities.
Notice that the axes are differently scaled.
}
\end{figure}

\begin{figure}[ht]
\centering
\psfrag{Density}{\Huge{density}}
\rotatebox{-90}{ \resizebox{2.0 in}{!}{ \includegraphics{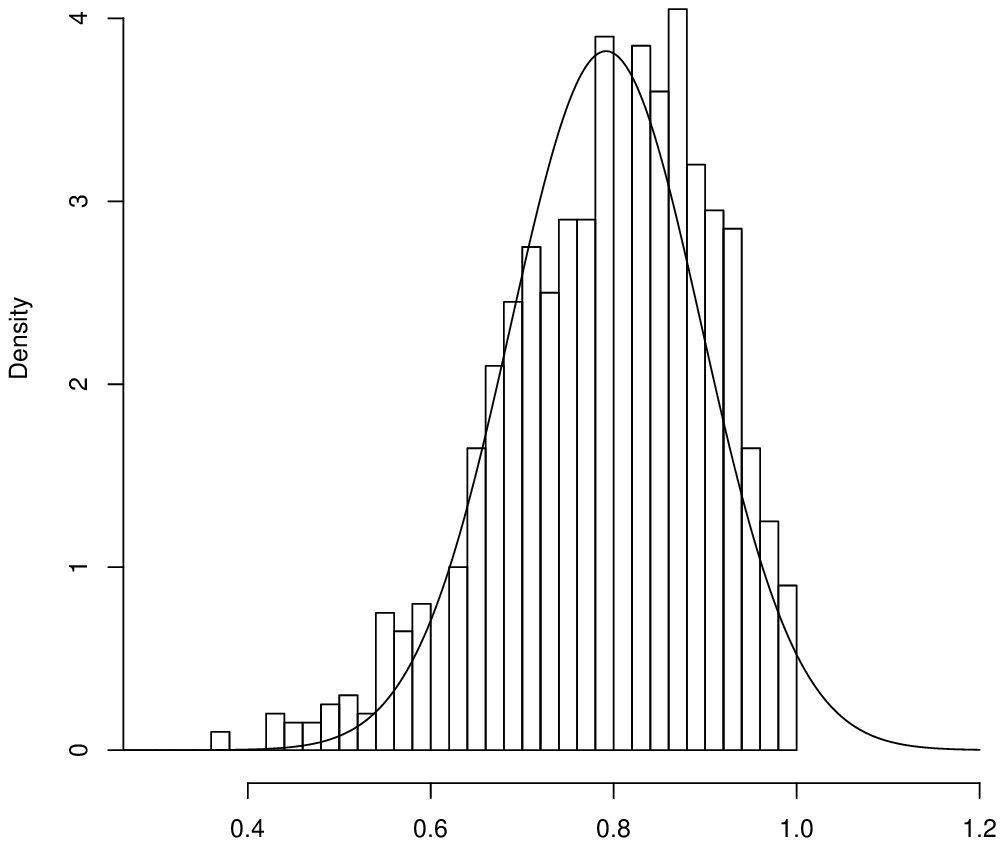} } }
\rotatebox{-90}{ \resizebox{2.0 in}{!}{ \includegraphics{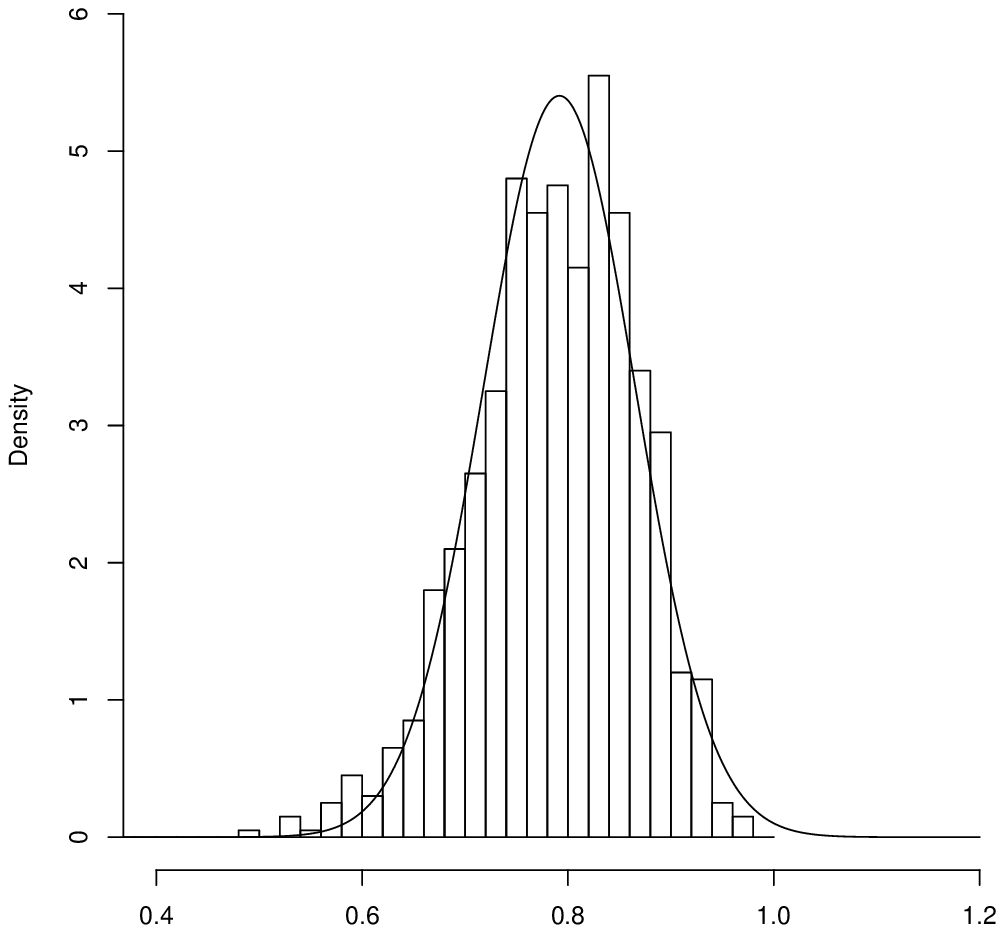} } }
\rotatebox{-90}{ \resizebox{2.0 in}{!}{ \includegraphics{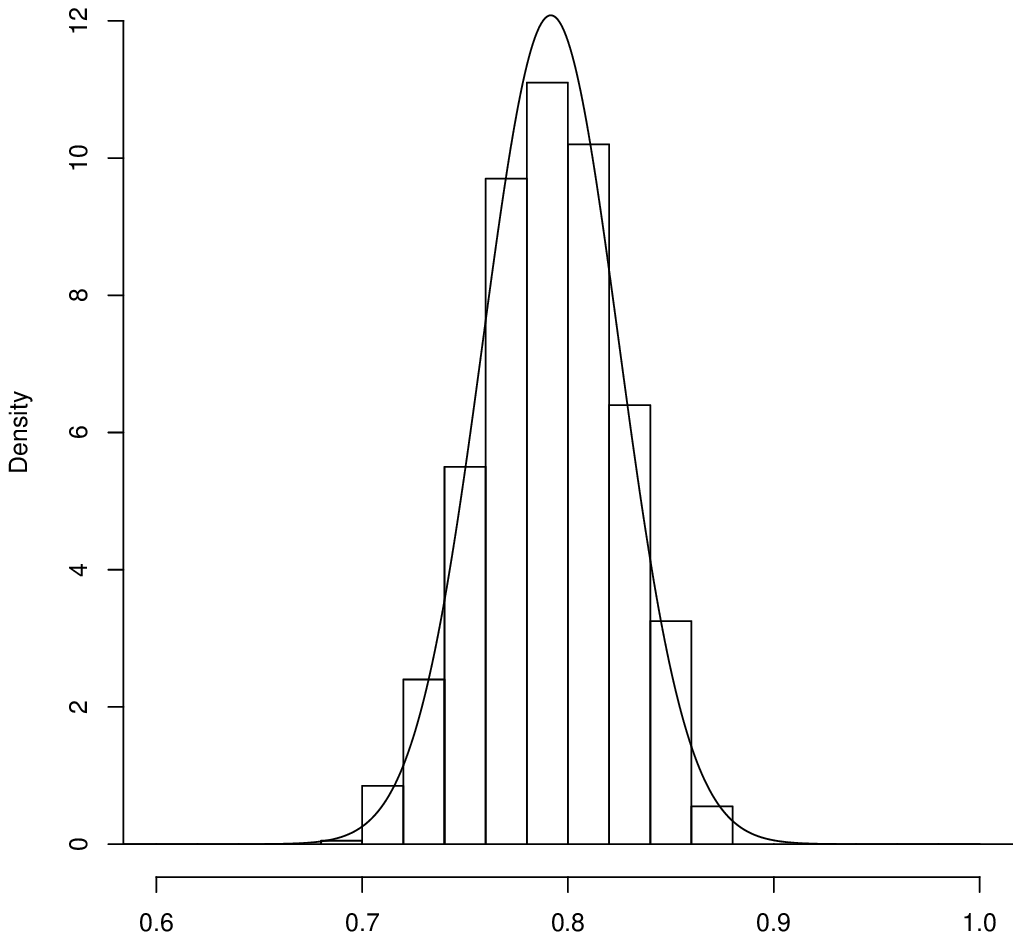} } }
\caption{
\label{fig:NormApprox2}
Depicted are $\rho^{\lo}_n(2) \stackrel{\text{approx}}{\sim} \mathcal{N}\left(\frac{19}{24},\frac{13189}{120960\,n}\right)$
for $n=10,20,100$ (left to right).
Histograms are based on 1000 Monte Carlo replicates.
Solid lines are the corresponding normal densities.
Notice that the axes are differently scaled.
}
\end{figure}

\begin{figure}[ht]
\centering
\psfrag{Density}{\Huge{density}}
\rotatebox{-90}{ \resizebox{2.0 in}{!}{ \includegraphics{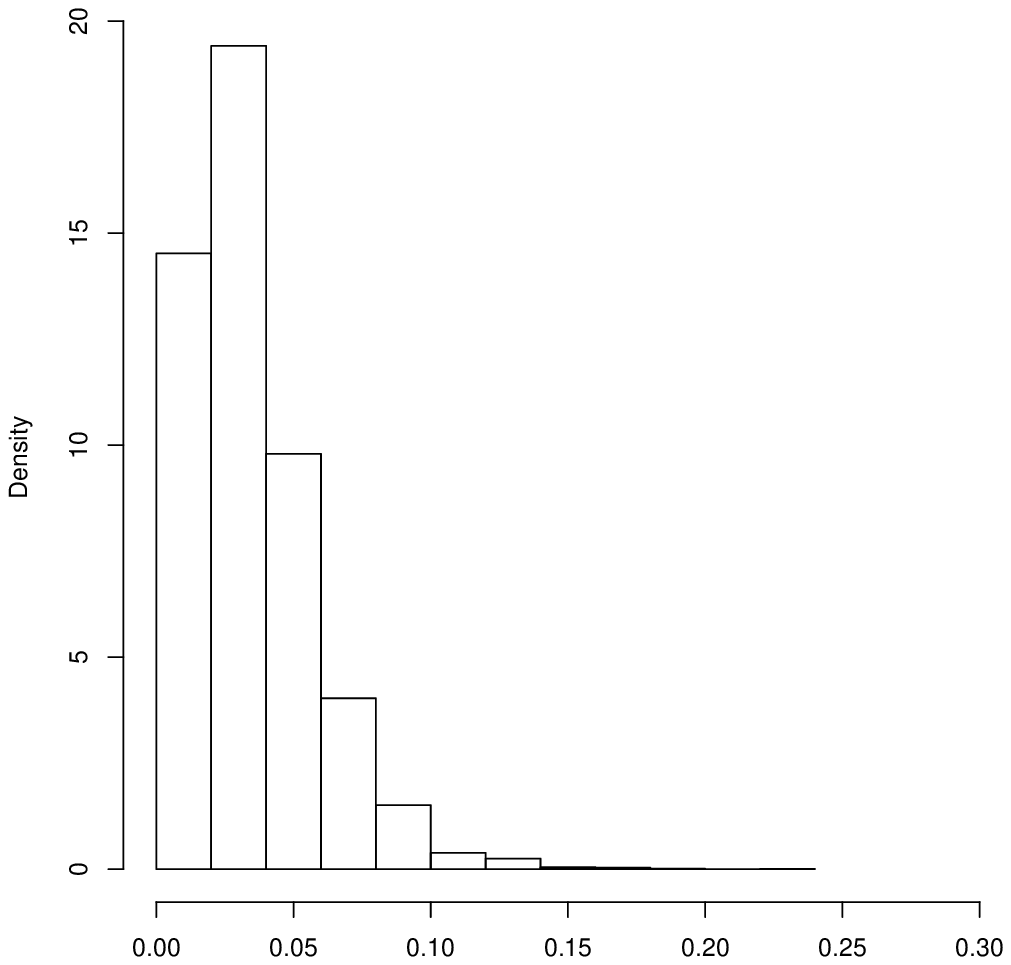} } }
\rotatebox{-90}{ \resizebox{2.0 in}{!}{ \includegraphics{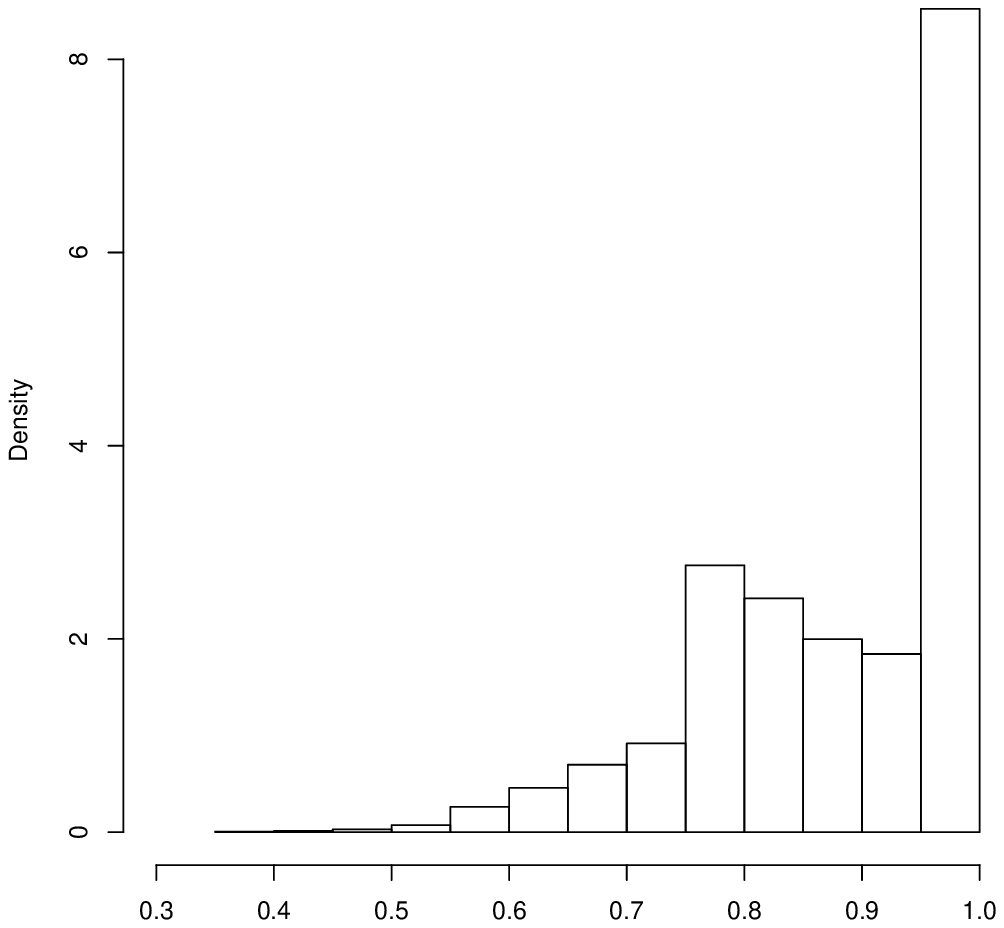} } }
\caption{
\label{fig:ANDskew}
Depicted are the histograms for 10000 Monte Carlo replicates of $\rho^{\la}_{10}(1.05)$ (left)
and $\rho^{\la}_{10}(5)$ (right) indicating severe small sample skewness for extreme values of $r$.
Notice that the vertical axes are differently scaled.
}
\end{figure}

\begin{figure}[ht]
\centering
\psfrag{Density}{\Huge{density}}
\rotatebox{-90}{ \resizebox{2.0 in}{!}{ \includegraphics{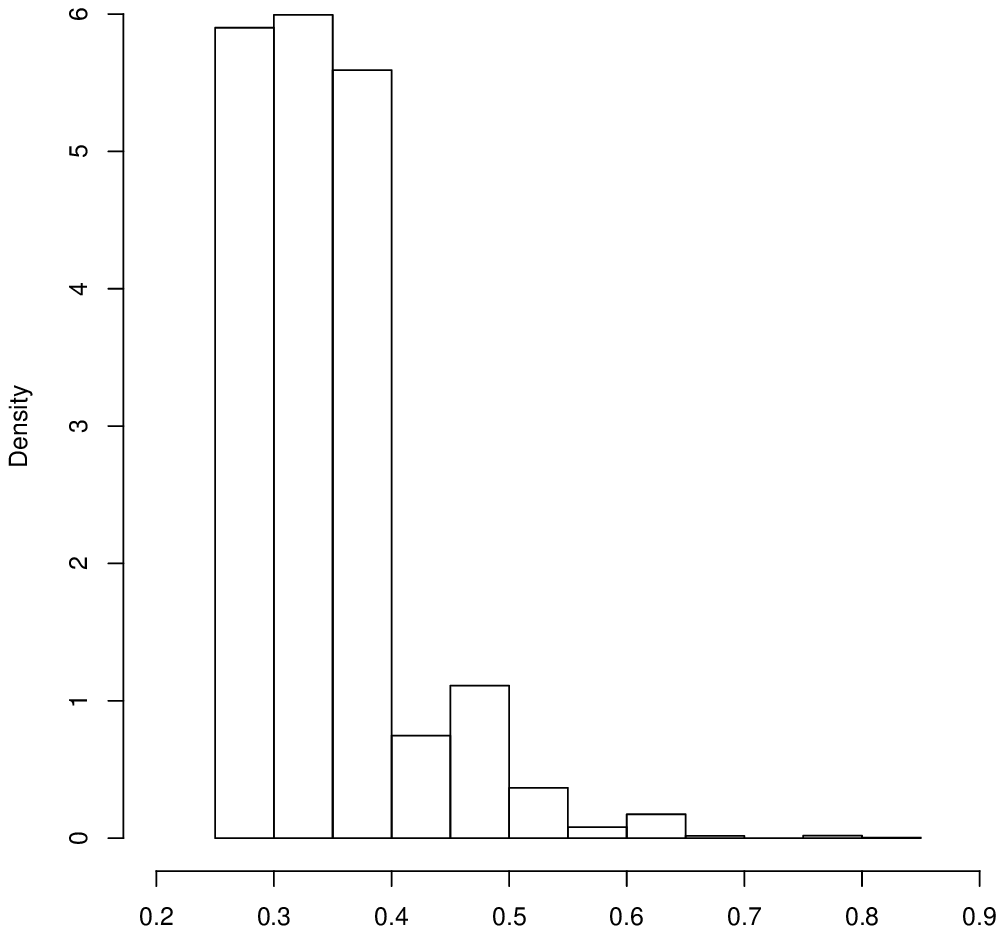} } }
\rotatebox{-90}{ \resizebox{2.0 in}{!}{ \includegraphics{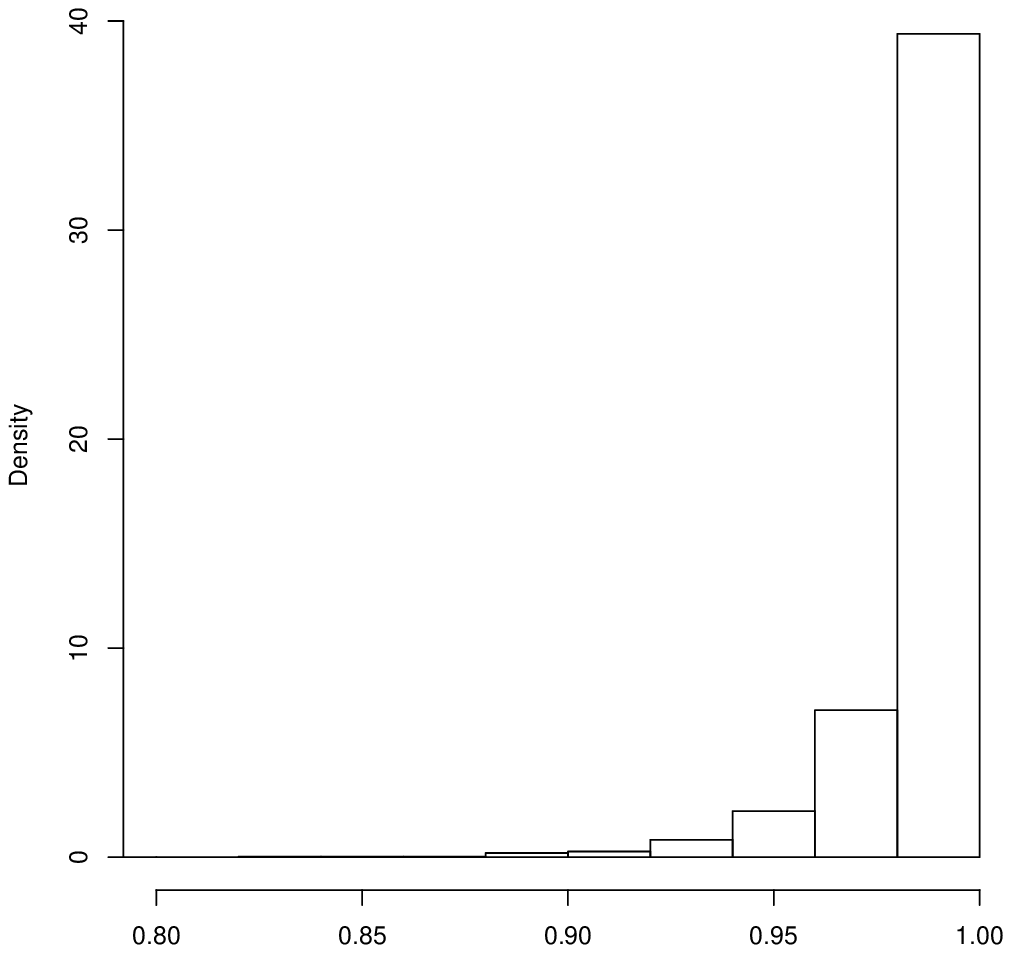} } }
\caption{
\label{fig:ORskew}
Depicted are the histograms for 10000 Monte Carlo replicates of $\rho^{\lo}_{10}(1)$ (left)
and $\rho^{\lo}_{10}(5)$ (right) indicating severe small sample skewness for extreme values of $r$.
Notice that the vertical axes are differently scaled.
}
\end{figure}

\begin{figure}
\centering
\psfrag{r}{\Huge{$r$}}
\rotatebox{-90}{ \resizebox{1.5 in}{!}{ \includegraphics{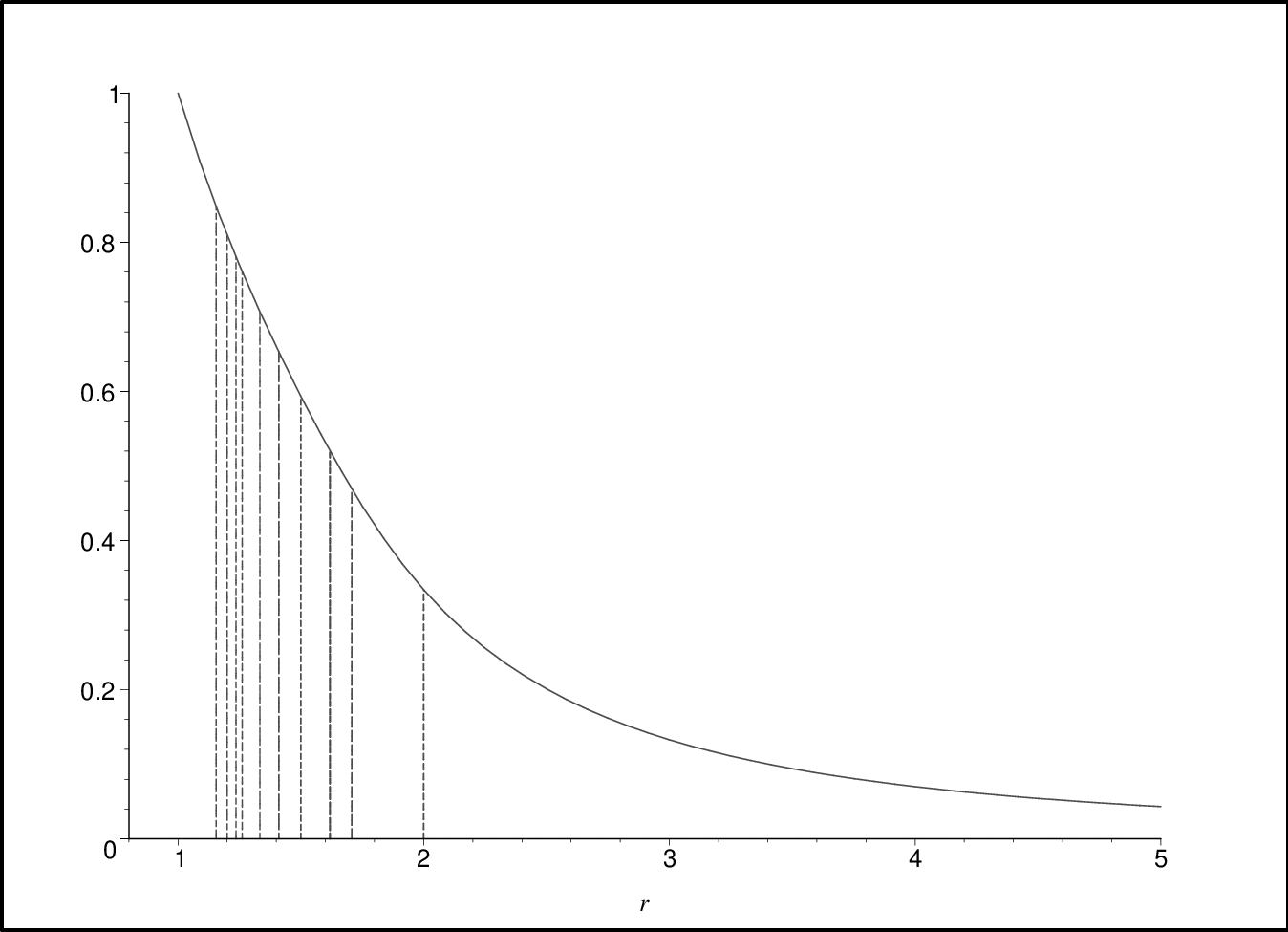} } }
\rotatebox{-90}{ \resizebox{1.5 in}{!}{ \includegraphics{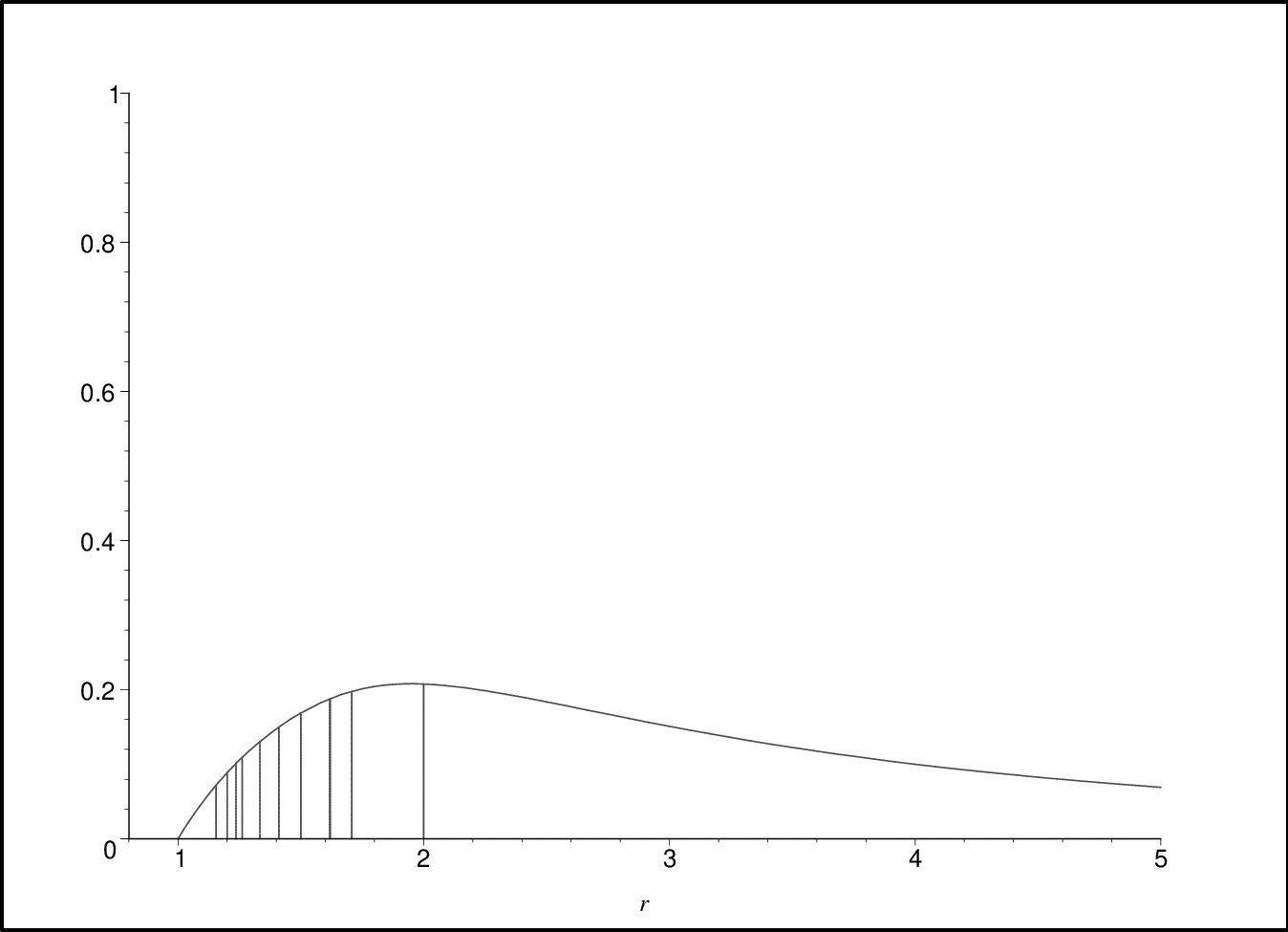} } }
\rotatebox{-90}{ \resizebox{1.5 in}{!}{ \includegraphics{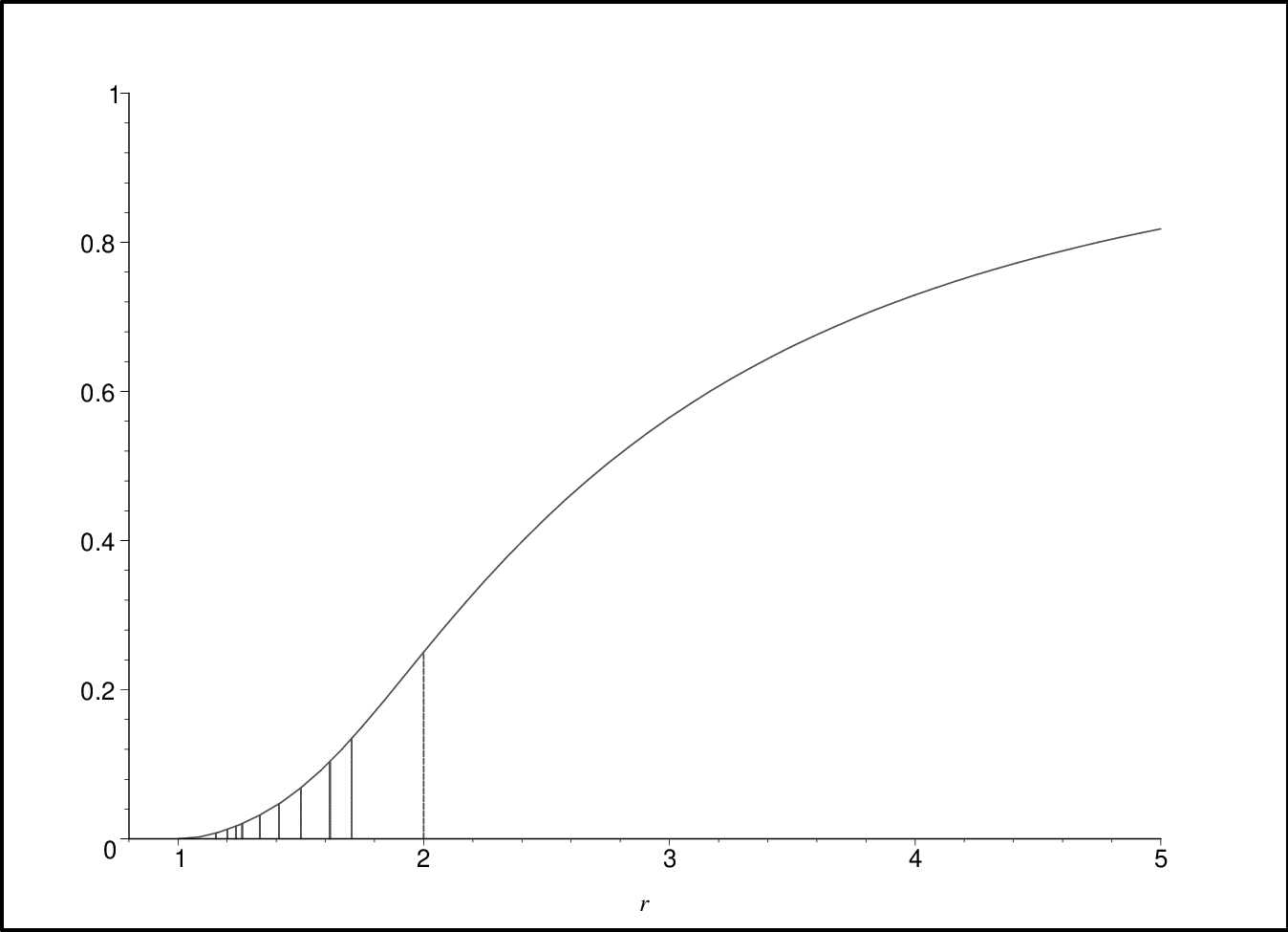} } }
\caption{
\label{fig:joint-dist-hand1213}
The plots for the joint distribution of $h^{\la}_{12}(r),h^{\la}_{13}(r)$ for $r \in [1,5]$.
Plotted are $P(h^{\la}_{12}(r),h^{\la}_{13}(r)=(0,0))$ (left),
$P(h^{\la}_{12}(r),h^{\la}_{13}(r)=(1,0))=P(h^{\la}_{12}(r),h^{\la}_{13}(r)=(0,1))$ (middle),
and
$P(h^{\la}_{12}(r),h^{\la}_{13}(r)=(1,1))$ (right).}
\end{figure}

\begin{figure}
\centering
\psfrag{r}{\Huge{$r$}}
\rotatebox{-90}{ \resizebox{1.5 in}{!}{ \includegraphics{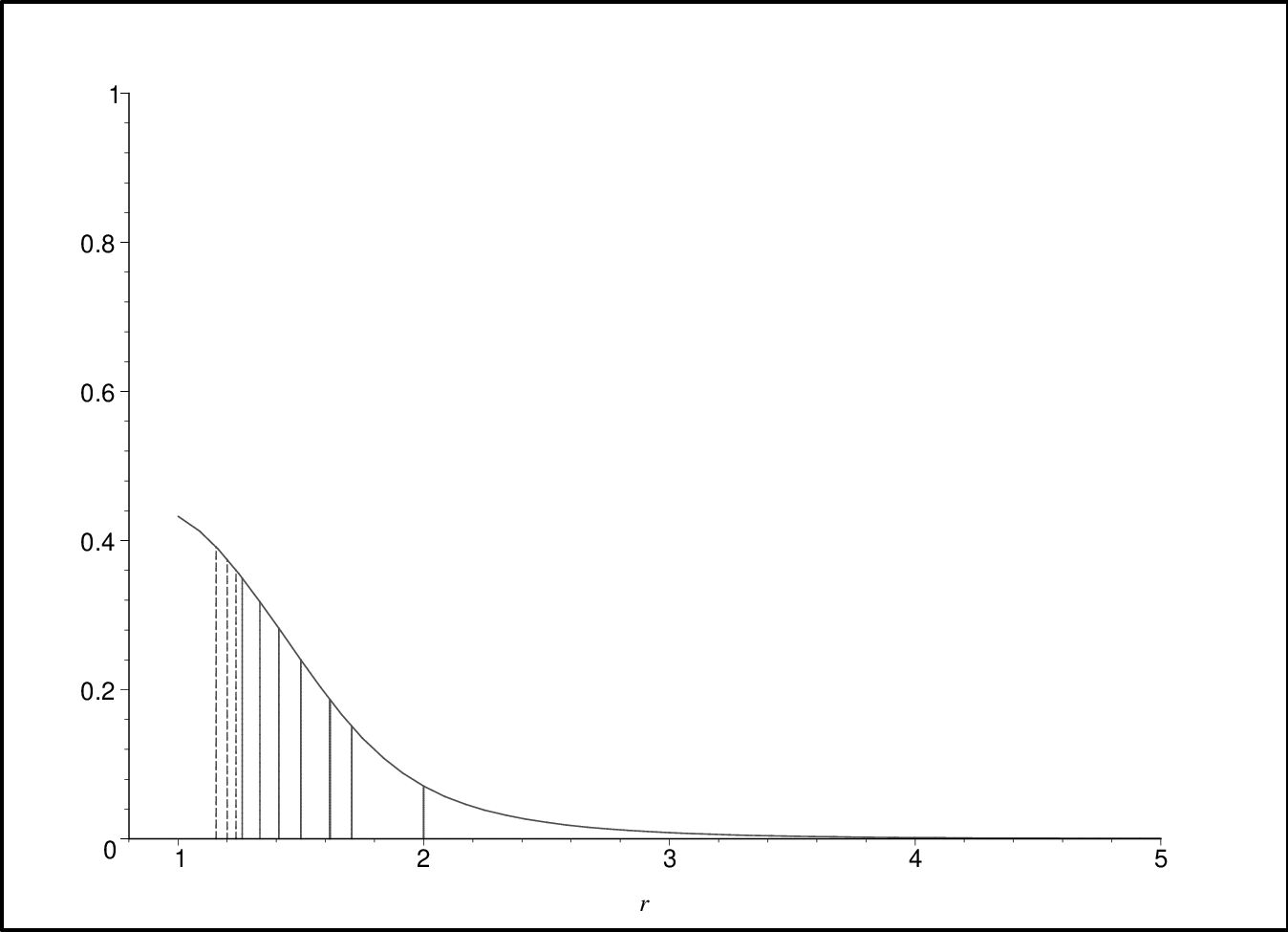} } }
\rotatebox{-90}{ \resizebox{1.5 in}{!}{ \includegraphics{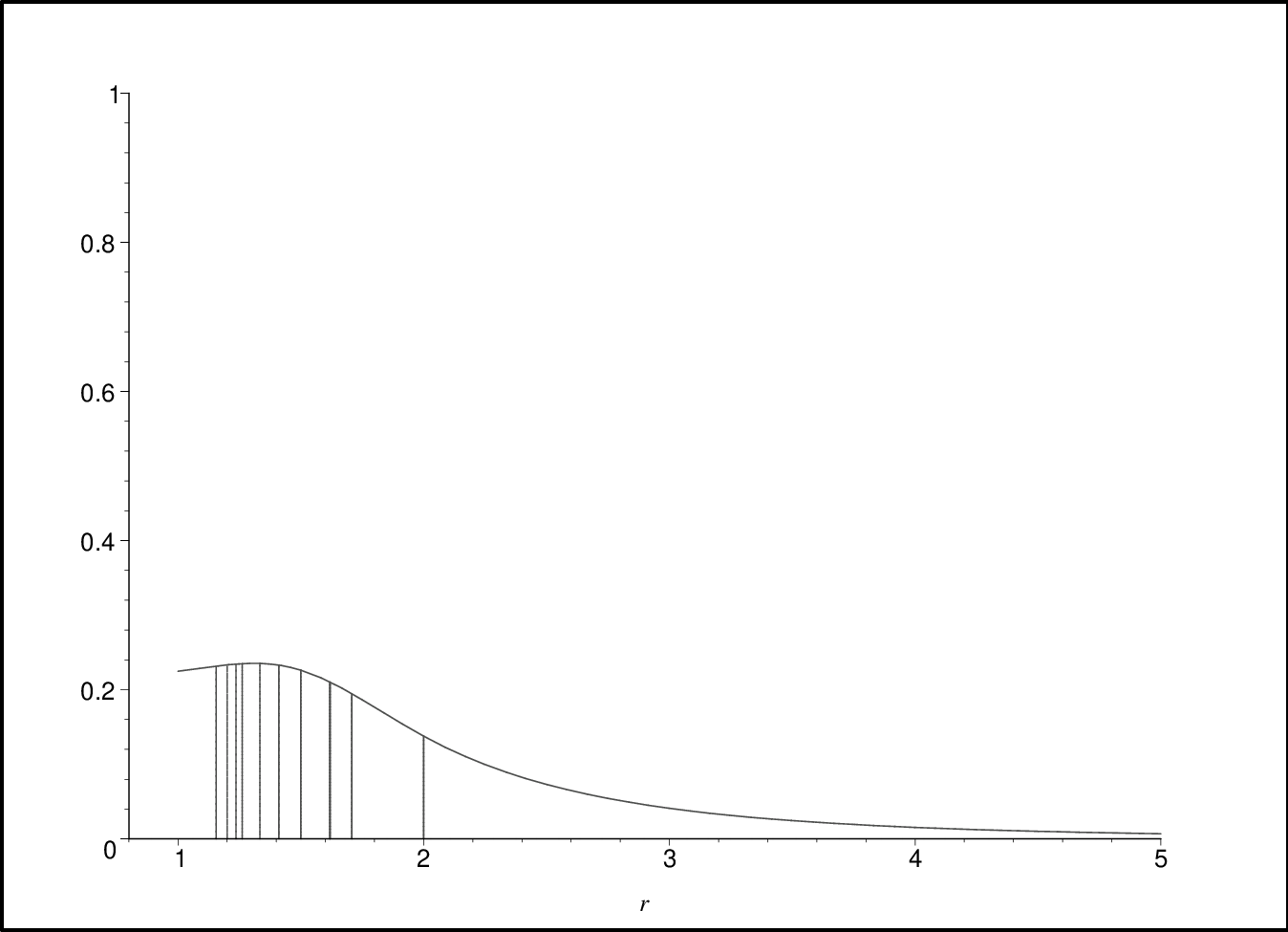} } }
\rotatebox{-90}{ \resizebox{1.5 in}{!}{ \includegraphics{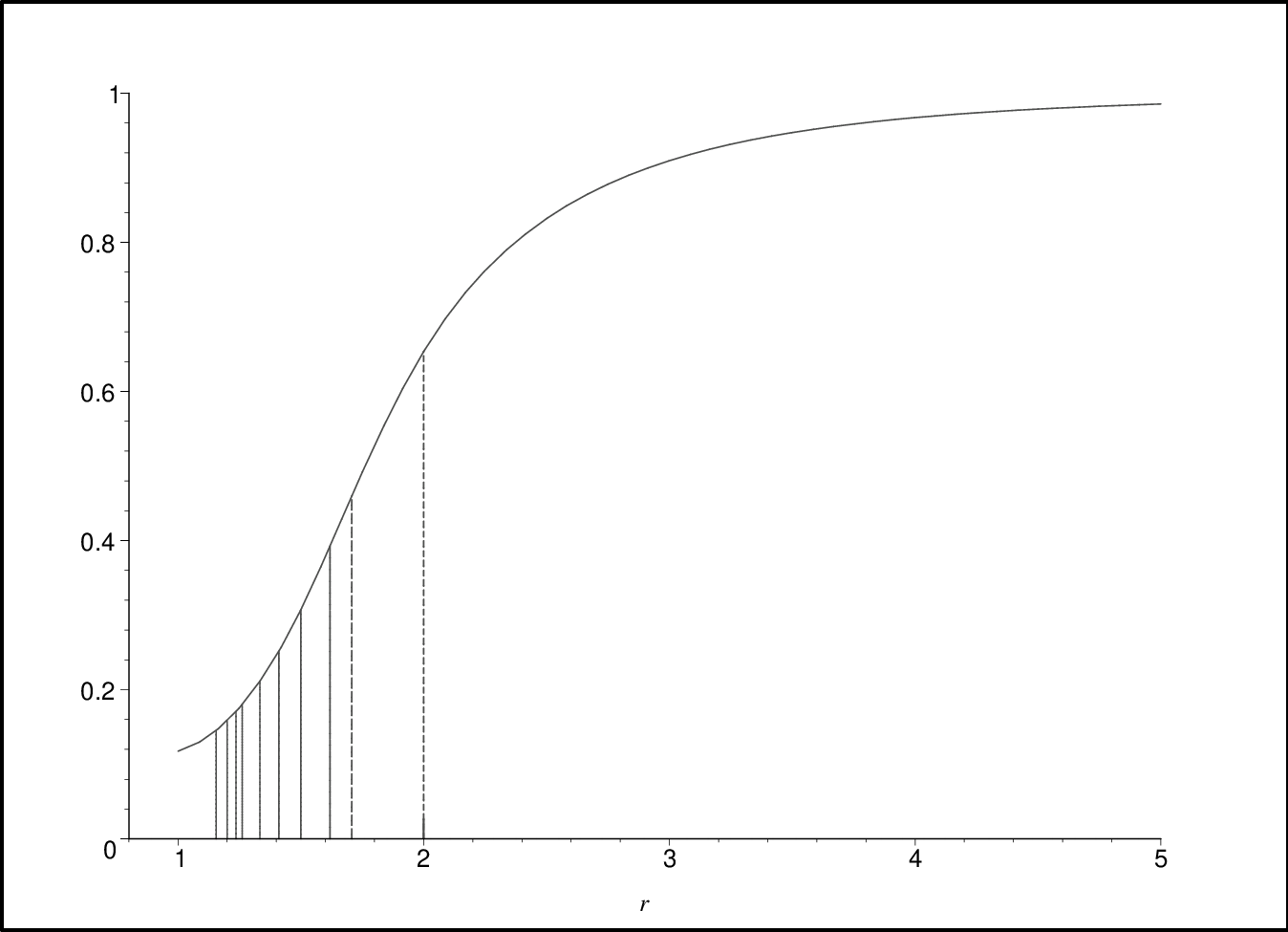} } }
\caption{
\label{fig:joint-dist-hor1213}
The plots for the joint distribution of $h^{\lo}_{12}(r),h^{\lo}_{13}(r)$ for $r \in [1,5]$.
Plotted are $P(h^{\lo}_{12}(r),h^{\lo}_{13}(r)=(0,0))$ (left),
$P(h^{\lo}_{12}(r),h^{\lo}_{13}(r)=(1,0))=P(h^{\lo}_{12}(r),h^{\lo}_{13}(r)=(0,1))$ (middle),
and
$P(h^{\lo}_{12}(r),h^{\lo}_{13}(r)=(1,1))$ (right).}
\end{figure}

We assess the accuracy of the asymptotic normality for finite sample data based on Monte Carlo simulations.
We generate $n$ $\X$ points independently uniformly in the standard equilateral triangle $T_e$.
For each data set generated,
we calculate the relative edge density values for the reflexivity and underlying
graphs based on the proportional-edge PCD with $r=2$.
We replicate the above process $N_{mc}=1000$ times for each of $n=10,20$, and 100.
We plot the histograms of the relative edge densities of the reflexivity and underlying graphs
using the simulated data and the corresponding (asymptotic) normal
curves in Figures \ref{fig:NormApprox_1} and \ref{fig:NormApprox2}, respectively.
Notice that, for $r=2$,
the normal approximation is accurate even for small $n$
although kurtosis may be indicated for $n=10$ in the reflexivity graph case,
and skewness may be indicated for $n=10$ in the underlying graph case.
We also investigate the behavior of the relative edge densities for extreme values of $n$ and $r$.
So we generate $n=10$ $\X$ points and calculate the relative edge densities for $r=1.05$ and $r=5$.
We repeat the above procedure $N_{mc}=10000$ times and plot
the histograms of the relative edge densities in Figures \ref{fig:ANDskew} and \ref{fig:ORskew},
which demonstrate that severe skewness is obtained for these extreme values of $n$ and $r$.
The finite sample variance and skewness may be derived analytically
in much the same way as was $4\,\nu_{\la}(r)$
(and $4\,\nu_{\lo}(r)]$)
for the asymptotic variance.
In fact,
the exact distribution of $\rho^{\la}_n(r)$ (and $\rho^{\lo}_n(r)$)
is, in principle, available
by successively conditioning on the values of the $X_i$.
Alas,
while the joint distribution of $h^{\la}_{12}(r),h^{\la}_{13}(r)$
(and $h^{\lo}_{12}(r),h^{\lo}_{13}(r)$) is available (see Figures \ref{fig:joint-dist-hand1213} and \ref{fig:joint-dist-hor1213}),
the joint distribution of $\{h^{\la}_{ij}(r)\}_{1 \leq i < j \leq n}$ (and $\{h^{\lo}_{ij}(r)\}_{1 \leq i < j \leq n}$),
and hence the calculation for the exact distribution of $\rho^{\la}_n(r)$ (and $\rho^{\lo}_n(r)$),
is extraordinarily tedious and lengthy for even small values of $n$.

\vspace*{0.05 in}

Let $\g_n(r)$ be the domination number of the proportional-edge PCD
based on $\X_n$ which is a random sample from $\U(\TY)$.
Additionally, let $\g^\la_n(r)$ and $\g^\lo_n(r)$ be
the domination numbers of the reflexivity and underlying graphs
based on the proportional-edge PCD, respectively.
Then we have the following stochastic ordering for the domination numbers.
\begin{theorem}
For all $r \in [1,\infty)$ and $n>1$,
$\g^\lo_n(r) <^{ST} \g_n(r) <^{ST} \g^\la_n(r)$.
\end{theorem}

\noindent \textbf{Proof:}
For all $x \in \TY $,
we have
$\left[ \NPE^r(x)\cap \G^r_1(x)\right] \subseteq \NPE^r(x) \subseteq \left[\NPE^r(x) \cup \G^r_1(x) \right]$.
For $X \sim \U(\TY)$,
we have
$\left[ \NPE^r(X)\cap \G^r_1\left( X \right)\right] \subsetneq \NPE^r(X) \subsetneq \left[\NPE^r(X) \cup \G^r_1\left( X \right)\right]$ a.s.
Moreover,
$\g_n(r)=1$ iff $\X_n \subset \NPE^r(X_i)$ for some $i$;
$\g^\la_n(r)=1$ iff $\X_n \subset \NPE^r(X_i)\cap \G^r_1\left( X_i \right)$ for some $i$;
and
$\g^\lo_n(r)=1$ iff $\X_n \subset \NPE^r(X_i)\cup \G^r_1\left( X_i \right)$ for some $i$.
So it follows that
$P(\g^\la_n(r)=1) < P(\g_n(r)=1) < P(\g^\lo_n(r)=1)$.
Similarly, for all $x, y \in \TY $,
we have
$\Bigl([\NPE^r(x)\cap \G^r_1(x)] \cup [\NPE^r(y)\cap \G^r_1(y)] \Bigr)  \subseteq
\left(\NPE^r(x) \cup \NPE^r(y) \right)\subseteq
\Bigl([\NPE^r(x) \cup \G^r_1(x)] \cup [\NPE^r(y) \cup \G^r_1(y)]\Bigr)$.
For $X,Y \stackrel{iid}{\sim} \U(\TY)$,
we have
$\Bigl( [\NPE^r(X)\cap \G^r_1(X)] \cup [\NPE^r(Y)\cap \G^r_1(Y)] \Bigr)
\subsetneq
\left( \NPE^r(X) \cup \NPE^r(Y) \right)
\subsetneq
\Bigl([\NPE^r(X) \cup \G^r_1(X)] \cup [\NPE^r(Y) \cup \G^r_1(Y)] \Bigr)$ a.s.
Moreover,
$\g_n(r) \le 2$ iff $\X_n \subset \NPE^r(X_i) \cup \NPE^r(X_j)$ for some $i \not= j$;
$\g^\la_n(r) \le 2$ iff $\X_n \subset \Bigl( [\NPE^r(X_i)\cap \G^r_1\left( X_i \right)] \cup
[\NPE^r(X_j)\cap \G^r_1\left( X_j \right)]\Bigr)$ for some $i \not= j$;
and
$\g^\lo_n(r)\le 2$ iff $\X_n \subset \Bigl( [\NPE^r(X_i)\cup \G^r_1\left( X_i \right)] \cup
[\NPE^r(X_j)\cup \G^r_1\left( X_j \right)] \Bigr)$ for some $i \not= j$.
So it follows that
$P(\g^\la_n(r) \le 2) < P(\g_n(r) \le 2) < P(\g^\lo_n(r) \le 2)$.
Since $P(\g_n(r)\le 3)=1$ (\cite{ceyhan:dom-num-NPE-SPL}),
it follows that $P(\g^\lo_n(r) \le 3)=1$ also holds since
$P(\g_n(r)\le 3) \le P(\g^\lo_n(r) \le 3)$.
Hence the desired stochastic ordering follows.
$\blacksquare$

Note the stochastic ordering in the above theorem
holds for any continuous distribution $F$ with support
being in $\TY$.
For $r=\infty$, we have
$\g^\lo_n(r) = \g_n(r) = \g^\la_n(r)=1$ a.s.

\section{Multiple Triangle Case}
\label{sec:multiple-triangle-case}
Suppose $\Y_m$ is a finite set of $m>3$ points in $\mathbb{R}^2$.
Consider the Delaunay triangulation (assumed to exist) of $\Y_m$.
Let $T_i$ denote the $i^{th}$ Delaunay triangle,
$J_m$ denote the number of triangles, and
$C_H(\Y_m)$ denote the convex hull of $\Y_m$.
For $X_i \stackrel{iid}{\sim} \mathcal{U}(C_H(\Y_m))$, $i=1,2,\ldots,n$,
we construct the proportional-edge PCD, $D_{n,m}(r)$,
using $\NPE^{r}(\cdot)$ as described in Section \ref{sec:prop-edge},
where for $X_i \in T_j$,
the three points in $\Y_m$ defining the
Delaunay triangle $T_j$ are used as $\Y_{[j]}$.
We investigate the relative edge densities of the underlying and reflexivity graphs
based on the proportional-edge PCD.
We consider various versions of the relative edge density
in the multiple triangle case.

\subsection{First Version of Relative Edge Density in the Multiple Triangle Case}
\label{sec:version-I-mult-tri}
For $J_m>1$,
as in Section \ref{sec:relative-density-PEPCD},
let $\rho^{\la}_{I,n}(r)=2\,\left|\mE_{\la}\right|/(n\,(n-1))$
and $\rho^{\lo}_n(r)=2\,\left|\mE_{\lo}\right|/(n\,(n-1))$.
Let $\mE^\la_{[i]}$ be the number of edges and
$\rho^{\la}_{{}_{[i]}}(r)$ be the relative edge density for triangle $i$ in the reflexivity graph case,
and $\mE^\lo_{[i]}$ and $\rho^{\lo}_{{}_{[i]}}(r)$ be similarly defined for underlying graph case.
Let $n_i$ be the number of $X$ points in $T_i$ for $i=1,2,\ldots,J_m$.
Letting $w_i = A(T_i) / A(C_H(\Y_m))$ with $A(\cdot)$ being the area functional,
we obtain the following as a corollary to Theorem \ref{thm:asy-norm-under}.

\begin{corollary}
\label{cor:MT-asy-norm-NYr}
For $r \in (1,\infty)$,
the asymptotic distribution for $\rho^{\la}_{I,n}(r)$ conditional on $\Y_m$
is given by
\begin{equation}
\sqrt{n}\left(\rho^{\la}_{I,n}(r)-\widetilde p_{\la}(m,r)\right)
\stackrel{\mathcal L}{\longrightarrow}\\
\mathcal{N}
 \left(
   0,
   4\,\widetilde \nu_{\la}(m,r)
 \right),
\end{equation}
as $n \rightarrow \infty$,
where $\widetilde p_{\la}(m,r)=p_{\la}(r) \left(\sum_{i=1}^{J_m}w_i^2\right)$
and
$$\widetilde \nu_{\la}(m,r)=
\left[  \nu_{\la}(r) \left(\sum_{i=1}^{J_m}w_i^3 \right)+
\left( p_{\la}(r) \right)^2\left(\sum_{i=1}^{J_m}w_i^3-\left(\sum_{j=1}^{J_m}w_i^2 \right)^2\right) \right]$$
with $p_{\la}(r)$ and $\nu_{\la}(r)$ being as in Equations \eqref{eqn:Asymean_and} and \eqref{eqn:Asyvar_and},
respectively.
The asymptotic distribution of $\rho^{\lo}_{I,n}(r)$ with $r \in [1,\infty)$ is similar.
\end{corollary}

The proof is provided in Appendix 4.
By an appropriate application of the Jensen's inequality,
we see that $\sum_{i=1}^{J_m}w_i^3 \ge \left(\sum_{i=1}^{J_m}w_i^2 \right)^2.$
So the covariance above is zero iff $\nu_{\la}(r)=0$ and
$\sum_{i=1}^{J_m}w_i^3=\left(\sum_{i=1}^{J_m}w_i^2 \right)^2$, so
asymptotic normality may hold even though $\nu_{\la}(r)=0$.
That is, $\rho^{\la}_{I,n}(r)$ has the asymptotic normality for $r \in \{1,\infty\}$ also
provided that $\sum_{i=1}^{J_m}w_i^3 > \left(\sum_{i=1}^{J_m}w_i^2 \right)^2$.
The same holds for the underlying graph case (for $r=\infty$).


\subsection{Other Versions of Relative Edge Density in the Multiple Triangle Case}
\label{sec:version-II-mult-tri}

Let $\displaystyle \Xi^{\la}_n(r):=\sum_{i=1}^{J_m}\frac{n_i\,(n_i-1)}{n\,(n-1)} \rho^{\la}_{{}_{[i]}}(r)$.
Then $\Xi^{\la}_n(r) = \rho^{\la}_{I,n}(r)$,
since $\displaystyle \Xi^{\la}_n(r)=\sum_{i=1}^{J_m}\frac{n_i\,(n_i-1)}{n\,(n-1)} \rho^{\la}_{{}_{[i]}}(r)=
\frac{\sum_{i=1}^{J_m} 2\,|\mE^\la_{[i]}|}{n\,(n-1)}=\frac{2\,\left|\mE_{\la}\right|}{n\,(n-1)}=\rho^{\la}_{I,n}(r)$.
Similarly, we have $\Xi^{\lo}_n(r) = \rho^{\lo}_n(r)$.

Furthermore,
let $\widehat{\Xi}^{\la}_n:= \sum_{i=1}^{J_m}w_i^2\,\rho^{\la}_{{}_{[i]}}(r)$
where $w_i$ is as in Section \ref{sec:version-I-mult-tri}.
So $\widehat{\Xi}^{\la}_n$ a mixture of the $\rho^{\la}_{{}_{[i]}}(r)$ values.
Since the $\rho^{\la}_{{}_{[i]}}(r)$ are asymptotically independent,
$\Xi^{\la}_n(r), \, \rho^{\la}_{I,n}(r)$ are asymptotically normal;
i.e., for large $n$ their distribution is approximately
$\mathcal N\left( \widetilde p_{\la}(m,r),4\,\widetilde \nu_{\la}(m,r)/n\right)$.
A similar result holds for the underlying graph case.

In Section \ref{sec:version-I-mult-tri},
the denominator of $\rho^{\la}_{I,n}(r)$ has
$n(n-1)/2$ as the maximum number of edges possible.
However, by definition,
given the $n_i$ values,
we can have a graph with at most $J_m$ complete components,
each with order $n_i$ for $i=1,2,\ldots,J_m$.
Then the maximum number of edges possible is $n_t:=\sum_{i=1}^{J_m}n_i\,(n_i-1)/2$
which suggests another version of edge density, namely,
$\displaystyle \rho^{\la}_{II,n}(r):=\frac{\left|E_{\la}\right|}{n_t}$.
Then $\displaystyle \rho^{\la}_{II,n}(r)=\frac{\sum_{i=1}^{J_m} |E^\la_{[i]}|}{n_t}=
\sum_{i=1}^{J_m}\frac{n_i\,(n_i-1)}{2\,n_t}\,\rho^{\la}_{{}_{[i]}}(r)$.
Since $\frac{n_i\,(n_i-1)}{2\,n_t} \ge 0$ for each $i$,
and $\displaystyle \sum_{i=1}^{J_m}\frac{n_i\,(n_i-1)}{2\,n_t}=1$,
$\rho^{\la}_{II,n}(r)$ is a mixture of the $\rho^{\la}_{{}_{[i]}}(r)$.
Then $\E\left[ \rho^{\la}_{II,n}(r) \right]=p_{\la}(r)$.
A similar result holds for the underlying graph case also.

\begin{theorem}
\label{thm:MT-asy-norm-II}
The asymptotic distribution for $\rho^{\la}_{II,n}(r)$ conditional on $\Y_m$
for $r \in (1,\infty)$
is given by
\begin{equation}
\sqrt{n}\left(\rho^{\la}_{II,n}(r)-p_{\la}(m,r)\right)
\stackrel{\mathcal L}{\longrightarrow}\\
\mathcal{N}
 \left(
   0,
   4\,\breve \nu_{\la}(m,r)
 \right),
\end{equation}
as $n \rightarrow \infty$,
where
$\breve \nu_{\la}(m,r)=
\left[  \nu_{\la}(r) \left(\sum_{i=1}^{J_m}w_i^3 \right)\Big/\left(\sum_{i=1}^{J_m}w_i^2 \right)^2 \right]$
with $p_{\la}(r)$ and $\nu_{\la}(r)$ being as in Equations \eqref{eqn:Asymean_and} and \eqref{eqn:Asyvar_and},
respectively.
The asymptotic distribution of $\rho^{\lo}_{II,n}(r)$ with $r \in [1,\infty)$ is similar.
\end{theorem}

The proof is provided in Appendix 5.
Notice that the covariance $\breve \nu_{\la}(m,r)$ is zero iff $\nu_{\la}(r)=0$.
The underlying graph case is similar.

\begin{remark}
\label{rem:comp-versions-mult-tri}
\textbf{Comparison of Versions of Relative Edge Density in the Multiple Triangle Case:}
Among the versions of the relative edge density we considered,
$\Xi^{\la}_n(r) = \rho^{\la}_{I,n}(r)$ for all $n>1$,
and
$\widehat{\Xi}^{\la}_n$ and $\rho^{\la}_{I,n}(r)$
are asymptotically equivalent (i.e., they have the
same asymptotic distribution).
However, $\rho^{\la}_{I,n}(r)$ and $\rho^{\la}_{II,n}(r)$
do not have the same distribution for finite or infinite $n$.
But we have $\rho^{\la}_{I,n}(r)=\frac{2\,n_t}{n(n-1)}\rho^{\la}_{II,n}(r)$
and
since $\sum_{i=1}^{J_m}w_i^2<1$,
it follows that $\widetilde p_{\la}(m,r) < \breve p_{\la}(m,r)=p_{\la}(r)$.
Furthermore,
since $\frac{2\,n_t}{n(n-1)}=\sum_{i=1}^{J_m} \frac{n_i(n_i-1)}{n(n-1)}
\longrightarrow \sum_{i=1}^{J_m} w_i^2$ as $n_i \rightarrow \infty$,
we have
$\lim_{n_i \rightarrow \infty} \Var[\sqrt{n}\rho^{\la}_{I,n}(r)]=
\left(\sum_{i=1}^{J_m} w_i^2\right)^2 \lim_{n_i \rightarrow \infty} \Var[\sqrt{n}\rho^{\la}_{I,n}(r)]$.
Hence $\breve \nu_{\la}(m,r) \ge \widetilde \nu_{\la}(m,r)$.
Therefore,
we recommend $\rho^{\la}_{I,n}(r)$
for use in spatial pattern analysis in the multiple triangle case.
Moreover, asymptotic normality might hold for $\rho^{\la}_{I,n}(r)$
even if $\nu_{\la}(r)=0$.
$\square$
\end{remark}

\subsection{Extension to Higher Dimensions}
\label{sec:extend-high-dim}
The extension to $\mathbb{R}^d$ for $d > 2$ is straightforward.
Let $\Y_{d+1} = \{\y_1,\y_2,\ldots,\y_{d+1}\}$ be $d+1$ non-coplanar points.
Denote the simplex formed by these $d+1$ points as $\mathfrak S (\Y_{d+1})$.
A simplex is the simplest polytope in $\mathbb{R}^d$
having $d+1$ vertices, $d\,(d+1)/2$ edges and $d+1$ faces of dimension $(d-1)$.
For $r \in [1,\infty]$, define the proportional-edge proximity map as follows.
Given a point $x$ in $\mathfrak S (\Y_{d+1})$,
let $\y := \arg\min_{\y \in \Y_{d+1}} \mbox{volume}(Q_\y(x))$
where $Q_\y(x)$ is the polytope with vertices being the $d\,(d+1)/2$ midpoints of the edges,
the vertex $\y$ and $x$.
That is, the vertex region for vertex $v$ is the polytope with vertices
given by $v$ and the midpoints of the edges.
Let $v(x)$ be the vertex in whose region $x$ falls.
If $x$ falls on the boundary of two vertex regions or at the center of mass, we assign $v(x)$ arbitrarily.
Let $\varphi(x)$ be the face opposite to vertex $v(x)$,
and $\eta(v(x),x)$ be the hyperplane parallel to $\varphi(x)$ which contains $x$.
Let $d(v(x),\eta(v(x),x))$ be the Euclidean distance from $v(x)$ to $\eta(v(x),x)$.
For $r \in [1,\infty)$, let $\eta_r(v(x),x)$ be the hyperplane parallel to $\varphi(x)$
such that
$$d(v(x),\eta_r(v(x),x))=r\,d(v(x),\eta(v(x),x)) \text{ and } d(\eta(v(x),x),\eta_r(v(x),x))< d(v(x),\eta_r(v(x),x)).$$
Let $\mathfrak S_r(x)$ be the polytope similar to and with the same orientation as $\mathfrak S$
having $v(x)$ as a vertex and $\eta_r(v(x),x)$ as the opposite face.
Then the proportional-edge proximity region $\NPE^r(x):=\mathfrak S_r(x) \cap \mathfrak S(\Y_{d+1})$.
Furthermore, let $\zeta_i(x)$ be the hyperplane such that
$\zeta_i(x) \cap \mathfrak S(\Y_{d+1}) \not=\emptyset$ and $r\,d(\y_i,\zeta_i(x))=d(\y_i,\eta(\y_i,x))$
for $i=1,2,\ldots,d+1$.
Then $\G_1^r(x)\cap R(\y_i)=\{z \in R(\y_i): d(\y_i,\eta(\y_i,z)) \ge d(\y_i,\zeta_i(x)\}$, for $i=1,2,3$.
Hence  $\G_1^r(x)=\bigcup_{j=1}^{d+1} (\G_1^r(x)\cap R(\y_i))$.  Notice that $r \ge 1$ implies $x \in \NPE^r(x)$ and $x \in \G_1^r(x)$.

Theorem \ref{thm:geo-inv-NYr-under} generalizes,
so that any simplex $\mathfrak S$ in $\mathbb{R}^d$
can be transformed into a regular polytope
(with edges being equal in length and faces being equal in volume)
preserving uniformity.
Delaunay triangulation becomes Delaunay tessellation in $\mathbb{R}^d$,
provided no more than $d+1$ points being cospherical
(lying on the boundary of the same sphere).
In particular, with $d=3$, the general simplex is a tetrahedron
(4 vertices, 4 triangular faces and 6 edges),
which can be mapped into a regular tetrahedron
(4 faces are equilateral triangles) with vertices
$(0,0,0)\,(1,0,0)\,(1/2,\sqrt{3}/2,0),\,(1/2,\sqrt{3}/4,\sqrt{3}/2)$.

Asymptotic normality of the $U$-statistic holds for $d>2$ in both underlying cases.

\section{Discussion and Conclusions}
\label{sec:discussion}
In this article,
we demonstrate that the relative edge density of random graphs
and relative arc density of random digraphs are one-sample $U$-statistics of degree 2.
Then,
we specify the conditions under which the asymptotic normality of the relative densities holds
for the random graphs and digraphs.
We consider the asymptotic distribution of the
relative edge density of the underlying and reflexivity graphs based on
(parameterized) proportional-edge proximity catch digraphs (PE-PCDs).
In particular, we consider the reflexivity and underlying graphs based on the proportional-edge PCD;
and derive the asymptotic distribution of the relative edge density
using the central limit theory of $U$-statistics.
We compute the asymptotic mean and variance of the limiting normal distribution
for uniform data
based on detailed geometric calculations.
Moreover, we compare the asymptotic distributions of the relative edge densities
of the underlying and reflexivity graphs and of the relative arc density of the PCDs.

The PCDs have applications in classification and spatial pattern analysis.
\cite{ceyhan:arc-density-PE} used that the relative (arc) density
of the PE-PCDs for testing bivariate spatial patterns.
The relative edge densities of the underlying and reflexivity graphs based on this PCD can be employed
for the same purpose.
More specifically,
the relative edge densities can be employed for testing the complete spatial randomness (CSR)
of two or more classes of points against the segregation or association
of the points from the classes.
\emph{CSR} is roughly defined as the lack of
spatial interaction between the points in a given study area.
In particular, the null hypothesis can be assumed to be CSR of $\X$ points,
i.e., the uniformness of $\X$ points in the convex hull of $\Y$ points.
\emph{Segregation} is the pattern in which points of one class tend to
cluster together, i.e., form one-class clumps.
On the other hand, \emph{association} is the pattern in which
the points of one class tend to occur more frequently around points from the other class.
Under the segregation alternative,
the $\X$ points will tend to be further away from $\Y$ points
and under the association alternative $\X$ points will tend to cluster around the $\Y$ points.
Such patterns can be detected by the test
statistics based on the relative edge densities,
since under segregation we expect them to be smaller,
and under association they tend to be larger.
The underlying and reflexivity graphs can also be used in pattern classification as outlined in \cite{priebe:2003b}.
Moreover, the methodology described here is also applicable to PCDs in higher dimensions.

\section*{Acknowledgments}
This research was supported by the European Commission under the Marie Curie International Outgoing Fellowship Programme
via Project \# 329370 titled PRinHDD.


\section*{APPENDIX}

\subsection*{Appendix 1: The Asymptotic Variance of Relative Edge Density for the Reflexivity Graph}
The variance of $h^{\la}_{12}(r)$,
denoted as $\Var_{\la}(r)=\Var\left[ h^{\la}_{12}(r) \right]$, is as follows:
$$\Var_{\la}(r)=
\varphi^{\la}_{1,1}(r)\I(r \in [1,4/3))+
\varphi^{\la}_{1,2}(r)\I(r \in [4/3,3/2))+
\varphi^{\la}_{1,3}(r)\I(r \in [3/2,2))+
\varphi^{\la}_{1,4}(r)\I(r \in [2,\infty))
$$
where
$\varphi^{\la}_{1,1}(r)=-{\frac{(5\,r^6-153\,r^5+393\,r^4-423\,r^3-54\,r^2+360\,r-128)(447\,r^4-261\,r^3+54\,r^2+5\,r^6-153\,r^5+360\,r-128)}{2916\,r^4(r+2)^2(r+1)^2}}$,\\
$\varphi^{\la}_{1,2}(r)=-{\frac{(101\,r^5-801\,r^4+1302\,r^3-732\,r^2-536\,r+672)(1518\,r^3-84\,r^2-104\,r+101\,r^5-801\,r^4+672)}{46656\,r^2(r+2)^2(r+1)^2}}$,\\
$\varphi^{\la}_{1,3}(r)=-{\frac{(r^8-13\,r^7+30\,r^6+148\,r^5-448\,r^4+264\,r^3+288\,r^2-368\,r+96)(22\,r^6+124\,r^5-464\,r^4+r^8-13\,r^7+264\,r^3+288\,r^2-368\,r+96)}{64\,r^8(r+2)^2(r+1)^2}}$,\\
$\varphi^{\la}_{1,4}(r)={\frac{(r^5+r^4-3\,r^3-3\,r^2+6\,r-2)(3\,r^3+3\,r^2-6\,r+2)}{r^8(r+1)^2}}$.
Note that $\Var_{\la}(r=1)=0$ and $\lim_{r \rightarrow \infty}\Var_{\la}(r)=0$ (at rate $O(r^{-2})$),
and $\argsup_{r \in [1,\infty)} \Var_{\la}(r) \approx 2.1126$ with $\sup \Var_{\la}(r)=0.25$.
\begin{figure}
\centering
\psfrag{r}{\normalsize{$r$}}
\epsfig{figure=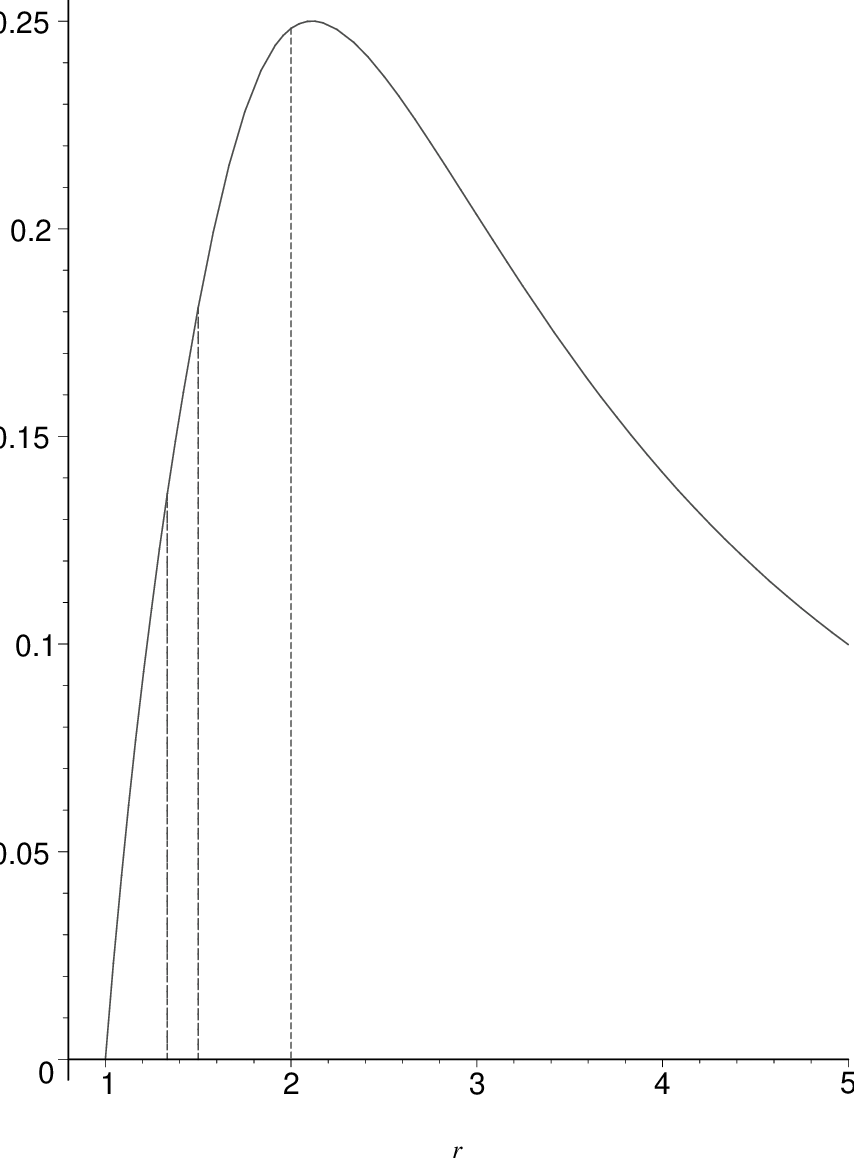, height=150pt, width=200pt}
\epsfig{figure=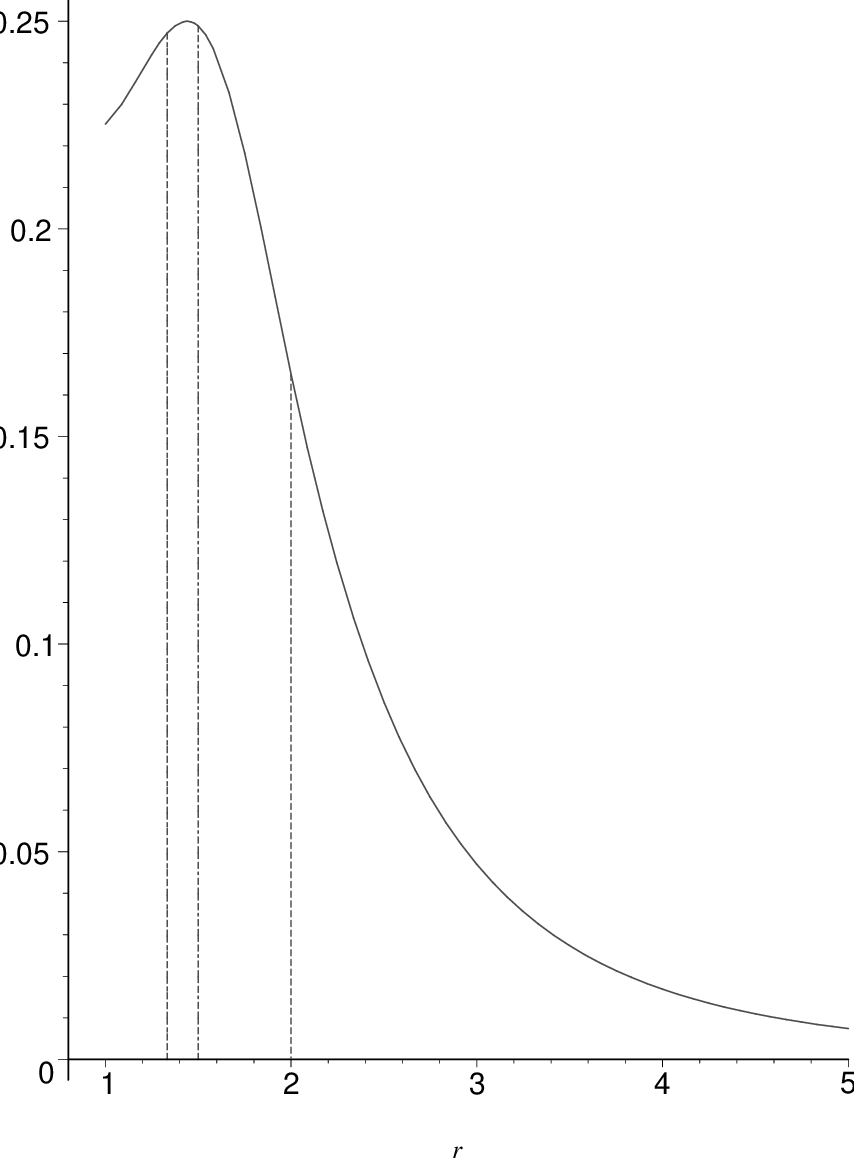, height=150pt, width=200pt}\\
\caption{
\label{var of nu under-graph}
$\Var\left[ h^{\la}_{12}(r) \right]$ (left)
and
$\Var\left[ h^{\lo}_{12}(r) \right]$ (right)
as a function of $r$ for $r \in [1,5] $.}
\end{figure}

The asymptotic variance for the reflexivity graph case is
$$\nu_{\la}(r):=\Cov\left[ h^{\la}_{12}(r),h^{\la}_{13}(r) \right]=\sum_{i=1}^{11}\vartheta^{\la}_i(r)\,\I(\mI_i)$$
where
{\small
\begin{multline*}
\vartheta^{\la}_1(r)=-\frac{1}{58320\,(2\,r^2+1)(r+2)^2(r+1)^3r^6}\, ((r-1)^2(972\,r^{19}+8748\,r^{18}+44456\,r^{17}+140328\,r^{16}+121371\,r^{15}\\
-412117\,r^{14}-27145\,r^{13}-4503501\,r^{12}+1336147\,r^{11}+10640999\,r^{10}-982009\,r^9-6677105\,r^8-2274458\,r^7\\
-1150162\,r^6+249126\,r^5+1232530\,r^4+1234372\,r^3+226776\,r^2-184944\,r-81920))
\end{multline*}
\begin{multline*}
\vartheta^{\la}_2(r)=-\frac{1}{116640\,(2\,r^2+1)(r+2)^2(r+1)^3r^6}\,(486\,r^{21}+3402\,r^{20}-269\,r^{19}-45155\,r^{18}-118850\,r^{17}+443518\,r^{16}\\
+3251855\,r^{15}-13836295\,r^{14}+13434672\,r^{13}+11140788\,r^{12}-27667544\,r^{11}+13293088\,r^{10}+7159710\,r^9-\\
13013598\,r^8 +4185440\,r^7+3262952\,r^6+586636\,r^5-1616444\,r^4-680120\,r^3-55952\,r^2+219936\,r+49152)
\end{multline*}
\begin{multline*}
\vartheta^{\la}_3(r)=-\frac{1}{116640\,(2\,r^2+1)(r+2)^2(r+1)^3r^6}\,(486\,r^{21}+3402\,r^{20}-269\,r^{19}-45155\,r^{18}-118850\,r^{17}+443518\,r^{16}\\
+2751855\,r^{15}-13736295\,r^{14}+18084672\,r^{13}+8770788\,r^{12}-43009544\,r^{11}+24604048\,r^{10}+27137438\,r^9-30889822\,r^8\\
-2832544\,r^7+11101160\,r^6-4168820\,r^5+2364868\,r^4+2305864\,r^3-3041936\,r^2+219936\,r+49152)
\end{multline*}
\begin{multline*}
\vartheta^{\la}_4(r)=-\frac{1}{58320\,(r+2)^3(r^2-2)(2\,r^2+1)(r+1)^3r^6}\,(3632\,r^{22}+25632\,r^{21}-60328\,r^{20}-441888\,r^{19}+1353430\,r^{18}\\
-297666\,r^{17}-4791125\,r^{16}+12849927\,r^{15}-10894618\,r^{14}-26295324\,r^{13}+62283823\,r^{12}-2280753\,r^{11}-81700012\,r^{10}\\
+32551926\,r^9+39974410\,r^8-11284026\,r^7-5806580\,r^6-9167580\,r^5-2004944\,r^4+4646688\,r^3+1931776\,r^2-489024\,r-98304)
\end{multline*}
\begin{multline*}
\vartheta^{\la}_5(r)=\vartheta^{\la}_6(r)=-\frac{1}{58320\,(r+2)^3(2\,r^2+1)(r^2+1)(r+1)^3r^6}\,(3632\,r^{22}+25632\,r^{21}-49432\,r^{20}-364992\,r^{19}+958940\,r^{18}\\
-1167012\,r^{17}+1200518\,r^{16}+5424126\,r^{15}-23566328\,r^{14}+23837088\,r^{13}+11797395\,r^{12}-41623065\,r^{11}+39261953\,r^{10}\\
-8239197\,r^9-30178496\,r^8+27901506\,r^7-4936170\,r^6+61038\,r^5+4719720\,r^4-5513952\,r^3+340736\,r^2+23328\,r+65536)
\end{multline*}
\begin{multline*}
\vartheta^{\la}_7(r)=\frac{1}{466560\,(r+2)^3(2\,r^2+1)(r^2+1)(r+1)^3r^5}\,(1562\,r^{21}-11142\,r^{20}-103099\,r^{19}+2105697\,r^{18}-9774118\,r^{17}+\\
10220280\,r^{16}+27825711\,r^{15}-69243129\,r^{14}+81624200\,r^{13}-76052574\,r^{12}-65530400\,r^{11}+262451196\,r^{10}-178092280\,r^9\\
-69106464\,r^8+158439568\,r^7-97568688\,r^6+12246288\,r^5+17591952\,r^4-21111616\,r^3+15628032\,r^2-2545664\,r+993024)
\end{multline*}
\begin{multline*}
\vartheta^{\la}_8(r)=-\frac{1}{1920\,(r+2)^3(r^2+1)(2\,r^2+1)(r+1)^3r^{10}}\,(2\,r^{26}-30\,r^{25}-2395\,r^{23}+281\,r^{24}+8770\,r^{22}
+29528\,r^{21}-268053\,r^{20}+\\
245667\,r^{19}+2066216\,r^{18}-5313494\,r^{17}-1589216\,r^{16}+18512684\,r^{15}-18946136\,r^{14}-2665248\,r^{13}+22789584\,r^{12}-\\
32987760\,r^{11}+20482512\,r^{10}+13109584\,r^9-28084416\,r^8+17326976\,r^7-3864576\,r^6-4579328\,r^5+6666240\,r^4-3576320\,r^3\\
+635904\,r^2-116736\,r+61440)
\end{multline*}
\begin{multline*}
\vartheta^{\la}_9(r)=-\frac{1}{1920\,(r+2)^3(r^2+1)(2\,r^2+1)(r+1)^3r^{10}}\,(2\,r^{26}-30\,r^{25}-2395\,r^{23}281\,r^{24}+8258\,r^{22}
+31064\,r^{21}-262677\,r^{20}+\\
225443\,r^{19}+2052136\,r^{18}-5219030\,r^{17}-1608928\,r^{16}+18337836\,r^{15}-18837080\,r^{14}-2598688\,r^{13}+22736336\,r^{12}-\\
32858736\,r^{11}+20384720\,r^{10}+12930896\,r^9-27988416\,r^8+17416832\,r^7-3862784\,r^6-4575488\,r^5+6638848\,r^4-3603200\,r^3\\
+640512\,r^2-107520\,r+63488)
\end{multline*}
\begin{multline*}
\vartheta^{\la}_{10}(r)=-\frac{1}{1920\,(r+2)^3(r-1)(r+1)^3(2\,r^2-1)r^{10}}\,(2\,r^{25}+307\,r^{23}-32\,r^{24}-2612\,r^{22}
+11572\,r^{21}+21934\,r^{20}-328867\,r^{19}+\\
524994\,r^{18}+2446870\,r^{17}-8676180\,r^{16}-437020\,r^{15}+36944680\,r^{14}-40677696\,r^{13}-44860384\,r^{12}+106256352\,r^{11}-\\
15515040\,r^{10}-98636848\,r^9+66358080\,r^8+27142272\,r^7-42614272\,r^6+7781120\,r^5+7327232\,r^4-3388672\,r^3+\\
430592\,r^2-171008\,r+63488)
\end{multline*}
\begin{multline*}
\vartheta^{\la}_{11}(r)=\frac{1}{15\,(2\,r^2-1)(r+1)^3r^{10}}\,(30\,r^{13}+90\,r^{12}-127\,r^{11}-621\,r^{10}+320\,r^9+1568\,r^8-858\,r^7
-1370\,r^6+909\,r^5+\\
295\,r^4-292\,r^3+44\,r^2+6\,r-2)
\end{multline*}
}
\noindent
and
$\mI_1=[1,2/\sqrt{3}),\;
\mI_2=[2/\sqrt{3},6/5),\;
\mI_3=[6/5,\sqrt{5}-1),\;
\mI_4=[\sqrt{5}-1,(6+2\,\sqrt{2})/7),\;
\mI_5=[(6+2\,\sqrt{2})/7,4/3),\;
\mI_6=[4/3,(6+\sqrt{15})/7),\;
\mI_7=[(6+\sqrt{15})/7,3/2),\;
\mI_8=[3/2,(1+\sqrt{5})/2),\;
\mI_9=[(1+\sqrt{5})/2,1+1/\sqrt{2}),\;
\mI_{10}=[1+1/\sqrt{2},2),\;
\mI_{11}=[2,\infty)$.
See Figure \ref{fig:asymptotics}.
Note that $\Cov_{\la}(r=1)=0$ and $\lim_{r \rightarrow \infty}\nu_{\la}(r)=0$ (at rate $O(r^{-2})$),
and $\argsup_{r \in [1,\infty)} \nu_{\la}(r)\approx 2.69$ with $\sup \nu_{\la}(r) \approx .0537$.


\begin{figure}
\centering
\psfrag{r}{\normalsize{$r$}}
\epsfig{figure=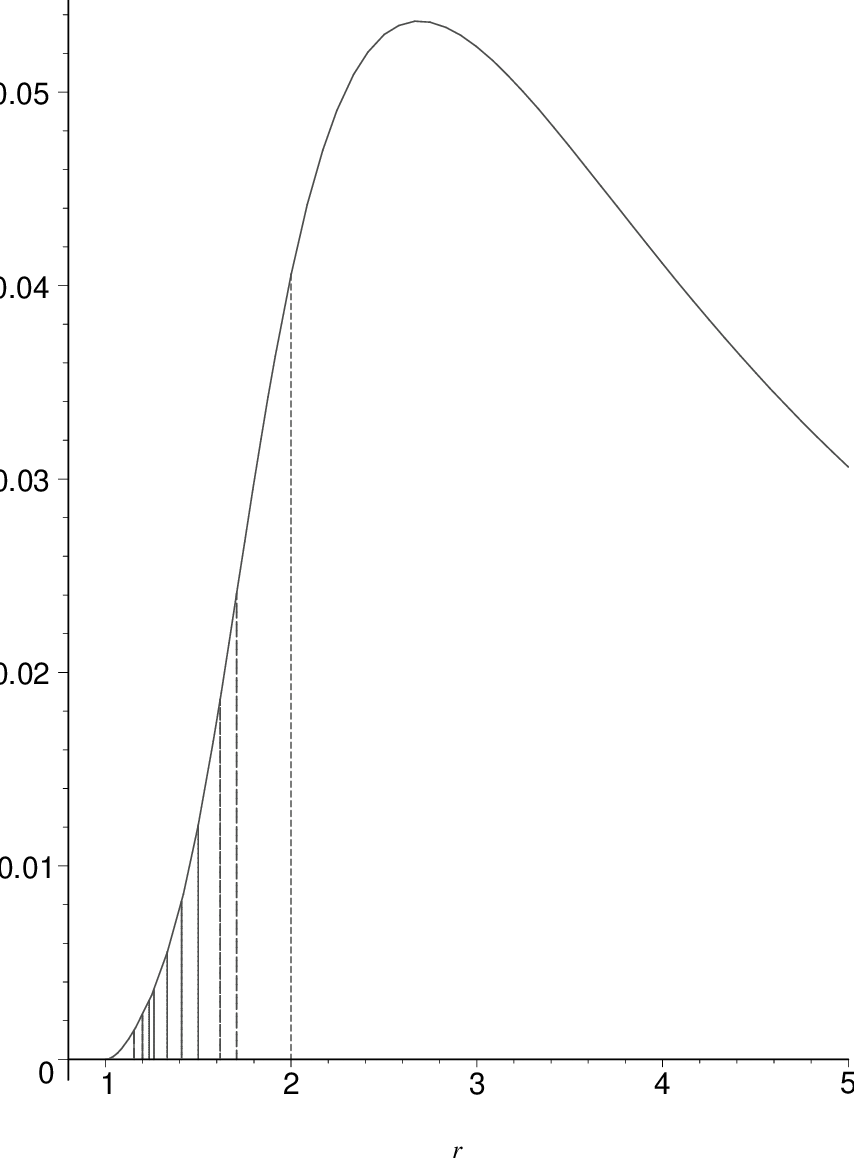, height=150pt, width=200pt}
\epsfig{figure=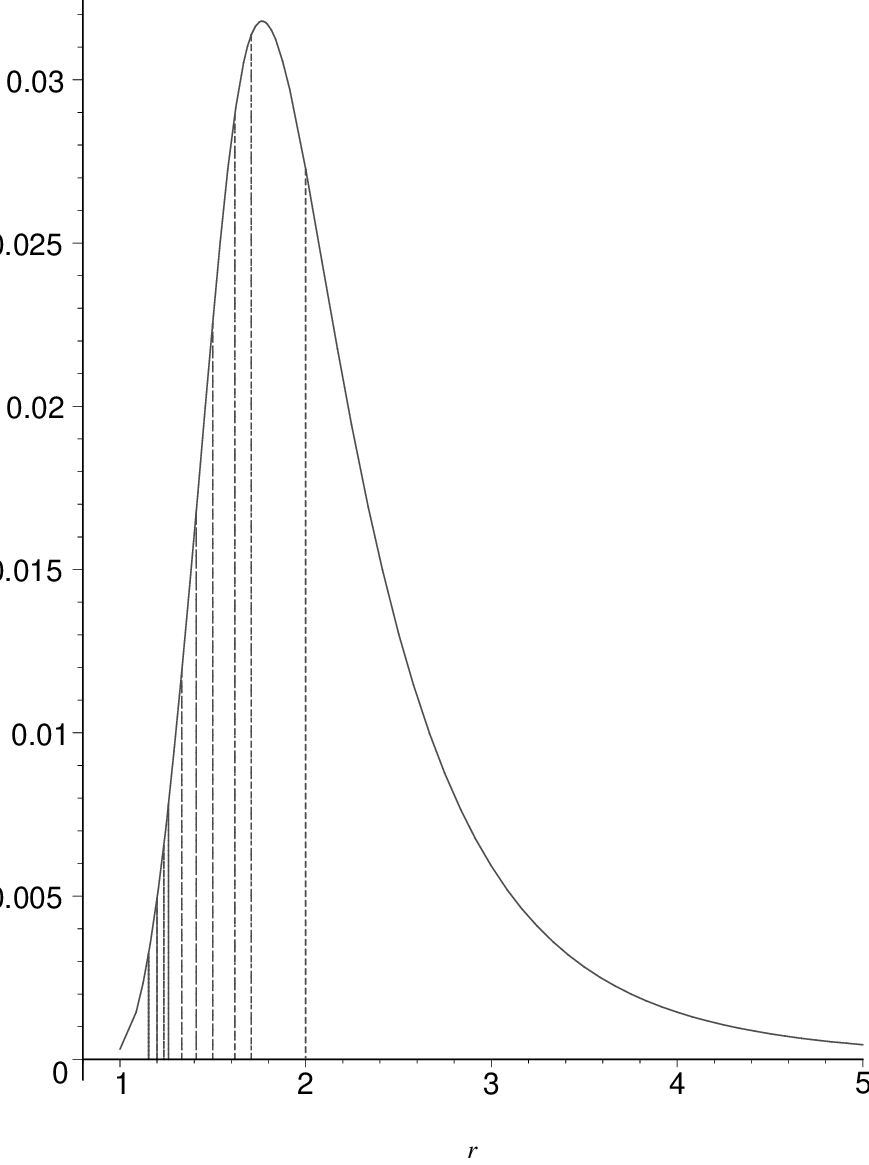, height=150pt, width=200pt}\\
\caption{
\label{cov of nu under-graph}
$\nu_{\la}(r)=\Cov\left[ h^{\la}_{12}(r),h^{\la}_{13}(r) \right]$ (left)
and
$\nu_{\lo}(r)=\Cov\left[ h^{\lo}_{12}(r),h^{\lo}_{13}(r) \right]$ (right)
as a function of $r$ for $r \in [1,5] $.}
\end{figure}

\subsection*{Appendix 2: The Asymptotic Variance of Relative Edge Density for the Underlying Graph}
The variance of $h^{\lo}_{12}(r)$,
denoted as $\Var_{\lo}(r)=\Var\left[ h^{\lo}_{12}(r) \right]$, is as follows:
$$\Var_{\lo}(r)=
\varphi^{\lo}_{1,1}(r)\I(r \in [1,4/3))+
\varphi^{\lo}_{1,2}(r)\I(r \in [4/3,3/2))+
\varphi^{\lo}_{1,3}(r)\I(r \in [3/2,2))+
\varphi^{\lo}_{1,4}(r)\I(r \in [2,\infty))
$$
where
$\varphi^{\lo}_{1,1}(r)=-\frac{(47\,r^6-195\,r^5+860\,r^4-846\,r^3-108\,r^2+720\,r-256)(752\,r^4-1170\,r^3-324\,r^2+47\,r^6-195\,r^5+720\,r-256)}{11664\,r^4(r+2)^2(r+1)^2}$,\\
$\varphi^{\lo}_{1,2}(r)=-\frac{(175\,r^5-579\,r^4+1450\,r^3-732\,r^2-536\,r+672)(1234\,r^3-1380\,r^2-968\,r+175\,r^5-579\,r^4+672)}{46656\,r^2(r+2)^2(r+1)^2}$,\\
$\varphi^{\lo}_{1,3}(r)=-\frac{(3\,r^8-7\,r^7-30\,r^6+84\,r^5-264\,r^4+304\,r^3+144\,r^2-368\,r+96)(-22\,r^6+108\,r^5-248\,r^4+3\,r^8-7\,r^7+304\,r^3+144\,r^2-368\,r+96)}{64\,r^8(r+2)^2(r+1)^2}$,\\
\\
$\varphi^{\lo}_{1,4}(r)=2\,\frac{(r^5+r^4-6\,r+2)(3\,r-1)}{r^8(r+1)^2}$.
See Figure \ref{var of nu under-graph}.

Note that $\Var_{\lo}(r=1)=2627/11664$ and $\lim_{r \rightarrow \infty}\Var_{\lo}(r)=0$ (at rate $O(r^{-4})$),
and $\argsup_{r \in [1,\infty)} \Var_{\lo}(r) \approx 1.44$ with $\sup \Var_{\lo}(r) \approx .25$.

The asymptotic variance for the underlying graph is
$$\nu_{\lo}(r):=\Cov[h^{\lo}_{12}(r),h^{\lo}_{13}(r)]=\sum_{i=1}^{11}\vartheta^{\lo}_i(r)\,\I(\mI_i)$$
where
{\small
\begin{multline*}
\vartheta^{\lo}_1(r)=-\frac{1}{58320\,(r^2+1)(2\,r^2+1)(r+1)^3(r+2)^3r^6}\,(1458\,r^{22}+13122\,r^{21}+50731\,r^{20}-84225\,r^{19}
-19193\,r^{18}-1823223\,r^{17}+\\
5576151\,r^{16}+2978697\,r^{15}-33432692\,r^{14}+37427862\,r^{13}+15883834\,r^{12}-60944766\,r^{11}+49876417\,r^{10}-1754523\,r^9-\\
36606859\,r^8+32338215\,r^7-10290256\,r^6-2234754\,r^5+7085471\,r^4-5608569\,r^3+1645826\,r^2-132876\,r+30824)
\end{multline*}
\begin{multline*}
\vartheta^{\lo}_2(r)=\vartheta^{\lo}_3(r)=-\frac{1}{116640\,(r^2+1)(2\,r^2+1)(r+1)^3(r+2)^3r^6}\,(1458\,r^{22}+13122\,r^{21}+62825\,r^{20}
-175011\,r^{19}+156014\,r^{18}-\\
3300900\,r^{17}+11053023\,r^{16}+5055135\,r^{15}-67685050\,r^{14}+75243552\,r^{13}+33155180\,r^{12}-120628524\,r^{11}+99831906\,r^{10}-\\
4883958\,r^9-74801558\,r^8+64360782\,r^7-19812000\,r^6-3667716\,r^5+14541630\,r^4-11254002\,r^3+3070468\,r^2-413208\,r+28880)
\end{multline*}
\begin{multline*}
\vartheta^{\lo}_4(r)=-\frac{1}{58320\,(r^2+1)(2\,r^2+1)(r^2-2)(r+2)^3(r+1)^3r^6}\,(972\,r^{24}+8748\,r^{23}+29590\,r^{22}-149106\,r^{21}-36820\,r^{20}-\\
986280\,r^{19}+5942884\,r^{18}+2883672\,r^{17}-47189711\,r^{16}+43450125\,r^{15}+85975304\,r^{14}-156173934\,r^{13}+27378901\,r^{12}+123606417\,r^{11}\\
-152209261\,r^{10}+64653597\,r^9+56621894\,r^8-88962768\,r^7+43754559\,r^6-5940597\,r^5-13006396\,r^4+17019366\,r^3-7037340\,r^2+\\
413208\,r-28880)
\end{multline*}
\begin{multline*}
\vartheta^{\lo}_5(r)=-\frac{1}{58320\,(r^2+1)(2\,r^2+1)(r+1)^3(r+2)^3r^6}\,(972\,r^{22}+8748\,r^{21}+31534\,r^{20}-131610\,r^{19}
+261546\,r^{18}-1552026\,r^{17}+\\
3745643\,r^{16}+4573731\,r^{15}-29416804\,r^{14}+26163354\,r^{13}+19600850\,r^{12}-43126062\,r^{11}+31497249\,r^{10}-7381467\,r^9-\\
22237963\,r^8+26778663\,r^7-9107024\,r^6-115074\,r^5+3136927\,r^4-5055609\,r^3+2292994\,r^2+14580\,r-1944)
\end{multline*}
\begin{multline*}
\vartheta^{\lo}_6(r)=\frac{1}{233280\,(r^2+1)(2\,r^2+1)(r+1)^3(r+2)^3r^6}\,(486\,r^{22}-7290\,r^{21}-181459\,r^{20}+1024401\,r^{19}
-2691213\,r^{18}+3921057\,r^{17}+\\
1844321\,r^{16}-33347697\,r^{15}+80028903\,r^{14}-29292735\,r^{13}-98093906\,r^{12}+125034492\,r^{11}-46658244\,r^{10}-57216612\,r^9+\\
88057996\,r^8-26383068\,r^7-12851392\,r^6+14179848\,r^5-8656508\,r^4+1593828\,r^3+134136\,r^2-58320\,r+7776)
\end{multline*}
\begin{multline*}
\vartheta^{\lo}_7(r)=\frac{1}{233280\,(r+2)^3(r^2+1)(2\,r^2+1)(r+1)^3(r-1)r^6}\,(486\,r^{23}-7776\,r^{22}-174169\,r^{21}
+1205860\,r^{20}-4656806\,r^{19}+\\
8763566\,r^{18}+7460036\,r^{17}-63559490\,r^{16}+91134324\,r^{15}+18516450\,r^{14}-122708655\,r^{13}+18577230\,r^{12}+80410332\,r^{11}-19357704\,r^{10}-\\
39129236\,r^9+75311048\,r^8-77449360\,r^7+4053376\,r^6+48283912\,r^5-40690240\,r^4+17736336\,r^3-4315680\,r^2+544320\,r-31104)
\end{multline*}
\begin{multline*}
\vartheta^{\lo}_8(r)=\frac{1}{960\,(r+2)^3(r^2+1)(2\,r^2+1)(r+1)^3r^8}\,(2\,r^{24}-30\,r^{23}-161\,r^{22}+107\,r^{21}+4137\,r^{20}
-10685\,r^{19}+8367\,r^{18}+\\
78713\,r^{17}-450859\,r^{16}+697707\,r^{15}+517846\,r^{14}-3723120\,r^{13}
+6565124\,r^{12}-1468692\,r^{11}-8695792\,r^{10}+9535720\,r^9-\\
6773160\,r^8+526744\,r^7+10691376\,r^6-7797264\,r^5+1137696\,r^4+523712\,r^3-2687872\,r^2+1701888\,r-245760)
\end{multline*}
\begin{multline*}
\vartheta^{\lo}_9(r)=\frac{1}{960\,(2\,r^2+1)(r+1)^2(r+2)^3(r^2+1)r^{10}}\,(2\,r^{25}-32\,r^{24}-129\,r^{23}+236\,r^{22}+4157\,r^{21}
-15610\,r^{20}+21289\,r^{19}+\\
67536\,r^{18}-511355\,r^{17}+1161830\,r^{16}-634128\,r^{15}-3001568\,r^{14}+9512164\,r^{13}-11014136\,r^{12}+2344968\,r^{11}+7126240\,r^{10}-\\
13850504\,r^9+14466592\,r^8-3823216\,r^7-4018976\,r^6+5155776\,r^5-4633984\,r^4+1959808\,r^3-244480\,r^2-3584\,r-1024)
\end{multline*}
\begin{multline*}
\vartheta^{\lo}_{10}(r)=\frac{1}{960\,(2\,r^2-1)(r+2)^3(r-1)(r+1)^2r^{10}}\,(2\,r^{24}-34\,r^{23}-101\,r^{22}+433\,r^{21}+5400\,r^{20}
-26982\,r^{19}+23049\,r^{18}+\\
166787\,r^{17}-717366\,r^{16}+1196092\,r^{15}+89468\,r^{14}-5130844\,r^{13}+12748688\,r^{12}-11274744\,r^{11}-12243496\,r^{10}+\\
33980568\,r^9-14886656\,r^8-19910592\,r^7+20667776\,r^6-1262208\,r^5-5402752\,r^4+2217088\,r^3-235776\,r^2-2560\,r-1024)
\end{multline*}
\begin{equation*}
\vartheta^{\lo}_{11}(r)=\frac{2}{15}\,{\frac
{180\,r^8-48\,r^7-648\,r^6+396\,r^5+214\,r^4-190\,r^3+39\,r^2-4\,r+1}{(2\,r^2-1)(r+1)^2r^{10}}}
\end{equation*}
}
\noindent
and
$\mI_1=[1,2/\sqrt{3}),\;
\mI_2=[2/\sqrt{3},6/5),\;
\mI_3=[6/5,\sqrt{5}-1),\;
\mI_4=[\sqrt{5}-1,(6+2\,\sqrt{2})/7),\;
\mI_5=[(6+2\,\sqrt{2})/7,4/3),\;
\mI_6=[4/3,(6+\sqrt{15})/7),\;
\mI_7=[(6+\sqrt{15})/7,3/2),\;
\mI_8=[3/2,(1+\sqrt{5})/2),\;
\mI_9=[(1+\sqrt{5})/2,1+1/\sqrt{2}),\;
\mI_{10}=[1+1/\sqrt{2},2),\;
\mI_{11}=[2,\infty)$.
See Figure \ref{fig:asymptotics}.
Note that $\Cov_{\lo}(r=1)=1/3240$ and $\lim_{r \rightarrow \infty}\nu_{\lo}(r)=0$ (at rate $O(r^{-6})$),
and $\argsup_{r \in [1,\infty)} \nu_{\lo}(r) \approx 1.765$ with $\sup \nu_{\lo}(r) \approx .0318$.

\begin{figure}[ht]
\centering
\rotatebox{-90}{ \resizebox{8.5 in}{!}{\includegraphics{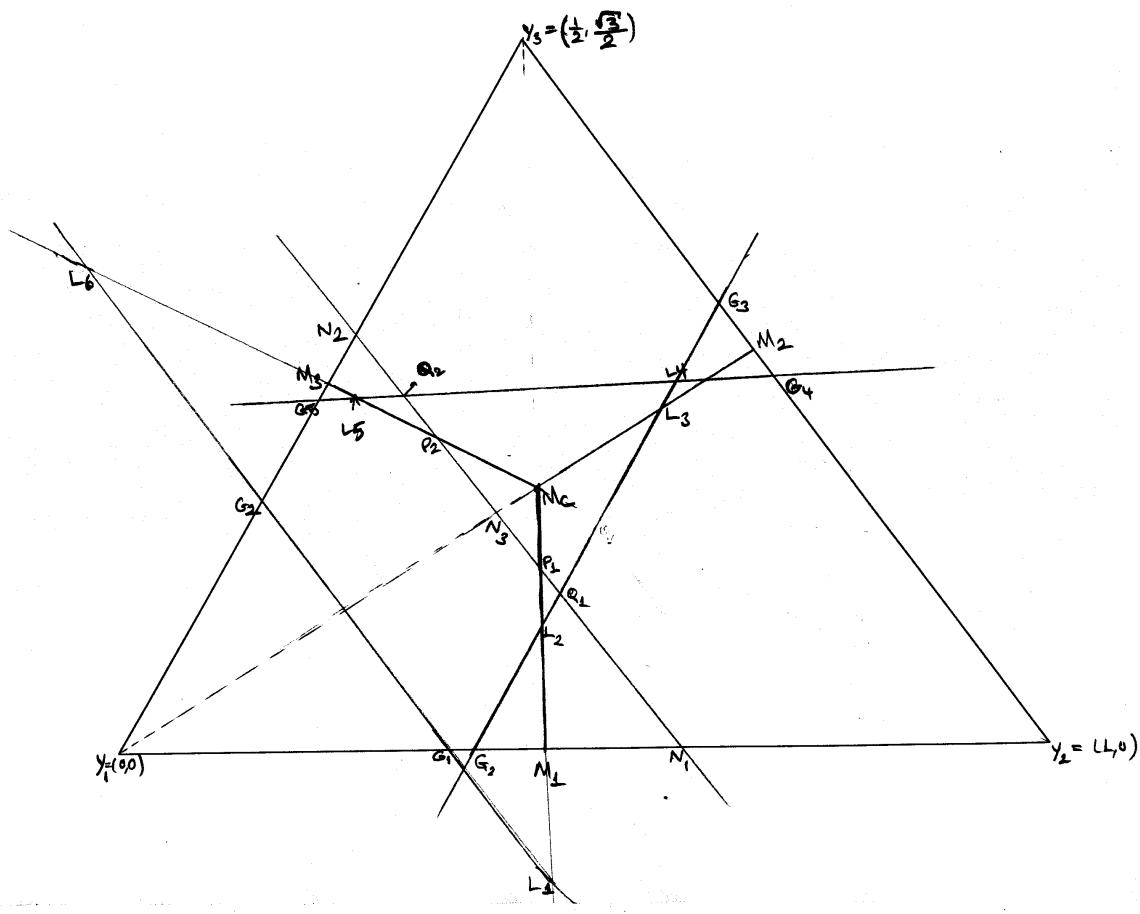}}}
\caption{
\label{fig:vertices-AND-OR}
An illustration of the vertices for possible types of $\NPE^r(x_1) \cap \G_1^r(x_1)$
for $x_1 \in T_s$.
}
\end{figure}
\begin{figure}[ht]
\centering
\rotatebox{-90}{ \resizebox{8.5 in}{!}{\includegraphics{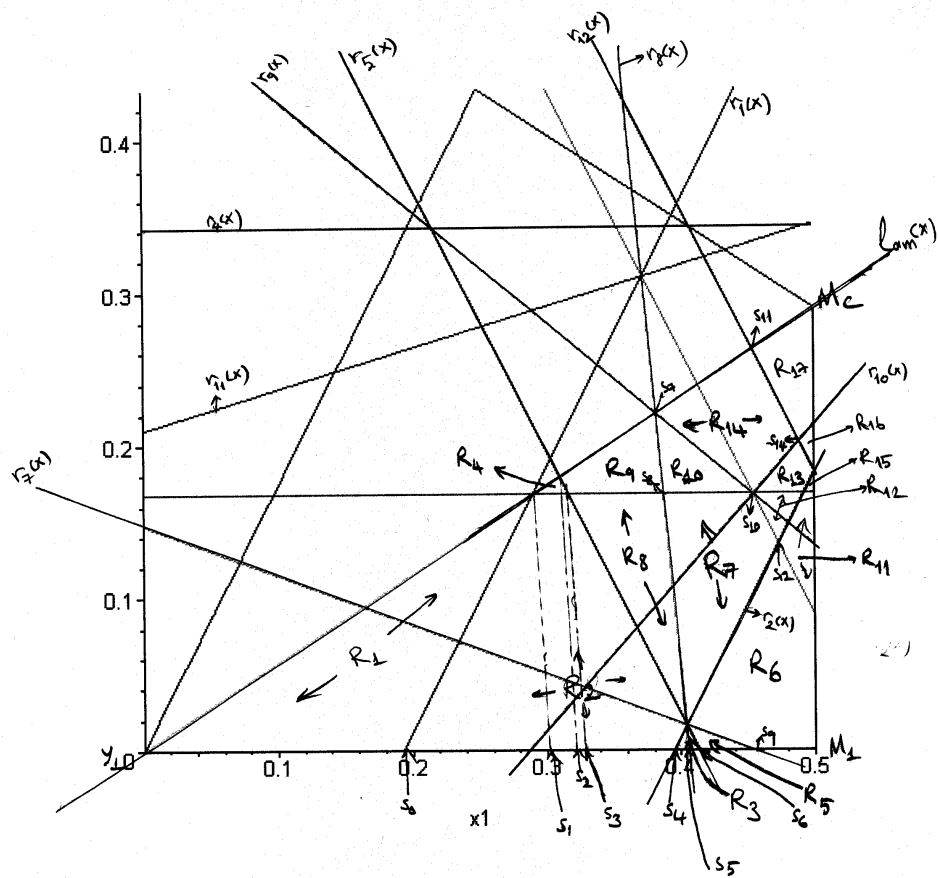}}}
\caption{
\label{fig:cases-AND-OR}
Prototype regions $R_i$ for various types of $\NPE^r(x_1) \cap \G_1^r(x_1)$
and the corresponding points whose $x$-coordinates are $s_k$ values.
}
\end{figure}

{\small
\section*{Appendix 3: Derivation of the Asymptotic Mean and Variance for Uniform Data}
In the standard equilateral triangle, let $\y_1=(0,0)$,
$\y_2=(1,0)$, $\y_3=\bigl( 1/2,\sqrt{3}/2 \bigr)$, $M_C$ be the
center of mass, $M_i$ be the midpoints of the edges $e_i$ for
$i=1,2,3$. Then $M_C=\bigl(1/2,\sqrt{3}/6\bigr)$,
$M_1=\bigl(3/4,\sqrt{3}/4 \bigr)$, $M_2=\bigl(1/4,\sqrt{3}/4\bigr)$,
$M_3=(1/2,0)$.
%
%
Let $\X_n$ be a random sample of size $n$ from $\U(\TY)$.
For $x_1=(u,v)$, $\ell_r(x_1)=r\,v+r\,\sqrt{3}\,u-\sqrt{3}\,x.$
Next, let $N_1:=\ell_r(x_1)\cap e_3$ and $N_2:=\ell_r(x_1)\cap e_2$.

\subsection*{Appendix 3.1: Derivation of $p_{\la}(r)$ and $\nu_{\la}(r)$ for Uniform Data}
\subsubsection*{Derivation of $\mu_\la(r)$ in Theorem \ref{thm:asy-norm-under}}
First we find $\mu_\la(r)$ for $r \in (1,\infty)$.
Observe that, by symmetry,
$$\mu_\la(r)=P\bigl( X_2 \in \NPE^r(X_1) \cap \G_1^r(X_1) \bigr)=
6\,P\bigl( X_2 \in \NY^r(X_1) \cap \G_1^r(X_1), X_1 \in T_s \bigr)$$
where $T_s$ is the triangle with vertices $\y_1$, $M_3$, and $M_C$.
Let $\ell_s(r,x)$ be the line such that $r\,d(\y_1,\ell_s(r,x))=d(\y_1,e_1)$,
so $\ell_s(r,x)=\sqrt{3}\,(1/r-x)$.
Then if $x_1 \in T_s$ is above $\ell_s(r,x)$ then $\NPE^r(x_1)=\TY$,
otherwise, $\NPE^r(x_1)\subsetneq \TY$.

\begin{figure} [ht]
\centering
\scalebox{.4}{\input{ls_lam_cases.pstex_t}}
\caption{
\label{fig:ls-lam-cases}
The cases for relative position of $\ell_s(r,x)$ with various $r$ values.
These are the prototypes for various types of $\NPE^r(x_1)$.
}
\end{figure}
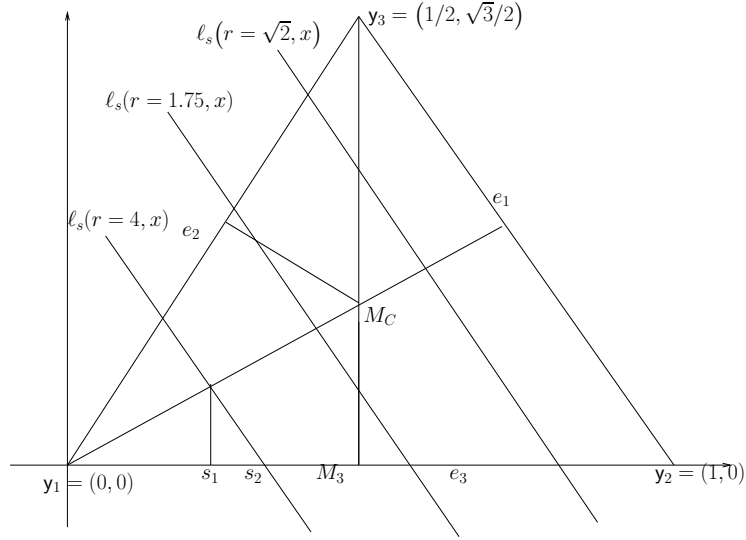

\begin{figure} [ht]
\centering
\scalebox{.27}{\input{G1ofxCase1.pstex_t}} 
\scalebox{.27}{\input{G1ofxCase2.pstex_t}}
\scalebox{.27}{\input{G1ofxCase3.pstex_t}}
\scalebox{.27}{\input{G1ofxCase4.pstex_t}}
\scalebox{.27}{\input{G1ofxCase5.pstex_t}}
\scalebox{.27}{\input{G1ofxCase6.pstex_t}}
\caption{
\label{fig:G_1-NYr-Cases-1}
The prototypes of the six cases of $\G^r_1\left(x\right)$ for $x \in T_s$ for $r \in [1,4/3)$.}
\end{figure}

To compute $\mu_\la(r)$,
we need to consider various cases for $\NPE^r(X_1)$ and $\G_1^r(X_1)$ given $X_1=(x,y) \in T_s$.
See Figures \ref{fig:ls-lam-cases} and \ref{fig:G_1-NYr-Cases-1}.
For any $x=(u,v) \in T(\Y)$, $\G^r_1(x)$ is a convex or nonconvex polygon.
Let $\xi_i(r,x)$ be the line  between $x$ and the vertex $\y_i$
parallel to the edge $e_i$ such that $r\,d(\y_i,\xi_i(r,x))=d(\y_i,\ell_r(x)) \text{ for } i=1,2,3.$
Then $\G^r_1(x)\cap R(\y_i)$ is bounded by $\xi_i(r,x)$ and the median lines.
For $x=(u,v)$,
$\xi_1(r,x)=-\sqrt{3}\,x+(v+\sqrt{3}\,u)/r,\; \xi_2(r,x)=
(v+\sqrt{3}r\,(x-1)+\sqrt{3}(1-u))/r \text{ and }
\xi_3(r,x)=(\sqrt{3}(r-1)+2\,v)/(2\,r).$
For $r\in \Bigl[6/5,\sqrt{5}-1)$,
there are six cases regarding $\G_1^r(x)$ and one case for $\NPE^r(x)$.
See Figure \ref{fig:G_1-NYr-Cases-1} for the prototypes of these six cases of $\G_1\left( x,\NY^r \right)$.
For the reflexivity graph case,
we determine the possible types of $\NPE^r(x_1) \cap \G_1^r(x_1)$
for $x_1 \in T_s$.
Depending on the location of $x_1$ and the value of the parameter $r$,
$\NPE^r(x_1) \cap \G_1^r(x_1)$ regions are polygons with various vertices.
See Figure \ref{fig:vertices-AND-OR} for the illustration of these vertices and
below for their explicit forms.

$G_1=\left( {\frac{\sqrt{3} y+3\, x}{3r}},0 \right)$,
$G_2=\left( -{\frac{\sqrt{3} y-3\,r+3-3\, x}{3r}},0\right) $,
$G_3=\left( -{\frac{\sqrt{3} y-6\,r+3-3\, x}{6r}},
-{\frac{\sqrt{3}\left( -\sqrt{3} y-3+3\, x \right) }{6r}}\right) $,
$G_4=\Bigl( {\frac{ \left( \sqrt{3}r+\sqrt{3}-2\, y \right) \sqrt{3}}{6 r}},\\
{\frac{\sqrt{3} \left( 3\,r -3+2\,\sqrt{3} y\right) }{6r}}\Bigr) $,
$G_5=\left( {\frac{ \left( \sqrt{3}r -\sqrt{3}+2\, y \right) \sqrt{3}}{6r}},
{\frac{\sqrt{3} \left( 3\,r -3+2\,\sqrt{3} y \right) }{6r}}\right) $,
$G_6=\left( {\frac{\sqrt{3} y+3\, x}{6r}},
{\frac{\sqrt{3} \left( \sqrt{3} y+3\, x \right) }{6r}}\right) $;

$P_1=\left( 1/2,\sqrt{3}/6\, \left( 2\,\sqrt{3}r\, y+6\,r\, x-3 \right) \right) $, and
$P_2=\left(-1/2+(\sqrt{3}r\, y+3\,r\, x)/2,
-\sqrt{3}/6\, \left( -3+\sqrt{3}r\, y+3\,r\, x \right)\right) $;

$L_1=\left( 1/2,{\frac{\sqrt{3} \left( 2\,\sqrt{3} y+6\, x-3\,r \right) }{6r}}\right) $,
$L_2=\left( 1/2,-{\frac{ \left( -2\,\sqrt{3} y-6+6\, x+3\,r \right) \sqrt{3}}{6r}}\right) $,
$L_3=\left( -{\frac{\sqrt{3} y-3\,r+3-3\, x}{2r}},
{\frac{\sqrt{3} \left(3\,r -\sqrt{3} y-3+3\, x \right) }{6r}}\right) $,
$L_4=\left( {\frac{3\,r -3+2\,\sqrt{3} y}{2r}},
{\frac{\sqrt{3} \left( 3\,r -3+2\,\sqrt{3} y \right) }{6r}}\right) $,
$L_5=\left( -{\frac{r -3+2\,\sqrt{3} y}{2r}},
{\frac{\sqrt{3} \left( 3\,r -3+2\,\sqrt{3} y \right) }{6r}}\right) $, and
$L_6=\left({\frac{-r+\sqrt{3} y+3\, x}{2r}},
-{\frac{\sqrt{3} \left( \sqrt{3} y+3\, x-3\,r \right) }{6r}}\right) $;
\noindent
$N_1=\left( \sqrt{3}r\, y/3+r\, x,0\right) $,
$N_2=\left( \sqrt{3}r\, y/6+r\, x/2,
\sqrt{3} \left( \sqrt{3} y/6+3\, x \right) r\right) $, and\\
$N_3=\left(\sqrt{3}r\, y/4+3\,r\, x/4,
\sqrt{3} \left( \sqrt{3} y/12+3\, x \right) r\right) $;
and
$Q_1=\left( {\frac{\sqrt{3}r^2 y+3\,r^2 x-\sqrt{
3} y+3\,r -3+3\, x}{6r}},
{\frac{ \left( \sqrt{3}r^2 y+3\,r^2 x+
\sqrt{3} y-3\,r+3-3\, x \right) \sqrt{3}}{6r}}\right) $,\\
and
$Q_2=\left( {\frac{2\,\sqrt{3}r^2 y+6\,r^2 x-3\,r+3-2\,\sqrt{3} y}{6r}},
{\frac{\sqrt{3} \left( 3\,r -3+2\,\sqrt{3} y \right) }{6r}}\right) $.

Let $\msP(a_1,a_2,\ldots, a_n)$ denote the polygon with vertices $a_1,a_2,\ldots, a_n $.
For $r \in \bigl[1,4/3\bigr)$,
there are 14 cases to consider for calculation of $p_{\la}(r)$ in the reflexivity graph version.
Each of these cases correspond to the regions in Figure \ref{fig:cases-AND-OR},
where Case 1 corresponds to $R_i$ for $i=1,2,3,4$,
and Case $j$ for $j >1$ corresponds to $R_{j+3}$ for $j=1,2,\ldots,14$.
These regions are bounded by various combinations of the lines defined below.

Let $\ell_{am}(x)$ be the line joining $\y_1$ to $M_C$,
then $\ell_{am}(x)=\sqrt{3}x/3$.
Let also
$r_1(x)=\sqrt{3} \left( 2\,r+3\, x-3 \right)/3$,
$r_2(x)=\sqrt{3}/2-\sqrt{3}r/3$,
$r_3(x)=\left( 2\, x-2+r \right) \sqrt{3}/2$,
$r_4(x)=\sqrt{3}/2-\sqrt{3}r/4$,
$r_5(x)=-{\frac{\sqrt{3} \left( 2\,r\, x-1 \right) }{2r}}$,
$r_6(x)=-{\frac{\sqrt{3} \left( -2+3\,r\, x \right) }{3r}}$,
$r_7(x)=-{\frac{ \left( 1+r^2 x-r - x \right) \sqrt{3}}{r^2+1}}$,
$r_8(x)=-{\frac{ \left( r^2 x-1+ x \right) \sqrt{3}}{r^2-1}}$,
$r_9(x)=-{\frac{ \left( r^2 x-1 \right) \sqrt{3}}{r^2+2}}$,
$r_{10}(x)=-{\frac{ \left( -2\,r+2+r^2 x \right) \sqrt{3}}{-4+r^2}}$,
$r_{11}(x)=-{\frac{ \left( -2\,r+2-2\, x+r^2 x \right) \sqrt{3}}{r^2+2}}$,
$r_{12}(x)=-\left( 2\, x-r \right) \sqrt{3}/2$, and
$r_{13}(x)=-\left( -1+ x \right) \sqrt{3}/3$.
Furthermore, to determine the integration limits,
we specify the $x$-coordinate of the boundaries of these regions using $s_k$  for $k=0,1,\ldots,14$.
See also Figure \ref{fig:cases-AND-OR} for an illustration of these points
whose explicit forms are provided below.

$s_0=1-2\,r/3$,
$s_1=3/2-r$,
$s_2=3/(8\,r)$,
$s_3={\frac{-3\,r+2\,r^2+3}{6r}}$,
$s_4=1-r/2$,
$s_5={\frac{2\,r-r^2+1}{4r}}$,
$s_6=1/(2\,r)$,
$s_7=\frac{3}{2\, \left( 2\,r^2+1 \right)}$,
$s_8={\frac{9-3\,r^2+2\,r^3-2\,r}{6(r^2+1)}}$,
$s_9=1/\left( r+1 \right)$,
$s_{10}={\frac{-3\,r+2\,r^2+4}{6r}}$,
$s_{11}=3\,r/8$,
$s_{12}={\frac{6\,r-3\,r^2+4}{12r}}$,
$s_{13}=3/2-5\,r/6$, and
$s_{14}=r-1/2-r^3/8$.

%

Below, we compute $P(X_2 \in \NPE^r(X_1) \cap \G_1^r(X_1), X_1 \in T_s)$ for each of the 14 cases:
{\small
\noindent \textbf{Case 1:}
\begin{multline*}
P(X_2 \in \NPE^r(X_1) \cap \G_1^r(X_1), X_1 \in T_s)=
\left(\int_{0}^{s_2}\int_{0}^{\ell_{am}(x)} + \int_{s_2}^{s_6}\int_{0}^{r_5(x)}\right)
\frac{A(\msP(G_1,N_1,N_2,G_6))}{A(\TY)^2}dydx=\\
{\frac{ \left( r -1 \right)  \left( r+1 \right)
 \left( r^2+1 \right) }{64\,r^6}}
\end{multline*}
where
$A(\msP(G_1,N_1,N_2,G_6))=
\sqrt{3}/36\, \left( \sqrt{3} y+3\, x \right)^2r^2-
{\frac{\sqrt{3} \left( \sqrt{3} y+3\, x
\right)^2}{36\,r^2}}
$.

\noindent \textbf{Case 2:}
\begin{multline*}
P(X_2 \in \NPE^r(X_1) \cap \G_1^r(X_1), X_1 \in T_s)=
\left(\int_{s_5}^{s_6}\int_{r_5(x)}^{r_7(x)} + \int_{s_6}^{s_9}\int_{0}^{r_7(x)}\right)
\frac{A(\msP(G_1,N_1,P_2,M_3,G_6))}{A(\TY)^2}dydx=\\
{\frac{ \left( 9\,r^5+23\,r^4+24\,r^3+24\,r^2+13\,r+3 \right)
\left( r -1 \right)^4}{96\, r^6 \left( r+1 \right)^3}}
\end{multline*}
where
$A(\msP(G_1,N_1,P_2,M_3,G_6))=
-{\frac{\sqrt{3} \left( -4\,r^3\sqrt{3} y-12\,r^3 x+
2\,r^4\,y^2+4\,r^4\sqrt{3}y\, x+6\,r^4 x^2+3\,r^2+2\,y^2+
4\,\sqrt{3} y\, x+6\, x^2 \right) }{24\,r^2}}$.

\noindent \textbf{Case 3:}
\begin{multline*}
P(X_2 \in \NPE^r(X_1) \cap \G_1^r(X_1), X_1 \in T_s)=\\
\left( \int_{s_5}^{s_9}\int_{r_7(x)}^{r_3(x)} +
\int_{s_9}^{s_{12}}\int_{0}^{r_3(x)} +
\int_{s_{12}}^{1/2}\int_{0}^{r_6(x)} \right)
\frac{A(\msP(G_1,G_2,Q_1,P_2,M_3,G_6))}{A(\TY)^2}dydx=\\
{\frac{324\,r^{11}-1620\,r^{10}-618\,r^9+4626\,r^8+990\,r^7-2454\,r^6+2703\,r^5-
5571\,r^4-3827\,r^3+1455\,r^2+3072\,r+1024}
{7776\, \left( r+1 \right)^3r^6}}
\end{multline*}
where
$A(\msP(G_1,G_2,Q_1,P_2,M_3,G_6))=
-\Bigl[\sqrt{3} \bigl( -4\,\sqrt{3}r\, y-12\,x+4\,y^2+4\,r^2\,y^2-12\,r+9\,r^2+12
\,r\, x+4\,r^4\,y^2-12\, x^2r^2-24\,r^3 x+12\,r^4 x^2+8\,r^4\sqrt{
3} y\, x+12\, x^2+12\,r^2 x+6-8\,r^3\sqrt{3} y+4\,\sqrt{3} y+4\,\sqrt{3}r^2 y \bigr)\Bigr]\Big/\Bigl[24\,r^2\Bigr]
$.

\noindent \textbf{Case 4:}
\begin{multline*}
P(X_2 \in \NPE^r(X_1) \cap \G_1^r(X_1), X_1 \in T_s)=\\
\left(\int_{s_8}^{s_5}\int_{r_8(x)}^{r_2(x)}+
\int_{s_5}^{s_{10}}\int_{r_3(x)}^{r_2(x)}+
\int_{s_{10}}^{s_{12}}\int_{r_3(x)}^{r_6(x)}\right)
\frac{A(\msP(G_1,M_1,L_2,Q_1,P_2,M_3,G_6))}{A(\TY)^2}dydx=\\
\Bigl[512+138240\,r^7+3654\,r^{12}-255
\,r^8+43008\,r^3-12369\,r^2-86387\,r^4-193581
\,r^6+148224\,r^5-100608\,r^9+94802\,r^{10}-\\
35328\,r^{11}\Bigr]\Big/
\Bigl[7776\, \left( r^2+1 \right)^3r^6\Bigr]
\end{multline*}
where
$A(\msP(G_1,M_1,L_2,Q_1,P_2,M_3,G_6))=
-\Bigl[\sqrt{3} \bigl( 6\, x+3\,r^2-2\,\sqrt{3} y+2\,\sqrt{3}r^2 y+2\,r^4\,y^2-4\,
r^3\sqrt{3} y+4\,\sqrt{3} y\, x+2\,r^2 y^2+4\,r^4\sqrt{3} y\, x-6\, x^2 r^2-
12\,r^3 x+6\,r^4 x^2+6\,r^2 x-3 \bigr) \Bigr]\Big/\Bigl[12\,r^2\Bigr]
$.

\noindent \textbf{Case 5:}
\begin{multline*}
P(X_2 \in \NPE^r(X_1) \cap \G_1^r(X_1), X_1 \in T_s)=
\left(\int_{s_3}^{s_8}\int_{r_5(x)}^{r_2(x)} +
\int_{s_8}^{s_5}\int_{r_5(x)}^{r_8(x)} \right)
\frac{A(\msP(G_1,M_1,P_1,P_2,M_3,G_6))}{A(\TY)^2}dydx=\\
-{\frac{ \left( 177\,r^8-648\,r^7+570\,r^6-360\,r^5+28\,r^4-24\,r^3+174\,r^2+
72\,r+27 \right)  \left( -12\,r+7\,r^2+3 \right)^2}{7776\,\left( r^2+1 \right)^3r^6}}
\end{multline*}
where
$A(\msP(G_1,M_1,L_2,Q_1,P_2,M_3,G_6))=
-{\frac{\sqrt{3} \left( -4\,r^3\sqrt{3} y-12\,r^3 x+
3\,r^2+6\,r^4\sqrt{3} y\, x+9\,r^4 x^2+3\,r^4\,y^2+\,y^2+2
\,\sqrt{3} y\, x+3\, x^2 \right) }{12\,r^2}}
$.

\noindent \textbf{Case 6:}
\begin{multline*}
P(X_2 \in \NPE^r(X_1) \cap \G_1^r(X_1), X_1 \in T_s)=\\
\left( \int_{s_2}^{s_3}\int_{r_5(x)}^{\ell_{am}(x)} +
\int_{s_3}^{s_7}\int_{r_2(x)}^{\ell_{am}(x)} +
\int_{s_7}^{s_8}\int_{r_2(x)}^{r_8(x)} \right)
\frac{A(\msP(G_1,M_1,P_1,P_2,M_3,G_6))}{A(\TY)^2}dydx=\\
\Bigl[137472\,r^{18}-952704\,r^{17}+2792712\,r^{16}-5116608\,r^{15}+7057828\,r^{14}-
7725792\,r^{13}+7022682\,r^{12}-5484816\,r^{11}+\\
3631995\,r^{10}-2213712\,r^9+1213271\,r^8-578976\,r^7+292518\,r^6-
101952\,r^5+36612\,r^4-11664\,r^3+3051\,r^2-1296\,r+243\Bigr]\Big/\\
\Bigl[ \left(15552\, r^2+1 \right)^3 \left( 2\,r^2+1 \right)^3r^6\Bigr]
\end{multline*}
where
$A(\msP(G_1,M_1,P_1,P_2,M_3,G_6))=
-{\frac{\sqrt{3} \left( -4\,r^3\sqrt{3} y-12\,r^3 x+
3\,r^2+6\,r^4\sqrt{3} y\, x+9\,r^4 x^2+3\,r^4\,y^2+\,y^2+2
\,\sqrt{3} y\, x+3\, x^2 \right) }{12\,r^2}}
$.

\noindent \textbf{Case 7:}
\begin{multline*}
P(X_2 \in \NPE^r(X_1) \cap \G_1^r(X_1), X_1 \in T_s)=
\left(\int_{s_7}^{s_8}\int_{r_8(x)}^{r_9(x)}+
\int_{s_8}^{s_{10}}\int_{r_2(x)}^{r_9(x)}\right)
\frac{A(\msP(G_1,M_1,L_2,Q_1,P_2,M_3,G_6))}{A(\TY)^2}dydx=\\
-{\frac{ 4\left( 100\,r^{11}-408\,r^{10}+454\,r^9-564\,r^8+283\,r^7-108\,r^6-34\,r^5+
204\,r^4-r^3+132\,r^2+26\,r+24 \right)  \left( 2\,r -1 \right)^2 \left( r -1 \right)^2}
{243 \left( r^2+1 \right)^3r^3 \left( 2\,r^2+1 \right)^3}}
\end{multline*}
where
$A(\msP(G_1,M_1,L_2,Q_1,P_2,M_3,G_6))=
-\Bigl[\sqrt{3} \bigl( 6\, x+3\,r^2-2\,\sqrt{3} y+2\,\sqrt{3}r^2 y+2\,r^4\,y^2-
4\,r^3\sqrt{3} y+4\,\sqrt{3} y\, x+2\,r^2 y^2+4\,r^4\sqrt{3} y\, x-
6\,x^2 r^2-12\,r^3 x+6\,r^4 x^2+6\,r^2 x-3 \bigr) \Bigr] \Big / \Bigl[12\,r^2\Bigr]
$.

\noindent \textbf{Case 8:}
\begin{multline*}
P(X_2 \in \NPE^r(X_1) \cap \G_1^r(X_1), X_1 \in T_s)=
\left( \int_{s_{12}}^{s_{13}}\int_{r_6(x)}^{r_3(x)} +
\int_{s_{13}}^{1/2}\int_{r_6(x)}^{r_2(x)} \right)
\frac{A(\msP(G_1,G_2,Q_1,N_3,M_C,M_3,G_6))}{A(\TY)^2}dydx=\\
\Bigl[ \left( -2+r \right)  \bigl( 2369\,r^{11}-11342\,r^{10}+29934\,r^9-50340\,r^8+
54056\,r^7-51824\,r^6+48320\,r^5-20864\,r^4-640\,r^3\\
-1280\,r^2+512\,r+1024 \bigr) \Bigr] \Big/ \Bigl[15552\,r^6\Bigr]
\end{multline*}
where
$A(\msP(G_1,G_2,Q_1,N_3,M_C,M_3,G_6))=
-\Bigl[\sqrt{3} \bigl( 4\,\sqrt{3}r^2 y-12\, x-12\,r+5\,r^2+12\,r\, x+4\,y^2-
12\,x^2r^2+4\,r^2\,y^2+r^4\,y^2+2\,r^4\sqrt{3} y\, x-4\,r^3\sqrt{3} y+
6-12\,r^3 x+3\,r^4 x^2+12\, x^2+12\,r^2 x-4\,\sqrt{3}r\, y+4\,\sqrt{3} y \bigr) \Bigr] \Big/ \Bigl[24\,r^2\Bigr]
$.

\noindent \textbf{Case 9:}
\begin{multline*}
P(X_2 \in \NPE^r(X_1) \cap \G_1^r(X_1), X_1 \in T_s)=
\left( \int_{s_{10}}^{s_{12}}\int_{r_6(x)}^{r_2(x)} +
\int_{s_{12}}^{s_{13}}\int_{r_3(x)}^{r_2(x)} \right)
\frac{A(\msP(G_1,M_1,L_2,Q_1,N_3,M_C,M_3,G_6))}{A(\TY)^2}dydx=\\
-{\frac{ \left( 49\,r^8-168\,r^7+354\,r^6-528\,r^5+236\,r^4-96\,r^3-224\,r^2
+384\,r+64 \right)  \left( -12\,r+7\,r^2+4 \right)^2}{15552\,r^6}}
\end{multline*}
where
$A(\msP(G_1,M_1,L_2,Q_1,N_3,M_C,M_3,G_6))=
-\Bigl[\sqrt{3} \bigl( 8\,\sqrt{3} y\, x+4\,\sqrt{3}r^2 y+12\, x+
2\,r^2-12\, x^2 r^2-4\,r^3\sqrt{3} y-12\,r^3 x+3\,r^4 x^2+
r^4\,y^2+2\,r^4\sqrt{3}y\,x+12\,r^2 x-6-4\,\sqrt{3} y+
4\,r^2\,y^2 \bigr) \Bigr] \Big/ \Bigl[24 \, r^2\Bigr]
$.

\noindent \textbf{Case 10:}
\begin{multline*}
P(X_2 \in \NPE^r(X_1) \cap \G_1^r(X_1), X_1 \in T_s)=\\
\left(\int_{s_{10}}^{s_{14}}\int_{r_2(x)}^{r_{10}(x)} +
\int_{s_{14}}^{s_{13}}\int_{r_2(x)}^{r_{12}(x)} +
\int_{s_{13}}^{1/2}\int_{r_3(x)}^{r_{12}(x)} \right)
\frac{A(\msP(G_1,M_1,L_2,Q_1,N_3,L_4,L_5,M_3,G_6))}{A(\TY)^2}dydx=\\
-{\frac{6144+195456\,r^6+324\,r^{11}-
76720\,r^7-801792\,r^2+217856\,r+946432\,r^3-
239904\,r^5-275328\,r^4+39408\,r^8-11849\,r^9}{31104\,r^3}}
\end{multline*}
where
$A(\msP(G_1,M_1,L_2,Q_1,N_3,L_4,L_5,M_3,G_6))=
-\Bigl[\sqrt{3} \bigl( 4\,\sqrt{3}r^2 y+8\,\sqrt{3} y\, x+
4\,r^2\,y^2-16\,\sqrt{3} r\, y-4\,r^3\sqrt{3} y-24\,y^2+
12\, x+24\,r -6\,r^2-12\, x^2r^2-12\,r^3 x+3\,r^4 x^2+12\,r^2 x+20\,\sqrt{3
} y+2\,r^4\sqrt{3} y\, x+r^4\,y^2-24 \bigr)\Bigr] \Big/ \Bigl[24\,r^2\Bigr]
$.

\noindent \textbf{Case 11:}
\begin{multline*}
P(X_2 \in \NPE^r(X_1) \cap \G_1^r(X_1), X_1 \in T_s)=\\
\left(\int_{s_7}^{s_{11}}\int_{r_9(x)}^{\ell_{am}(x)} +
\int_{s_{11}}^{s_{10}}\int_{r_9(x)}^{r_{12}(x)} +
\int_{s_{10}}^{s_{14}}\int_{r_{10}(x)}^{r_{12}(x)} \right)
\frac{A(\msP(G_1,M_1,L_2,Q_1,Q_2,L_5,M_3,G_6))}{A(\TY)^2}dydx=\\
\Bigl[ \left( r -1 \right)  \bigl( 1080\,r^{16}+1080\,r^{15}-17820\,r^{14}-
540\,r^{13}+65394\,r^{12}-46926\,r^{11}+105435\,r^{10}-261765\,r^9+229286\,r^8-180586\,r^7+\\
101638\,r^6+40774\,r^5-46112\,r^4+24448\,r^3-20224\,r^2+10496\,r -6144 \bigr) \Bigr]
\Big/ \Bigl[10368 \, r^3 \left( 2\,r^2+1 \right)^3\Bigr]
\end{multline*}
where
$A(\msP(G_1,M_1,L_2,Q_1,Q_2,L_5,M_3,G_6))=
-\Bigl[\sqrt{3} \bigl( 6\, x+3\,r^2-4\,r^2 x\,\sqrt{3} y-4\,y^2-
6\, x^2r^2+2\,r^4\sqrt{3} y\, x+4\,\sqrt{3} y\,x-2\,r^2\,y^2-
4\,r^3\sqrt{3} y+r^4\,y^2-12\,r^3 x+3\,r^4 x^2+12\,r^2 x-6+
4\,\sqrt{3}r^2 y+2\,\sqrt{3} y \bigr) \Bigr]\Big/\Bigl[12\,r^2\Bigr]
$.

\noindent \textbf{Case 12:}
\begin{multline*}
P(X_2 \in \NPE^r(X_1) \cap \G_1^r(X_1), X_1 \in T_s)=
\int_{s_{13}}^{1/2}\int_{r_2(x)}^{r_3(x)} \frac{A(\msP(G_1,G_2,Q_1,N_3,L_4,L_5,M_3,G_6))}{A(\TY)^2}dydx=\\
-{\frac{ \left( 49\,r^6-204\,r^5+476\,r^4-768\,r^3-8\,r^2+768\,r -288 \right)  \left( -6+
5\,r \right)^2}{7776\,r^2}}
\end{multline*}
where
$A(\msP(G_1,G_2,Q_1,N_3,L_4,L_5,M_3,G_6))=
-\Bigl[\sqrt{3} \bigl( -12\, x+12\,r -3\,r^2+12\,r\, x-20\,\sqrt{3}r\, y-
12\, x^2r^2+4\,\sqrt{3}r^2 y-12\,r^3 x+3\,r^4 x^2+28\,\sqrt{3} y+
12\, x^2+12\,r^2 x-12-20\,y^2+4\,r^2\,y^2-4\,r^3 \sqrt{3} y+
r^4\,y^2+2\,r^4\sqrt{3} y\, x \bigr) \Bigr] \Big/ \Bigl[24\,r^2\Bigr]
$.

\noindent \textbf{Case 13:}
\begin{multline*}
P(X_2 \in \NPE^r(X_1) \cap \G_1^r(X_1), X_1 \in T_s)=
\int_{s_{14}}^{1/2}\int_{r_{12}(x)}^{r_{10}(x)} \frac{A(\msP(L_1,L_2,Q_1,N_3,L_4,L_5,L_6))}{A(\TY)^2}dydx=\\
{\frac{ \left( 4\,r^7+8\,r^6-37\,r^5-58\,r^4-84\,r^3+168\,r^2+336\,r -352 \right)
\left( -2+r \right)  \left( r^2+2\,r -4 \right)^2}{384\, \left( r+2 \right)^2r^2}}
\end{multline*}
where
$A(\msP(L_1,L_2,Q_1,N_3,L_4,L_5,L_6))=
-\Bigl[\sqrt{3} \bigl( -4\,r^3\sqrt{3} y-8\,\sqrt{3}r\, y+12\, x+24\,r -
8\,\sqrt{3} y\,x-12\,r^2+24\,r\, x-24-12\, x^2r^2+4\,\sqrt{3}r^2 y-
32\,y^2-12\,r^3 x+3\,r^4 x^2+20\,\sqrt{3} y-24\, x^2+
12\,r^2 x+2\,r^4\sqrt{3} y\, x+r^4\,y^2+4\,r^2\,y^2 \bigr) \Big/ \Bigl[24\,r^2\Bigr]
$.

\noindent \textbf{Case 14:}
\begin{multline*}
P(X_2 \in \NPE^r(X_1) \cap \G_1^r(X_1), X_1 \in T_s)=
\left( \int_{s_{11}}^{s_{14}}\int_{r_{12}(x)}^{\ell_{am}(x)} +
\int_{s_{14}}^{1/2}\int_{r_{10}(x)}^{\ell_{am}(x)} \right)
\frac{A(\msP(L_1,L_2,Q_1,Q_2,L_5,L_6))}{A(\TY)^2}dydx=\\
-\Bigl[ \bigl( 135\,r^{11}+675\,r^{10}-1350\,r^9-9450\,r^8+702\,r^7+39150\,r^6+24272
\,r^5-47432\,r^4-135040\,r^3+57088\,r^2+204800\,r -\\
134144 \bigr)
\left( r -1 \right) \Bigr]\Big/ \Bigl[10368 \left( r+2 \right)^2 r^2\Bigr]
\end{multline*}
where
$A(\msP(L_1,L_2,Q_1,Q_2,L_5,L_6))=
-\Bigl[\sqrt{3} \bigl( -4\,r^3\sqrt{3} y+4\,\sqrt{3}r\, y+r^4\,y^2+6\, x-
4\,\sqrt{3} y\, x+2\,r^4\sqrt{3} y\, x+12\,r\,x-
4\,r^2 x\,\sqrt{3} y-6\, x^2 r^2+4\,\sqrt{3}r^2 y-12\,r^3 x+3\,r^4 x^2+
2\,\sqrt{3} y-12\, x^2+12\,r^2 x-6-8\,y^2-2\,r^2\,y^2 \bigr) \Bigr]\Big/ \Bigl[12\,r^2\Bigr]$.
}

Adding up the $P(X_2 \in \NPE^r(X_1) \cap \G_1^r(X_1), X_1 \in T_s)$
values in the 14 possible cases above, and multiplying by 6
we get for $r \in [1,4/3)$,
$$\mu_\la(r)=
-{\frac{ \left( r -1 \right)  \left( 5\,r^5-
148\,r^4+245\,r^3-178\,r^2-232\,r+128 \right) }
{54\,r^2 \left( r+2 \right)  \left( r+1 \right) }}.$$
The $\mu_\la(r)$ values for the other intervals can be calculated similarly.
For $r=\infty $, $\mu_\la(r)=1$ follows trivially.

\subsubsection*{Derivation of $\nu_\la(r)$ in Theorem \ref{thm:asy-norm-under}}

By symmetry,
$P(\{X_2,X_3\} \subset \NPE^r(X_1)\cap \G_1^r(X_1))=
6\,P(\{X_2,X_3\} \subset \NPE^r(X_1)\cap \G_1^r(X_1),\; X_1 \in T_s)$.

For $r \in \bigl[6/5,\sqrt{5}-1\bigr)$,
there are 14 cases to consider for calculation of $\nu_{\la}(r)$ in the reflexivity graph version:
{\small
\noindent \textbf{Case 1:}
\begin{multline*}
P(\{X_2,X_3\} \subset \NPE^r(X_1)\cap \G_1^r(X_1),\; X_1 \in T_s)=
\left(\int_{0}^{s_2}\int_{0}^{\ell_{am}(x)}+
\int_{s_2}^{s_6}\int_{0}^{r_5(x)} \right)
\frac{A(\msP(G_1,N_1,N_2,G_6))^2}{A(\TY)^3}dydx=\\
{\frac{ \left( r^2+1 \right)^2 \left( r+1 \right)^2 \left( r -1 \right)^2}{384\,r^{10}}}
\end{multline*}
where
$A(\msP(G_1,N_1,N_2,G_6))=
\sqrt{3} \left( \sqrt{3} y+3\, x \right)^2 r^2/36-
{\frac{ \left( \sqrt{3} y+3\, x \right)^2\sqrt{3}}{36\,r^2}}
$.

\noindent \textbf{Case 2:}
\begin{multline*}
P(\{X_2,X_3\} \subset \NPE^r(X_1)\cap \G_1^r(X_1),\; X_1 \in T_s)=
\left(\int_{s_5}^{s_6}\int_{r_5(x)}^{r_7(x)} +
\int_{s_6}^{s_9}\int_{0}^{r_7(x)} \right)
\frac{A(\msP(G_1,N_1,P_2,M_3,G_6))^2}{A(\TY)^3}dydx=\\
{\frac { \left( 5+38\,r+137\,r^2+320\,r^3+
552\,r^4+736\,r^5+792\,r^6+640\,r^7+407\,r^8+178\,r^9+
35\,r^{10} \right)  \left( -1+r \right)^5}{960\,r^{10} \left( r+1 \right)^5}}
\end{multline*}
where
$A(\msP(G_1,N_1,P_2,M_3,G_6))=
-{\frac {\sqrt{3} \left( -4\,r^3\sqrt{3}y-12\,r^3 x+
2\,r^4 y^2+4\,r^4\sqrt{3}y\,x+6\,r^4x^2+3\,r^2+
2\,y^2+ 4\,\sqrt{3}y\,x+6\,x^2 \right) }{24\,r^2}}
$.

\noindent \textbf{Case 3:}
\begin{multline*}
P(\{X_2,X_3\} \subset \NPE^r(X_1)\cap \G_1^r(X_1),\; X_1 \in T_s)=\\
\left(\int_{s_5}^{s_9}\int_{r_7(x)}^{r_3(x)} +
\int_{s_9}^{s_{12}}\int_{0}^{r_3(x)} +
\int_{s_{12}}^{1/2}\int_{0}^{r_6(x)}\right)
\frac{A(\msP(G_1,G_2,Q_1,P_2,M_3,G_6))^2}{A(\TY)^3}dydx=\\
-\Bigl[17496\,r^{19}-122472\,r^{18}+139968\,r^{17}+524880\,r^{16}-553095\,r^{15}-
595971\,r^{14}+368826\,r^{13}-724758\,r^{12}-543876\,r^{11}+\\
1416996\,r^{10}+1646470\,r^9+92870\,r^8+523048\,r^7-768368\,r^6-1729902\,r^5-1434990\,r^4+122185\,r^3+ 941941\,r^2+\\
573440\,r+114688\Bigr] \Big/ \Bigl[2099520\, \left( r+1 \right)^5r^{10}\Bigr]
\end{multline*}
where
$A(\msP(G_1,G_2,Q_1,P_2,M_3,G_6))=
-\Bigl[\sqrt{3} \bigl( 4\,\sqrt{3}r^2y-8\,r^3\sqrt{3}y+4\,r^2 y^2+4\,r^4 y^2+
4\,y^2+8\,r^4\sqrt{3}y\,x+6-12\,x^2 r^2-12\,x-12\,r-24\,r^3 x+
12\,r^4 x^2+9\,r^2+12\,r\,x-4\,\sqrt{3}r\,y+12\,x^2+
4\,\sqrt{3}y+12\,r^2 x \bigr)\Bigr]\Big/\Bigl[24\,r^2\Bigr]
$.

\noindent \textbf{Case 4:}
\begin{multline*}
P(\{X_2,X_3\} \subset \NPE^r(X_1)\cap \G_1^r(X_1),\; X_1 \in T_s)=\\
\left( \int_{s_8}^{s_5}\int_{r_8(x)}^{r_2(x)} +
\int_{s_5}^{s_{10}}\int_{r_3(x)}^{r_2(x)} +
\int_{s_{10}}^{s_{12}}\int_{r_3(x)}^{r_6(x)} \right)
\frac{A(\msP(G_1,M_1,L_2,Q_1,P_2,M_3,G_6))^2}{A(\TY)^3}dydx=\\
-\Bigl[32768-409264128\,r^7+1455989508\,r^{12}+680709729\,r^8-4423680\,r^3+155509\,r^2+
22889801\,r^4+202936917\,r^6+\\
6011901\,r^{20}+1060982949\,r^{16}-614739456\,r^{17}+240330993\,r^{18}-
56097792\,r^{19}-77783040\,r^5-999857664\,r^9+\\
1299257316\,r^{10}-1461851136\,r^{11}-1407624192\,r^{13}+
1414729905\,r^{14}-1352392704\,r^{15}\Bigr]\Big/\Bigl[2099520\, \left( r^2+1 \right)^5r^{10}\Bigr]
\end{multline*}
where
$A(\msP(G_1,M_1,L_2,Q_1,P_2,M_3,G_6))=
-\Bigl[\sqrt{3} \bigl( -6\,x^2r^2-3+6\,x-12\,r^3 x+6\,r^4 x^2-
4\,r^3\sqrt{3} y+4\,\sqrt{3}y\,x+4\,r^4\sqrt{3} y\,x+2\,r^4y^2+
3\,r^2+2\,\sqrt{3}r^2 y-2\,\sqrt{3} y+2\,r^2 y^2+6\,r^2 x \bigr)\Bigr]\Big/\Bigl[12 \, r^2\Bigr]
$.

\noindent \textbf{Case 5:}
\begin{multline*}
P(\{X_2,X_3\} \subset \NPE^r(X_1)\cap \G_1^r(X_1),\; X_1 \in T_s)=
\left(\int_{s_3}^{s_8}\int_{r_5(x)}^{r_2(x)} +
\int_{s_8}^{s_5}\int_{r_5(x)}^{r_8(x)}\right)
\frac{A(\msP(G_1,M_1,P_1,P_2,M_3,G_6))^2}{A(\TY)^3}dydx=\\
\Bigl[\bigl( 35361\,r^{16}-229392\,r^{15}+602820\,r^{14}-858384\,r^{13}+778848\,r^{12}-
460368\,r^{11}+277740\,r^{10}-258768\,r^9+160594\,r^8-62256\,r^7-\\
5892\,r^6-17712\,r^5+19224\,r^4+11664\,r^3+5076\,r^2+1296\,r+405 \bigr)
\left( -12\,r+7\,r^2+3 \right)^2\Bigr]\Big/\Bigl[699840\,r^{10} \left( r^2+1 \right)^5\Bigr]
\end{multline*}
where
$A(\msP(G_1,M_1,P_1,P_2,M_3,G_6))=
-{\frac {\sqrt{3} \left( -4\,r^3\sqrt{3}y-12\,r^3x+3\,r^2+6\,r^4\sqrt{3}y\,x+
9\,r^4x^2+3\,r^4y^2+y^2+2\,\sqrt{3}y\,x+3\,x^2 \right) }{12\,r^2}}
$.

\noindent \textbf{Case 6:}
\begin{multline*}
P(\{X_2,X_3\} \subset \NPE^r(X_1)\cap \G_1^r(X_1),\; X_1 \in T_s)=\\
\left(\int_{s_2}^{s_3}\int_{r_5(x)}^{\ell_{am}(x)} +
\int_{s_3}^{s_7}\int_{r_2(x)}^{\ell_{am}(x)} +
\int_{s_7}^{s_8}\int_{r_2(x)}^{r_8(x)} \right)
\frac{A(\msP(G_1,M_1,P_1,P_2,M_3,G_6))^2}{A(\TY)^3}dydx=\\
-\Bigl[3645-17496\,r+5003898912\,r^{28}+31646646384\,r^{26}+110098944\,r^{30}-1090803456\,r^{29}-
14630751360\,r^{27}+66339\,r^2-\\
99072645696\,r^{23}+79269457632\,r^{24}+66073158\,r^8-
4870743552\,r^{13}-168073488\,r^9+535086\,r^4-262440\,r^3-1737936\,r^5-\\
18592416\,r^7-107383563504\,r^{21}-41219053272\,r^{17}+58981892347\,r^{18}-78265758888\,r^{19}+95887286866\,r^{20}+\\
109053166552\,r^{22}+5500548\,r^6+466565130\,r^{10}-1070573040\,r^{11}+
2380992104\,r^{12}+9191633420\,r^{14}-16312513248\,r^{15}+\\
26801184917\,r^{16}-54759787776\,r^{25}\Bigr]\Big/\Bigl[1399680\, \left( r^2+1 \right)^5 \left( 2\,r^2+1 \right)
^5r^{10}\Bigr]
\end{multline*}
where
$A(\msP(G_1,M_1,P_1,P_2,M_3,G_6))=
-{\frac {\sqrt{3} \left( -4\,r^3\sqrt{3}y-12\,r^3 x+
3\,r^2+6\,r^4\sqrt{3}y\,x+9\,r^4 x^2+3\,r^4 y^2+y^2+
2\,\sqrt{3}y\,x+3\,x^2 \right) }{12\,r^2}}
$.

\noindent \textbf{Case 7:}
\begin{multline*}
P(\{X_2,X_3\} \subset \NPE^r(X_1)\cap \G_1^r(X_1),\; X_1 \in T_s)=\\
\left(\int_{s_7}^{s_8}\int_{r_8(x)}^{r_9(x)}+
\int_{s_8}^{s_{10}}\int_{r_2(x)}^{r_9(x)} \right)
\frac{A(\msP(G_1,M_1,L_2,Q_1,P_2,M_3,G_6))^2}{A(\TY)^3}dydx=\\
\Bigl[4\,\bigl( 162576\,r^{22}-1083456\,r^
{21}+3368016\,r^{20}-6969888\,r^{19}+11578088\,r^{18}-
15664080\,r^{17}+18796852\,r^{16}-19984824\,r^{15}+\\
19534445\,r^{14}-18170472\,r^{13}+15507752\,r^{12}-
13150464\,r^{11}+9987958\,r^{10}-7448736\,r^9+5016464\,r^8-2991768\,r^7+1857485\,r^6-\\
749160\,r^5+481804\,r^4-96720\,r^3+76160\,r^2-4032\,r+4320 \bigr)
\left( 2\,r -1 \right)^2 \left( r -1 \right)^2\Bigr]\Big/\Bigl[32805\, \left( r^2+1 \right)^5r^6 \left( 2\,r^2+1 \right)^5\Bigr]
\end{multline*}
where
$A(\msP(G_1,M_1,L_2,Q_1,P_2,M_3,G_6))=
-\Bigl[\sqrt{3} \bigl( -6\,x^2r^2-3+6\,x-12\,r^3x+6\,r^4x^2-
4\,r^3\sqrt{3}y+4\,\sqrt{3}y\,x+4\,r^4\sqrt{3}y\,x+2\,r^4y^2+
3\,r^2+2\,\sqrt{3}r^2 y-2\,\sqrt{3}y+2\,r^2 y^2+6\,r^2 x \bigr)\Bigr]\Big/\Bigl[12\,r^2\Bigr]
$.

\noindent \textbf{Case 8:}
\begin{multline*}
P(\{X_2,X_3\} \subset \NPE^r(X_1)\cap \G_1^r(X_1),\; X_1 \in T_s)=\\
\left(\int_{s_{12}}^{s_{13}}\int_{r_6(x)}^{r_3(x)} +
\int_{s_{13}}^{1/2}\int_{r_6(x)}^{r_2(x)} \right)
\frac{A(\msP(G_1,G_2,Q_1,N_3,M_C,M_3,G_6))^2}{A(\TY)^3}dydx=\\
-\Bigl[-458752+811008\,r^2+329205504\,r^8-582626304\,r^{13}-489563136\,r^9-65536\,r^4-
168708096\,r^7-57883680\,r^{17}+\\
18009258\,r^{18}-3623400\,r^{19}+352563\,r^{20}+41502720\,r^6+659111904\,r^{10}-
761846400\,r^{11}+725173376\,r^{12}+409477188\,r^{14}-\\
254829600\,r^{15}+135968852\,r^{16}\Bigr]\Big/\Bigl[8398080\,r^{10}\Bigr]
\end{multline*}
where
$A(\msP(G_1,G_2,Q_1,N_3,M_C,M_3,G_6))=
-\Bigl[\sqrt{3} \bigl( -12\,x^2r^2-12\,x-12\,r-12\,r^3x+3\,r^4x^2+
4\,\sqrt{3}r^2 y+5\,r^2+12\,r\,x+12\,x^2+2\,r^4\sqrt{3}y\,x+4\,r^2y^2-
4\,r^3\sqrt{3}y+6+4\,y^2+r^4 y^2+4\,\sqrt{3}y+12\,r^2 x-4\,\sqrt{3}r\,y \bigr)\Bigr]\Big/\Bigl[24\,r^2\Bigr]
$.

\noindent \textbf{Case 9:}
\begin{multline*}
P(\{X_2,X_3\} \subset \NPE^r(X_1)\cap \G_1^r(X_1),\; X_1 \in T_s)=\\
\left(\int_{s_{10}}^{s_{12}}\int_{r_6(x)}^{r_2(x)} +
\int_{s_{12}}^{s_{13}}\int_{r_3(x)}^{r_2(x)} \right)
\frac{A(\msP(G_1,M_1,L_2,Q_1,N_3,M_C,M_3,G_6))^2}{A(\TY)^3}dydx=
\Bigl[\bigl( 7203\,r^{16}-49392\,r^{15}+\\
170226\,r^{14}-392112\,r^{13}+680784\,r^{12}-1040256\,r^{11}+
1385628\,r^{10}-1337760\,r^9+816224\,r^8-253824\,r^7+\\
469088\,r^6-1029888\,r^5+820992\,r^4-488448\,r^3+190976\,r^2+49152\,r+8192 \bigr)
 \left( -12\,r+7\,r^2+4 \right)^2\Bigr]\Big/\Bigl[8398080\,r^{10}\Bigr]
\end{multline*}
where
$A(\msP(G_1,M_1,L_2,Q_1,N_3,M_C,M_3,G_6))=
-\Bigl[\sqrt{3} \bigl( -12\,x^2r^2-6+12\,x-12\,r^3x+3\,r^4x^2+2\,r^2+
2\,r^4\sqrt{3}y\,x+r^4y^2+8\,\sqrt{3}y\,x+4\,r^2 y^2-
4\,\sqrt{3} y+4\,\sqrt{3}r^2y+12\,r^2x-4\,r^3\sqrt{3}y \bigr)\Bigr]\Big/\Bigl[24\,r^2\Bigr]
$.

\noindent \textbf{Case 10:}
\begin{multline*}
P(\{X_2,X_3\} \subset \NPE^r(X_1)\cap \G_1^r(X_1),\; X_1 \in T_s)=\\
\left(\int_{s_{10}}^{s_{14}}\int_{r_2(x)}^{r_{10}(x)} +
\int_{s_{14}}^{s_{13}}\int_{r_2(x)}^{r_{12}(x)}+
\int_{s_{13}}^{1/2}\int_{r_3(x)}^{r_{12}(x)} \right)
\frac{A(\msP(G_1,M_1,L_2,Q_1,N_3,L_4,L_5,M_3,G_6))^2}{A(\TY)^3}dydx=\\
\Bigl[4423680-4627454976\,r^6+511684992\,r^{11}+2163142656\,r^7-660127744\,r^2-31555584\,r+
3534520320\,r^3+7647989760\,r^5+\\
7785504\,r^{15}-1313880\,r^{16}+19683\,r^{18}-7240624128\,r^4-1511047552\,r^8+
1204122240\,r^9-796453824\,r^{10}-282583320\,r^{12}+\\
107804736\,r^{13}-30362052\,r^{14}\Bigr]\Big/\Bigl[16796160\,r^6\Bigr]
\end{multline*}
where
$A(\msP(G_1,M_1,L_2,Q_1,N_3,L_4,L_5,M_3,G_6))=
-\Bigl[\sqrt{3} \bigl( -16\,\sqrt{3}r\,y+20\,\sqrt{3}y-24\,y^2-12\,x^2 r^2+
12\,x+24\,r-12\,r^3 x+3\,r^4 x^2-6\,r^2-24+4\,\sqrt{3}r^2 y+
8\,\sqrt{3}y\,x-4\,r^3\sqrt{3}y+4\,r^2y^2+r^4 y^2+2\,r^4\sqrt{3}y\,x+
12\,r^2 x \bigr) \Bigr]\Big/\Bigl[24\,r^2\Bigr]
$.

\noindent \textbf{Case 11:}
\begin{multline*}
P(\{X_2,X_3\} \subset \NPE^r(X_1)\cap \G_1^r(X_1),\; X_1 \in T_s)=\\
\left(\int_{s_7}^{s_{11}}\int_{r_9(x)}^{\ell_{am}(x)} +
\int_{s_{11}}^{s_{10}}\int_{r_9(x)}^{r_{12}(x)} +
\int_{s_{10}}^{s_{14}}\int_{r_{10}(x)}^{r_{12}(x)} \right
)\frac{A(\msP(G_1,M_1,L_2,Q_1,Q_2,L_5,M_3,G_6))^2}{A(\TY)^3}dydx=\\
-\Bigl[ \left( r -1 \right)  \bigl( -1474560+
8847360\,r+111456\,r^{26}+111456\,r^{27}-27738112\,r^2
+23311152\,r^{23}-167184\,r^{24}-808889416\,r^8-\\
2228253688\,r^{13}+366739256\,r^9-207619072\,r^4+
98557952\,r^3+397199360\,r^5+802401664\,r^7-34733448\,r^{21}-624736557\,r^{17}+\\
400615470\,r^{18}-134938386\,r^{19}+39014136\,r^{20}-18026064\,r^{22}-640058432\,r^6+
407655352\,r^{10}-1227078728\,r^{11}+\\
1996721576\,r^{12}+2033409092\,r^{14}-1681870468\,r^{15}+1064030499\,r^{16}-
2842128\,r^{25} \bigr) \Bigr]\Big/\Bigl[1866240\, \left( 2\,r^2+1 \right)^5r^6\Bigr]
\end{multline*}
where
$A(\msP(G_1,M_1,L_2,Q_1,Q_2,L_5,M_3,G_6))=
-\Bigl[\sqrt{3} \bigl( 4\,\sqrt{3}r^2y+4\,\sqrt{3}y\,x-2\,r^2 y^2-
4\,r^3\sqrt{3}y-4\,y^2-4\,\sqrt{3}r^2y\,x-6\,x^2 r^2+6\,x-12\,r^3x+
3\,r^4 x^2+3\,r^2+2\,r^4\sqrt{3}y\,x+r^4 y^2+2\,\sqrt{3}y+12\,r^2 x-6 \bigr) \Bigr]\Big/\Bigl[12\,r^2\Bigr]
$.

\noindent \textbf{Case 12:}
\begin{multline*}
P(\{X_2,X_3\} \subset \NPE^r(X_1)\cap \G_1^r(X_1),\; X_1 \in T_s)=
\int_{s_{13}}^{1/2}\int_{r_2(x)}^{r_3(x)} \frac{A(\msP(G_1,G_2,Q_1,N_3,L_4,L_5,M_3,G_6))^2}{A(\TY)^3}dydx=\\
\Bigl[ \bigl( 2322432-7554816\,r+9510912\,{
r}^2+1046068\,r^8-558720\,r^9+2444224\,r^4-5799360\,r^3-2134656\,r^5-1608672\,r^7+\\
2169696\,r^6+216300\,r^{10}-55440\,r^{11}+7095\,r^{12} \bigr)
 \left( -6+5\,r \right)^2\Bigr]\Big/\Bigl[4199040\,r^4\Bigr]
\end{multline*}
where
$A(\msP(G_1,G_2,Q_1,N_3,L_4,L_5,M_3,G_6))=
-\Bigl[\sqrt{3} \bigl( -12\,x^2r^2-12\,x+12\,r-12\,r^3x+3\,r^4x^2-3\,r^2+
12\,r\,x+28\,\sqrt{3}y+12\,x^2-20\, y^2+12\,r^2 x+r^4y^2+4\,r^2 y^2-
4\,r^3\sqrt{3}y+2\,r^4\sqrt{3} y\,x+4\,\sqrt{3}r^2 y-20\,\sqrt{3}r\, y-12 \bigr) \Bigr]\Big/\Bigl[24\,r^2\Bigr]
$.

\noindent \textbf{Case 13:}
\begin{multline*}
P(\{X_2,X_3\} \subset \NPE^r(X_1)\cap \G_1^r(X_1),\; X_1 \in T_s)=
\int_{s_{14}}^{1/2}\int_{r_{12}(x)}^{r_{10}(x)} \frac{A(\msP(L_1,L_2,Q_1,N_3,L_4,L_5,L_6))^2}{A(\TY)^3}dydx=\\
-\Bigl[\bigl( 9\,r^{14}+36\,r^{13}-132\,r^{12}-576\,r^{11}+164\,r^{10}+2512\,r^9+
4976\,r^8-1536\,r^7-13888\,r^6-17536\,r^5-3072\,r^4+79360\,r^3+\\
9216\,r^2-120832\,r+61440 \bigr)  \left( -2+r \right)
\left( r^2+2\,r -4 \right)^2\Bigr]\Big/\Bigl[7680\, \left( r+2 \right)^3r^4\Bigr]
\end{multline*}
where
$A(\msP(L_1,L_2,Q_1,N_3,L_4,L_5,L_6))=
-\Bigl[\sqrt{3} \bigl( r^4 y^2-8\,\sqrt{3} r\,y-8\,\sqrt{3}y\,x+4\,r^2y^2-
4\,r^3\sqrt{3}y-32\,y^2+2\,r^4\sqrt{3}y\,x-12\,x^2 r^2+12\,x+
24\,r-12\,r^3 x+3\,r^4 x^2-12\,r^2+4\,\sqrt{3}r^2 y+24\,r\,x-24\,x^2-24+
20\,\sqrt{3}y+12\,r^2 x \bigr) \Bigr]\Big/\Bigl[24\,r^2\Bigr]
$.

\noindent \textbf{Case 14:}
\begin{multline*}
P(\{X_2,X_3\} \subset \NPE^r(X_1)\cap \G_1^r(X_1),\; X_1 \in T_s)=\\
\left(\int_{s_{11}}^{s_{14}}\int_{r_{12}(x)}^{\ell_{am}(x)} +
\int_{s_{14}}^{1/2}\int_{r_{10}(x)}^{\ell_{am}(x)}\right)
\frac{A(\msP(L_1,L_2,Q_1,Q_2,L_5,L_6))^2}{A(\TY)^3}dydx=\\
\Bigl[\left( r -1 \right)  \bigl( 3483\,r^{18}+24381\,r^{17}-34830\,r^{16}-529416\,r^{15}-
265680\,r^{14}+4274208\,r^{13}+4999320\,r^{12}-15227352\,r^{11}-\\
25751336\,r^{10}+19466488\,r^9+62834064\,r^8+17452256\,r^7-53339200\,r^6-117114624\,r^5-
51206656\,r^4+270430208\,r^3+\\
58073088\,r^2-296222720\,r+122159104 \bigr) \Bigr]\Big/\Bigl[ 1866240\, \left( r+2 \right)^3r^4\Bigr]
\end{multline*}
where
$A(\msP(L_1,L_2,Q_1,Q_2,L_5,L_6))=
-\Bigl[\sqrt{3} \bigl( -4\,\sqrt{3}y\,x-2\,r^2 y^2+4\,\sqrt{3}r\,y-4\,r^3\sqrt{3}y-
8\,y^2-4\,\sqrt{3}r^2 y\,x-6\,x^2 r^2+6\,x-12\,r^3 x+3\,r^4 x^2+
4\,\sqrt{3}r^2y+12\,r\,x-12\, x^2+2\,r^4\sqrt{3}y\,x+r^4 y^2+
2\,\sqrt{3}y+12\,r^2x-6 \bigr) \Bigr]\Big/\Bigl[12\,r^2\Bigr]
$.
}

Adding up the $P(\{X_2,X_3\} \subset \NPE^r(X_1)\cap \G_1^r(X_1),\; X_1 \in T_s)$
values in the 14 possible cases above, and multiplying by 6
we get for $r \in \bigl[6/5,\sqrt{5}-1\bigr)$,
\begin{multline*}
\nu_\la(r)=
-\Bigl[219936\,r -3041936\,r^2-30889822
\,r^8+18084672\,r^{13}+27137438\,r^9+2364868\,r^4+2305864\,r^3-4168820\,r^5-\\
2832544\,r^7+486\,r^{21}-118850\,r^{17}-45155\,r^{18}-269\,r^{19}+3402\,r^{20}+
11101160\,r^6+24604048\,r^{10}-43009544\,r^{11}+8770788\,r^{12}-\\
13736295\,r^{14}+2751855\,r^{15}+443518\,r^{16}+49152\Bigr]
\Big/\Bigl[116640\,r^6 \left( r+2 \right)^2 \left( 2\,r^{
2}+1 \right)  \left( r+1 \right)^3\Bigr].
\end{multline*}
The $\nu_\la(r)$ values for the other intervals can be calculated similarly.

\section*{Appendix 3.2: Derivation of $\mu_\lo(r)$ and $\nu_\lo(r)$ for Uniform Data}
\subsubsection*{Derivation of $\mu_\lo(r)$ in Theorem \ref{thm:asy-norm-under}}
First we find $\mu_\lo(r)$ for $r \in \Bigl[1,\infty)$.
Observe that, by symmetry,
$$\mu_\lo(r)=P\bigl( X_2 \in \NPE^r(X_1) \cup \G_1^r(X_1) \bigr)=
6\,P\bigl( X_2 \in \NY^r(X_1) \cup \G_1^r(X_1), X_1 \in T_s \bigr).$$

For $r \in [1,4/3)$,
there are 17 cases to consider for calculation of $\nu_{\lo}(r)$ in the underlying graph case.
Each Case $j$ correspond to $R_i$ for $i=1,2,\ldots,17$ in Figure \ref{fig:cases-AND-OR}.
{\small
\noindent \textbf{Case 1:}
\begin{multline*}
P(X_2 \in \NPE^r(X_1) \cup \G_1^r(X_1), X_1 \in T_s)=
\left(\int_{0}^{s_0}\int_{0}^{\ell_{am}(x)} +
\int_{s_0}^{s_1}\int_{r_1(x)}^{\ell_{am}(x)} \right)
\frac{A(\msP(A,M_1,M_C,M_3))}{A(\TY)^2}dydx=\\
{\frac{4}{27}}\,r^2-4\,r/9+1/3
\end{multline*}
where
$A(\msP(A,M_1,M_C,M_3))=\sqrt{3}/12
$.

\noindent \textbf{Case 2:}
\begin{multline*}
P(X_2 \in \NPE^r(X_1) \cup \G_1^r(X_1), X_1 \in T_s)=\\
\left(\int_{s_0}^{s_1}\int_{0}^{r_1(x)}+
\int_{s_1}^{s_3}\int_{0}^{r_2(x)}+
\int_{s_3}^{s_4}\int_{0}^{r_5(x)} +
\int_{s_4}^{s_5}\int_{r_3(x)}^{r_5(x)} \right)
\frac{A(\msP(A,M_1,L_2,L_3,M_C,M_3))}{A(\TY)^2}dydx=\\
-{\frac{ \left( r -1 \right)  \left( 1817\,r^7-7807\,r^6+14157\,r^5-14067\,r^4+7893\,r^3-
2475\,r^2+405\,r -27 \right) }{864\,r^6}}
\end{multline*}
where
$A(\msP(A,M_1,L_2,L_3,M_C,M_3))=
{\frac{\sqrt{3} \left( -4\,\sqrt{3}r\, y-12\,r+12\,r\, x+
5\,r^2+3\,y^2+6\,\sqrt{3} y-6\,\sqrt{3} y\, x+9-18\, x+9\, x^2 \right) }
{12\,r^2}}
$.

\noindent \textbf{Case 3:}
\begin{multline*}
P(X_2 \in \NPE^r(X_1) \cup \G_1^r(X_1), X_1 \in T_s)=
\left(\int_{s_4}^{s_5}\int_{0}^{r_3(x)} +
\int_{s_5}^{s_6}\int_{0}^{r_5(x)} \right)
\frac{A(\msP(A,G_2,G_3,M_2,M_C,M_3))}{A(\TY)^2}dydx=\\
{\frac{ \left( 13\,r^4-4\,r^3+4\,r -1-2\,r^2 \right)  \left( r -1 \right)^4}{96\,r^6}}
\end{multline*}
where
$A(\msP(A,G_2,G_3,M_2,M_C,M_3))=
-{\frac{\sqrt{3} \left( \,y^2+2\,\sqrt{3} y-2\,\sqrt{3} y\, x+
3-6\, x+3\, x^2-2\,r^2 \right) }{12\,r^2}}
$.

\noindent \textbf{Case 4:}
\begin{multline*}
P(X_2 \in \NPE^r(X_1) \cup \G_1^r(X_1), X_1 \in T_s)=
\left(\int_{s_1}^{s_2}\int_{r_2(x)}^{\ell_{am}(x)} +
\int_{s_2}^{s_3}\int_{r_2(x)}^{r_5(x)} \right)
\frac{A(\msP(A,M_1,L_2,L_3,L_4,L_5,M_3))}{A(\TY)^2}dydx=\\
{\frac{ \left( 9-72\,r+192\,r^2-192\,r^3+76\,r^4 \right)
\left( 4\,r -3+\sqrt{3} \right)^2 \left( 4\,r -3-\sqrt{3} \right)^2}{10368\,r^6}}
\end{multline*}
where
$A(\msP(A,M_1,L_2,L_3,L_4,L_5,M_3))=
{\frac{\sqrt{3} \left( 4\,\sqrt{3}r\, y+9\,r^2-24\,r+12\,r\, x+15\,y^2-6\,\sqrt{3} y-
6\,\sqrt{3} y\, x+18-18\, x+9\, x^2 \right)}{12\,r^2}}
$.

\noindent \textbf{Case 5:}
\begin{multline*}
P(X_2 \in \NPE^r(X_1) \cup \G_1^r(X_1), X_1 \in T_s)=
\left(\int_{s_5}^{s_6}\int_{r_5(x)}^{r_7(x)} +
\int_{s_6}^{s_9}\int_{0}^{r_7(x)} \right)
\frac{A(\msP(A,G_2,G_3,M_2,M_C,P_2,N_2))}{A(\TY)^2}dydx=\\
{\frac{ \left( -1+2\,r+6\,r^2-6\,r^3+
22\,r^5+17\,r^6 \right)  \left( r -1 \right)^3}{96\,r^6 \left( r+1 \right)^3}}
\end{multline*}
where
$A(\msP(A,G_2,G_3,M_2,M_C,P_2,N_2))=
\Bigl[\sqrt{3} \bigl( -2\,y^2-4\,\sqrt{3} y+4\,\sqrt{3} y\, x-6+12\, x-6\, x^2+
7\,r^2-4\,r^3\sqrt{3} y-12\,r^3 x+8\,r^4\sqrt{3} y\, x+12\,r^4 x^2+
4\,r^4\,y^2 \bigr) \Bigr]\Big/\Bigl[24\,r^2\Bigr]
$.

\noindent \textbf{Case 6:}
\begin{multline*}
P(X_2 \in \NPE^r(X_1) \cup \G_1^r(X_1), X_1 \in T_s)=\\
\left(\int_{s_5}^{s_9}\int_{r_7(x)}^{r_3(x)} +
\int_{s_9}^{s_{12}}\int_{0}^{r_3(x)} +
\int_{s_{12}}^{1/2}\int_{0}^{r_6(x)} \right)
\frac{A(\msP(A,N_1,Q_1,G_3,M_2,M_C,P_2,N_2))}{A(\TY)^2}dydx=\\
-{\frac{81\,r^9-189\,r^8+561\,r^7-45\,r^6-1894\,r^5-18\,r^4+1912\,r^3+
224\,r^2-384\,r -128}{1296\, \left( r+1 \right)^3 r^4}}
\end{multline*}
where
$A(\msP(A,N_1,Q_1,G_3,M_2,M_C,P_2,N_2))=
\Bigl[\sqrt{3} \bigl( 4\,r\,y^2-4\,\sqrt{3} y+12\, x+13\,r -12+18\,r^3 x^2+
12\,r\, x-12\,r\, x^2-8\,\sqrt{3}r^2 y+4\,\sqrt{3}r\, y-24\,r^2 x+
12\,\sqrt{3}r^3 y\, x+6\,r^3\,y^2 \bigr) \Bigr]\Big/\Bigl[24\,r\Bigr]
$.

\noindent \textbf{Case 7:}
\begin{multline*}
P(X_2 \in \NPE^r(X_1) \cup \G_1^r(X_1), X_1 \in T_s)=\\
\left(\int_{s_8}^{s_5}\int_{r_8(x)}^{r_2(x)} +
\int_{s_5}^{s_{10}}\int_{r_3(x)}^{r_2(x)} +
\int_{s_{10}}^{s_{12}}\int_{r_3(x)}^{r_6(x)} \right)
\frac{A(\msP(A,N_1,Q_1,L_3,M_C,P_2,N_2))}{A(\TY)^2}dydx=\\
-\Bigl[128-1536\,r -302592\,r^7+11753\,r^{12}+346171\,r^8-28416\,r^3+
8384\,r^2+69760\,r^4+220201\,r^6-135936\,r^5-305664\,r^9+\\
186683\,r^{10}-69120\,r^{11}\Bigr]\Big/\Bigl[1944\, \left( r^2+1 \right)^3r^6\Bigr]
\end{multline*}
where
$A(\msP(A,N_1,Q_1,L_3,M_C,P_2,N_2))=
\Bigl[\sqrt{3} \bigl( -4\,\sqrt{3}r\, y+2\,\sqrt{3}
r^2 y-12\, x-12\,r+8\,r^2+12\,r\, x-6\, x^2r^2+2\,r^2\,y^2-
4\,\sqrt{3} y\, x+3\,r^4\,y^2-4\,r^3\sqrt{3} y-12\,r^3 x+9\,r^4 x^2+
4\,\sqrt{3} y+6\,r^4\sqrt{3} y\, x+6\, x^2+6\,r^2 x+6+2\,y^2 \bigr) \Bigr]\Big/\Bigl[12\,r^2\Bigr]
$.

\noindent \textbf{Case 8:}
\begin{multline*}
P(X_2 \in \NPE^r(X_1) \cup \G_1^r(X_1), X_1 \in T_s)=
\left(\int_{s_3}^{s_8}\int_{r_5(x)}^{r_2(x)} +
\int_{s_8}^{s_5}\int_{r_5(x)}^{r_8(x)} \right)
\frac{A(\msP(A,N_1,P_1,L_2,L_3,M_C,P_2,N_2))}{A(\TY)^2}dydx=\\
{\frac{ \left( 895\,r^8-2472\,r^7+3363\,r^6-2880\,r^5+2220\,r^4-1296\,r^3+675\,r^2-
216\,r+27 \right)  \left( -12\,r+7\,r^2+3 \right)^2}{7776\, \left( r^2+1 \right)^3r^6}}
\end{multline*}
where
$A(\msP(A,N_1,P_1,L_2,L_3,M_C,P_2,N_2))=
\Bigl[\sqrt{3} \bigl( 4\,r^4\,y^2+8\,r^4 \sqrt{3} y\, x+12\,r^4 x^2-
4\,r^3 \sqrt{3} y-12\,r^3 x-4\,\sqrt{3}r\, y-12\,r+12\,r\, x+8\,r^2+
3\,y^2+6\,\sqrt{3}y-6\,\sqrt{3} y\, x+9-18\, x+9\, x^2 \bigr) \Bigr]\Big/\Bigl[ 12\, r^2\Bigr]
$.

\noindent \textbf{Case 9:}
\begin{multline*}
P(X_2 \in \NPE^r(X_1) \cup \G_1^r(X_1), X_1 \in T_s)=\\
\left(\int_{s_2}^{s_3}\int_{r_5(x)}^{\ell_{am}(x)} +
\int_{s_3}^{s_7}\int_{r_2(x)}^{\ell_{am}(x)} +
\int_{s_7}^{s_8}\int_{r_2(x)}^{r_8(x)} \right)
\frac{A(\msP(A,N_1,P_1,L_2,L_3,L_4,L_5,P_2,N_2))}{A(\TY)^2}dydx=\\
-\Bigl[355328\,r^{18}-2204160\,r^{17}+
6591792\,r^{16}-13254912\,r^{15}+20639832\,r^{14}-26417664
\,r^{13}+28578916\,r^{12}-26760576\,r^{11}+\\
21960774\,r^{10}-15877152\,r^9+10180620\,r^8-5753232\,r^7+
2856483\,r^6-1222128\,r^5+438777\,r^4-128304\,r^3+28107\,r^2-\\
3888\,r+243\Bigr]\Big/\Bigl[7776\, \left( r^2+1 \right)^3 \left( 2\,r^2+1 \right)^3r^6\Bigr]
\end{multline*}
where
$A(\msP(A,N_1,P_1,L_2,L_3,L_4,L_5,P_2,N_2))=
\Bigl[\sqrt{3} \bigl( 18+4\,\sqrt{3}r\, y-18\,x-24\,r+12\,r^2+12\,r\, x-6\,\sqrt{3} y+
8\,r^4\sqrt{3} y\, x-12\,r^3 x+12\,r^4 x^2+9\, x^2+15\,y^2+4\,r^4 y^2-
4\,r^3\sqrt{3} y-6\,\sqrt{3} y\,x \bigr) \Bigr]\Big/\Bigl[12\,r^2\Bigr]
$.

\noindent \textbf{Case 10:}
\begin{multline*}
P(X_2 \in \NPE^r(X_1) \cup \G_1^r(X_1), X_1 \in T_s)=
\left(\int_{s_7}^{s_8}\int_{r_8(x)}^{r_9(x)}+
\int_{s_8}^{s_{10}}\int_{r_2(x)}^{r_9(x)} \right)
\frac{A(\msP(A,N_1,Q_1,L_3,L_4,L_5,P_2,N_2))}{A(\TY)^2}dydx=\\
\Bigl[ 8\bigl( 288\,r^{12}-864\,r^{11}+1486\,r^{10}-1896\,r^9+2056\,r^8-1608\,r^7+1189\,r^6-
654\,r^5+317\,r^4-132\,r^3+44\,r^2-12\,r+2 \bigr)  \left( 2\,r -1 \right)^2 \\
\left( r -1 \right)^2\Bigr]\Big/\Bigl[243\, \left( r^2+1 \right)^3 \left( 2\,r^2+1 \right)^3 r^4\Bigr]
\end{multline*}
where
$A(\msP(A,N_1,Q_1,L_3,L_4,L_5,P_2,N_2))=
\Bigl[\sqrt{3} \bigl( 4\,\sqrt{3}r\, y+2\,\sqrt{3}{
r}^2 y-8\,\sqrt{3} y-12\, x-24\,r+12\,r^2+12\,r\, x-6\, x^2r^2+15-
12\,r^3 x+9\,r^4 x^2+6\, x^2+6\,r^2 x+6\,r^4\sqrt{3} y\, x+2\,r^2\,y^2-
4\,\sqrt{3} y\, x+3\,r^4\,y^2-4\,r^3\sqrt{3} y+14\,y^2 \bigr) \Bigr]\Big/\Bigl[12\,r^2\Bigr]
$.

\noindent \textbf{Case 11:}
\begin{multline*}
P(X_2 \in \NPE^r(X_1) \cup \G_1^r(X_1), X_1 \in T_s)=
\left(\int_{s_{12}}^{s_{13}}\int_{r_6(x)}^{r_3(x)} +
\int_{s_{13}}^{1/2}\int_{r_6(x)}^{r_2(x)} \right)
\frac{A(\msP(A,N_1,Q_1,G_3,M_2,N_3,N_2))}{A(\TY)^2}dydx=\\
-{\frac{1536-6528\,r^2+133834\,r^8-48240\,r^9+95616\,r^4-20736\,r^3-158976\,r^5-
200064\,r^7+196680\,r^6+7107\,r^{10}}{15552\,r^4}}
\end{multline*}
where
$A(\msP(A,N_1,Q_1,G_3,M_2,N_3,N_2))=
\Bigl[\sqrt{3} \bigl( 4\,r\,y^2+12\, x+9\,r -12+9\,r^3 x^2+12\,r\, x-
12\,r\,x^2-4\,\sqrt{3}r^2 y+4\,\sqrt{3}r\, y+6\,\sqrt{3}r^3 y\, x+
3\,r^3\,y^2-12\,r^2 x-4\,\sqrt{3} y \bigr) \Bigr]\Big/\Bigl[24\,r\Bigr]
$.

\noindent \textbf{Case 12:}
\begin{multline*}
P(X_2 \in \NPE^r(X_1) \cup \G_1^r(X_1), X_1 \in T_s)=
\left(\int_{s_{10}}^{s_{13}}\int_{r_6(x)}^{r_2(x)} +
\int_{s_{12}}^{s_{13}}\int_{r_3(x)}^{r_2(x)} \right)
\frac{A(\msP(A,N_1,Q_1,L_3,N_3,N_2))}{A(\TY)^2}dydx=\\
{\frac{ \left( 147\,r^8-504\,r^7+530
\,r^6-336\,r^5+876\,r^4-1056\,r^3+896\,r^2-384\,r+64 \right)  \left( -12\,r+7\,r^2+4 \right)^2}{15552\,r^6}}
\end{multline*}
where
$A(\msP(A,N_1,Q_1,L_3,N_3,N_2))=
\Bigl[\sqrt{3} \bigl( 4\,y^2-8\,\sqrt{3} y\, x-24\, x-24\,r+8\,\sqrt{3} y+12\,r^2+
4\,\sqrt{3}r^2 y+6\,r^4\sqrt{3} y\, x+24\,r\, x-4\,r^3\sqrt{3} y+3\,r^4\,y^2-
8\,\sqrt{3}r\, y-12\, x^2r^2-12\,r^3 x+9\,r^4 x^2+12\, x^2+12\,r^2 x+
4\,r^2\,y^2+12 \bigr) \Bigr]\Big/\Bigl[24\,r^2\Bigr]
$.

\noindent \textbf{Case 13:}
\begin{multline*}
P(X_2 \in \NPE^r(X_1) \cup \G_1^r(X_1), X_1 \in T_s)=\\
\left(\int_{s_{10}}^{s_{14}}\int_{r_2(x)}^{r_{10}(x)} +
\int_{s_{14}}^{s_{13}}\int_{r_2(x)}^{r_{12}(x)} +
\int_{s_{13}}^{1/2}\int_{r_3(x)}^{r_{12}(x)} \right)
\frac{A(\msP(A,N_1,Q_1,L_3,N_3,N_2))}{A(\TY)^2}dydx=\\
\Bigl[1024-12288\,r+295680\,r^7+1053\,r^{12}-197140\,r^8+626688\,r^3-100864\,r^2-
1294848\,r^4-686528\,r^6+1282560\,r^5+\\
114336\,r^9-30930\,r^{10} \Bigr]\Big/\Bigl[31104\,r^4\Bigr]
\end{multline*}
where
$A(\msP(A,N_1,Q_1,L_3,N_3,N_2))=
\Bigl[\sqrt{3} \bigl( 4\,y^2-8\,\sqrt{3} y\, x-24\, x-24\,r+8\,\sqrt{3} y+12\,r^2+
4\,\sqrt{3}r^2 y+6\,r^4\sqrt{3} y\, x+24\,r\, x-4\,r^3\sqrt{3} y+3\,r^4\,y^2-
8\,\sqrt{3}r\, y-12\, x^2r^2-12\,r^3 x+9\,r^4 x^2+12\, x^2+12\,r^2 x+
4\,r^2\,y^2+12 \bigr)  \Bigr]\Big/\Bigl[24\,r^2\Bigr]
$.

\noindent \textbf{Case 14:}
\begin{multline*}
P(X_2 \in \NPE^r(X_1) \cup \G_1^r(X_1), X_1 \in T_s)=\\
\left(\int_{s_7}^{s_{11}}\int_{r_9(x)}^{\ell_{am}(x)} +
\int_{s_{11}}^{s_{10}}\int_{r_9(x)}^{r_{12}(x)} +
\int_{s_{10}}^{s_{14}}\int_{r_{10}(x)}^{r_{12}(x)} \right)
\frac{A(\msP(A,N_1,Q_1,L_3,L_4,Q_2,N_2))}{A(\TY)^2}dydx=\\
-\Bigl[ \left( r -1 \right)  \bigl( 1512\,r^{17}+1512\,r^{16}-16740\,r^{15}+540\,r^{14}+
84078\,r^{13}-83538\,r^{12}-164835\,r^{11}+401085\,r^{10}-487872\,r^9+\\
535728\,r^8-463124\,r^7+335596\,r^6-197440\,r^5+64640\,r^4-7936\,r^3-1792\,r^2+5632\,r-
512 \bigr)\Bigr]\Big/\Bigl[5184\, \left( 2\,r^2+1 \right)^3r^4\Bigr]
\end{multline*}
where
$A(\msP(A,N_1,Q_1,L_3,L_4,Q_2,N_2))=
\Bigl[\sqrt{3} \bigl( -6\, x-12\,r+6\,r^2+6\,r\, x+2\,\sqrt{3}r^2 y-r^2\,y^2-
2\,\sqrt{3} y\, x+r^4\,y^2+5\,y^2-2\,r^2 x\,\sqrt{3} y+2\,r^4\sqrt{3} y\, x+
2\,\sqrt{3}r\, y-2\,r^3\sqrt{3} y-3\,x^2r^2-6\,r^3 x+3\,r^4 x^2-2\,\sqrt{3} y+
3\, x^2+6\,r^2 x+6 \bigr) \Bigr]\Big/\Bigl[ 6\,r^2\Bigr]
$.

\noindent \textbf{Case 15:}
\begin{multline*}
P(X_2 \in \NPE^r(X_1) \cup \G_1^r(X_1), X_1 \in T_s)=
\int_{s_{13}}^{1/2}\int_{r_2(x)}^{r_3(x)} \frac{A(\msP(A,N_1,Q_1,G_3,M_2,N_3,N_2))}{A(\TY)^2}dydx=\\
{\frac{ \left( 147\,r^5-612\,r^4+980\,r^3-768\,r^2+744\,r -288 \right)  \left( -6+5\,r \right)^2}{7776\,r}}
\end{multline*}
where
$A(\msP(A,N_1,Q_1,L_3,L_4,Q_2,N_2))=
\Bigl[\sqrt{3} \bigl( 4\,r\,y^2+12\, x+9\,r -12+9\,r^3 x^2+12\,r\, x-12\,r\,x^2-
4\,\sqrt{3}r^2 y+4\,\sqrt{3}r\, y+6\,\sqrt{3}r^3 y\, x+3\,r^3\,y^2-12\,r^2 x-
4\,\sqrt{3} y \bigr) \Bigr]\Big/\Bigl[ 24\,r\Bigr]
$.

\noindent \textbf{Case 16:}
\begin{multline*}
P(X_2 \in \NPE^r(X_1) \cup \G_1^r(X_1), X_1 \in T_s)=
\int_{s_{14}}^{1/2}\int_{r_{12}(x)}^{r_{10}(x)} \frac{A(\msP(A,N_1,Q_1,L_3,N_3,N_2))}{A(\TY)^2}dydx=\\
-{\frac{ \left( 13\,r^8+52\,r^7+10\,r^6-184\,r^5+60\,r^4+624\,r^3-48\,r^2-832\,r+448 \right)
\left( -2+r \right)  \left( r^2+2\,r -4 \right)^2}{384\, \left( r+2 \right)^3r^2}}
\end{multline*}
where
$A(\msP(A,N_1,Q_1,L_3,N_3,N_2))=
\Bigl[\sqrt{3} \bigl( 4\,y^2-8\,\sqrt{3} y\, x-24\, x-24\,r+8\,\sqrt{3} y+12\,r^2+
4\,\sqrt{3}r^2 y+6\,r^4\sqrt{3} y\, x+24\,r\, x-4\,r^3\sqrt{3} y+3\,r^4\,y^2-
8\,\sqrt{3}r\, y-12\, x^2r^2-12\,r^3 x+9\,r^4 x^2+12\, x^2+12\,r^2 x+
4\,r^2\,y^2+12 \bigr)\Bigr]\Big/\Bigl[24\,r^2\Bigr]
$.

\noindent \textbf{Case 17:}
\begin{multline*}
P(X_2 \in \NPE^r(X_1) \cup \G_1^r(X_1), X_1 \in T_s)=
\left(\int_{s_{11}}^{s_{14}}\int_{r_{12}(x)}^{\ell_{am}(x)} +
\int_{s_{14}}^{1/2}\int_{r_{10}(x)}^{\ell_{am}(x)}\right)
\frac{A(\msP(A,N_1,Q_1,L_3,L_4,Q_2,N_2))}{A(\TY)^2}dydx=\\
\Bigl[\bigl( 189\,r^{12}+1323\,r^{11}+1026\,r^{10}-10692\,r^9-14364\,r^8+51732\,r^7+
64664\,r^6-183952\,r^5-153504\,r^4+398080\,r^3+124928\,r^2-\\
470528\,r+197632 \bigr)  \left( r -1 \right) \Bigr]\Big/\Bigl[5184\,r^2 \left( r+2 \right)^3\Bigr]
\end{multline*}
where
$A(\msP(A,N_1,Q_1,L_3,N_3,N_2))=
\Bigl[\sqrt{3} \bigl( -6\, x-12\,r+6\,r^2+6\,r\, x+2\,\sqrt{3}r^2 y-r^2\,y^2-
2\,\sqrt{3} y\, x+r^4\,y^2+5\,y^2-2\,r^2 x\,\sqrt{3} y+2\,r^4\sqrt{3} y
\, x+2\,\sqrt{3}r\, y-2\,r^3\sqrt{3} y-3\, x^2r^2-6\,r^3 x+3\,r^4 x^2-
2\,\sqrt{3} y+3\, x^2+6\,r^2 x+6 \bigr)\Bigr]\Big/\Bigl[6\,r^2\Bigr]
$.
}

Adding up the $P(X_2 \in \NPE^r(X_1) \cup \G_1^r(X_1), X_1 \in T_s)$
values in the 17 possible cases above, and multiplying by 6
we get for $r \in [1,4/3)$,
$$\nu_\lo(r)=
{\frac{860\,r^4-195\,r^5-256+720\,r -
846\,r^3-108\,r^2+47\,r^6}{108\,r^2 \left( r+2 \right)  \left( r+1 \right) }}.$$
The $\nu_\lo(r)$ values for the other intervals can be calculated similarly.

\subsection*{Derivation of $\nu_\lo(r)$ in Theorem \ref{thm:asy-norm-under}}

By symmetry,
$P(\{X_2,X_3\} \subset \NPE^r(X_1)\cup \G_1^r(X_1))=
6\,P(\{X_2,X_3\} \subset \NPE^r(X_1)\cup \G_1^r(X_1),\; X_1 \in T_s)$.
For $r \in \bigl[6/5,\sqrt{5}-1\bigr)$,
there are 17 cases to consider for calculation of $\nu_{\lo}(r)$ in the underlying graph case
(see also Figure \ref{fig:cases-AND-OR}):
{\small
\noindent \textbf{Case 1:}
\begin{multline*}
P(\{X_2,X_3\} \subset \NPE^r(X_1)\cup \G_1^r(X_1),\; X_1 \in T_s)=
\left(\int_{0}^{s_0}\int_{0}^{\ell_{am}(x)} +
\int_{s_0}^{s_1}\int_{r_1(x)}^{\ell_{am}(x)} \right)
\frac{A(\msP(A,M_1,M_C,M_3))^2}{A(\TY)^3}dydx=\\
{\frac{4}{81}}\,r^2-{\frac{4}{27}}\,r+1/9
\end{multline*}
where
$A(\msP(A,M_1,M_C,M_3))=
1/12\,\sqrt{3}
$.

\noindent \textbf{Case 2:}
\begin{multline*}
P(\{X_2,X_3\} \subset \NPE^r(X_1)\cup \G_1^r(X_1),\; X_1 \in T_s)=\\
\left(\int_{s_0}^{s_1}\int_{0}^{r_1(x)} +
\int_{s_1}^{s_3}\int_{0}^{r_2(x)} +
\int_{s_3}^{s_4}\int_{0}^{r_5(x)} +
\int_{s_4}^{s_5}\int_{r_3(x)}^{r_5(x)}\right)
\frac{A(\msP(A,M_1,L_2,L_3,M_C,M_3))^2}{A(\TY)^3}dydx=\\
-\Bigl[ \left( r -1 \right)  \bigl( 119155\,r^{11}-845345\,r^{10}+2724777\,r^9-5206743\,r^8+
6475257\,r^7-5454855\,r^6+3155193\,r^5-1249479\,r^4+\\
332181\,r^3-56619\,r^2+5589\,r -243 \bigr) \Bigr]\Big/\Bigl[25920\,r^{10}\Bigr]
\end{multline*}
where
$A(\msP(A,M_1,L_2,L_3,M_C,M_3))=
{\frac {\sqrt{3} \left( -4\,\sqrt{3}r\,y-12\,r+12\,r\,x+
5\,r^2+3\,y^2+6\,\sqrt{3}y-6\,\sqrt{3}y\,x+9-18\,x+9\,x^2 \right)}
{12\,r^2}}
$.

\noindent \textbf{Case 3:}
\begin{multline*}
P(\{X_2,X_3\} \subset \NPE^r(X_1)\cup \G_1^r(X_1),\; X_1 \in T_s)=
\left(\int_{s_4}^{s_5}\int_{0}^{r_3(x)} +
\int_{s_5}^{s_6}\int_{0}^{r_5(x)}\right)
\frac{A(\msP(A,G_2,G_3,M_2,M_C,M_3))^2}{A(\TY)^3}dydx=\\
{\frac{ \left( 215\,r^8-136\,r^7-56\,r^6+172\,r^5-55\,r^4-60\,r^3+66\,r^2-24\,r+
3 \right)  \left( r -1 \right)^4}{2880\,r^{10}}}
\end{multline*}
where
$A(\msP(A,G_2,G_3,M_2,M_C,M_3))=
-{\frac {\sqrt{3} \left( y^2+2\,\sqrt{3}y-2\,\sqrt{3}y\,x+3-
6\,x+3\,x^2-2\,r^2 \right) }{12\,r^2}}
$.

\noindent \textbf{Case 4:}
\begin{multline*}
P(\{X_2,X_3\} \subset \NPE^r(X_1)\cup \G_1^r(X_1),\; X_1 \in T_s)=\\
\left(\int_{s_1}^{s_2}\int_{r_2(x)}^{\ell_{am}(x)} +
\int_{s_2}^{s_3}\int_{r_2(x)}^{r_5(x)} \right)
\frac{A(\msP(A,M_1,L_2,L_3,L_4,L_5,M_3))^2}{A(\TY)^3}dydx=\\
\Bigl[ \bigl( 37072\,r^8-195072\,r^7+453120\,r^6-589248\,r^5+460728\,r^4-
217728\,r^3+60480\,r^2-9072\,r+567 \bigr) \\
\left( 4\,r -3+\sqrt{3} \right)^2 \left( 4\,r -3-\sqrt{3} \right)^2\Bigr]\Big/\Bigl[1866240\,r^{10}\Bigr]
\end{multline*}
where
$A(\msP(A,M_1,L_2,L_3,L_4,L_5,M_3))=
{\frac {\sqrt{3} \left( 4\,\sqrt{3}r\,y+9\,r^2-
24\,\nu+12\,r\,x+15\,y^2-6\,\sqrt{3}y-6\,
\sqrt{3}y\,x+18-18\,x+9\,x^2 \right)}{12\,r^2}}
$.

\noindent \textbf{Case 5:}
\begin{multline*}
P(\{X_2,X_3\} \subset \NPE^r(X_1)\cup \G_1^r(X_1),\; X_1 \in T_s)=\\
\left(\int_{s_5}^{s_6}\int_{r_5(x)}^{r_7(x)} +
\int_{s_6}^{s_9}\int_{0}^{r_7(x)} \right)
\frac{A(\msP(A,G_2,G_3,M_2,M_C,P_2,N_2))^2}{A(\TY)^3}dydx=\\
{\frac{ \left( 3-12\,r -15\,r^2+84\,r^3+18\,r^4-232\,r^5+130\,r^6+504\,r^7-108\,r^8-
288\,r^9+623\,r^{10}+920\,r^{11}+373\,r^{12} \right)
\left( r -1 \right)^3}{2880\,r^{10} \left( r+1 \right)^5}}
\end{multline*}
where
$A(\msP(A,G_2,G_3,M_2,M_C,P_2,N_2))=
\Bigl[\sqrt{3} \bigl( -2\,y^2-4\,\sqrt{3}y+4\,\sqrt{3}y\,x-6+12\,x-6\,x^2+
7\,r^2-4\,r^3\sqrt{3}y-12\,r^3x+8\,r^4\sqrt{3}y\,x+12\,r^4x^2+
4\,r^4 y^2 \bigr) \Bigr]\Big/\Bigl[24\,r^2\Bigr]
$.

\noindent \textbf{Case 6:}
\begin{multline*}
P(\{X_2,X_3\} \subset \NPE^r(X_1)\cup \G_1^r(X_1),\; X_1 \in T_s)=\\
\left(\int_{s_5}^{s_9}\int_{r_7(x)}^{r_3(x)} +
\int_{s_9}^{s_{12}}\int_{0}^{r_3(x)} +
\int_{s_{12}}^{1/2}\int_{0}^{r_6(x)} \right)
\frac{A(\msP(A,N_1,Q_1,G_3,M_2,M_C,P_2,N_2))^2}{A(\TY)^3}dydx=\\
-\Bigl[19683\,r^{15}-59049\,r^{14}+83106\,r^{13}+167670\,r^{12}-211626\,r^{11}+
344466\,r^{10}-142614\,r^9-2573586\,r^8-128853\,r^7+\\
3465675\,r^6+1103824\,r^5-1473304\,r^4-730880\,r^3+107776\,r^2+
158720\,r+31744\Bigr]\Big/\Bigl[1049760\, \left( r+1 \right)^5r^6\Bigr]
\end{multline*}
where
$A(\msP(A,N_1,Q_1,G_3,M_2,M_C,P_2,N_2))=
\Bigl[\sqrt{3} \bigl( 4\,r\,y^2+12\,x+13\,r+12\,r\,x-4\,\sqrt{3}y-12+
4\,\sqrt{3}r\,y-8\,\sqrt{3}r^2y+18\,x^2r^3-12\,r\,x^2+6\,r^3y^2-
24\,r^2 x+12\,\sqrt{3}r^3y\,x \bigr)\Bigr]\Big/\Bigl[24\,r\Bigr]
$.

\noindent \textbf{Case 7:}
\begin{multline*}
P(\{X_2,X_3\} \subset \NPE^r(X_1)\cup \G_1^r(X_1),\; X_1 \in T_s)=\\
\left(\int_{s_8}^{s_5}\int_{r_8(x)}^{r_2(x)} +
\int_{s_5}^{s_{10}}\int_{r_3(x)}^{r_2(x)} +
\int_{s_{10}}^{s_{12}}\int_{r_3(x)}^{r_6(x)} \right)
\frac{A(\msP(A,N_1,Q_1,L_3,M_C,P_2,N_2))^2}{A(\TY)^3}dydx=\\
-\Bigl[6144-110592\,r -310846464\,r^7+2127553557\,r^{12}+570050560\,r^8-
5031936\,r^3+936960\,r^2+19526656\,r^4+147203072\,r^6+\\
7627473\,r^{20}+1419072042\,r^{16}-762467328\,r^{17}+288811029\,r^{18}
-68327424\,r^{19}-59166720\,r^5-923627520\,r^9+\\
1340817105\,r^{10}-1765251072\,r^{11}-2350015488\,r^{13}+
2339575338\,r^{14}-2016377856\,r^{15}\Bigr]\Big/\Bigl[ 262440\,\left( r^2+1 \right)^5r^{10}\Bigr]
\end{multline*}
where
$A(\msP(A,N_1,Q_1,L_3,M_C,P_2,N_2))=
\Bigl[\sqrt{3} \bigl( -4\,\sqrt{3}r\,y+2\,\sqrt{3}r^2 y-6\,x^2r^2-12\,x-12\,r-
12\,r^3 x+9\,r^4x^2+8\,r^2+12\,r\,x+6\,x^2+6\,r^4\sqrt{3}y\,x+2\,r^2 y^2-
4\,\sqrt{3}y\,x+3\,r^4 y^2-4\,r^3\sqrt{3}y+4\,\sqrt{3}y+2\,y^2+
6\,r^2x+6 \bigr)\Bigr]\Big/\Bigl[12\,r^2\Bigr]
$.

\noindent \textbf{Case 8:}
\begin{multline*}
P(\{X_2,X_3\} \subset \NPE^r(X_1)\cup \G_1^r(X_1),\; X_1 \in T_s)=\\
\left(\int_{s_3}^{s_8}\int_{r_5(x)}^{r_2(x)} +
\int_{s_8}^{s_5}\int_{r_5(x)}^{r_8(x)} \right)
\frac{A(\msP(A,N_1,P_1,L_2,L_3,M_C,P_2,N_2))^2}{A(\TY)^3}dydx=\\
\Bigl[\bigl( 426497\,r^{16}-2443992\,r^{15}+6726107\,r^{14}-11753232\,r^{13}+
15220771\,r^{12}-16367448\,r^{11}+15754449\,r^{10}-13773024\,r^9+\\
10839672\,r^8-7552440\,r^7+4592889\,r^6-2374272\,r^5+
1018899\,r^4-344088\,r^3+81891\,r^2-11664\,r+729 \bigr) \\
\left( -12\,r+7\,r^2+3 \right)^2\Bigr]\Big/\Bigl[699840\, \left( r^2+1 \right)^5r^{10}\Bigr]
\end{multline*}
where
$A(\msP(A,N_1,P_1,L_2,L_3,M_C,P_2,N_2))=
\Bigl[\sqrt{3} \bigl( -4\,r^3\sqrt{3}y-12\,r^3 x+8\,r^4\sqrt{3}y\,x+12\,r^4 x^2+
4\,r^4y^2-4\,\sqrt{3}r\,y-12\,r+12\,r\,x+3\,y^2+6\,\sqrt{3}y-6\,\sqrt{3}y\,x+
8\,r^2+9-18\,x+9\,x^2 \bigr)\Bigr]\Big/\Bigl[12\,r^2\Bigr]
$.

\noindent \textbf{Case 9:}
\begin{multline*}
P(\{X_2,X_3\} \subset \NPE^r(X_1)\cup \G_1^r(X_1),\; X_1 \in T_s)=\\
\left(\int_{s_2}^{s_3}\int_{r_5(x)}^{\ell_{am}(x)} +
\int_{s_3}^{s_7}\int_{r_2(x)}^{\ell_{am}(x)} +
\int_{s_7}^{s_8}\int_{r_2(x)}^{r_8(x)} \right)
\frac{A(\msP(A,N_1,P_1,L_2,L_3,L_4,L_5,P_2,N_2))^2}{A(\TY)^3}dydx=\\
-\Bigl[15309-367416\,r+60475010560\,r^{28}+437704472832\,r^{26}+1444872192\,r^{30}-
13250101248\,r^{29}-185909870592\,r^{27}+\\
4148739\,r^2-2027754648576\,r^{23}+1397612375040\,r^{24}+20429177589\,r^8-
677278256112\,r^{13}-49656902904\,r^9+\\
159963012\,r^4-30005640\,r^3-681714144\,r^5-7515142416\,r^7-3097406755584\,r^{21}-2609245249920\,r^{17}+\\
3051035360256\,r^{18}-3315184235136\,r^{19}+3337272236928\,r^{20}+2631941507968\,r^{22}+
2435971806\,r^6+\\
109069315047\,r^{10}-218273842152\,r^{11}+400534503738\,r^{12}+1059615993384\,r^{14}-
1538314485120\,r^{15}+\\
2076627064432\,r^{16}-845838600192\,r^{25}\Bigr]\Big/\Bigl[1399680\, \left(
r^2+1 \right)^5 \left( 2\,r^2+1 \right)^5r^{10}\Bigr]
\end{multline*}
where
$A(\msP(A,N_1,P_1,L_2,L_3,L_4,L_5,P_2,N_2))=
\Bigl[\sqrt{3} \bigl( 18-18\,x-24\,r-12\,r^3 x+12\,r^4x^2+12\,r^2+12\,r\,x+
4\,\sqrt{3}r\,y-4\,r^3\sqrt{3}y+4\,r^4 y^2-6\,\sqrt{3}y\,x+8\,r^4\sqrt{3} y\,x+
9\,x^2+15\,y^2-6\,\sqrt{3}y \bigr) \Bigr]\Big/\Bigl[12\,r^2\Bigr]
$.

\noindent \textbf{Case 10:}
\begin{multline*}
P(\{X_2,X_3\} \subset \NPE^r(X_1)\cup \G_1^r(X_1),\; X_1 \in T_s)=\\
\left(\int_{s_7}^{s_8}\int_{r_8(x)}^{r_9(x)} +
\int_{s_8}^{s_{10}}\int_{r_2(x)}^{r_9(x)} \right)
\frac{A(\msP(A,N_1,Q_1,L_3,L_4,L_5,P_2,N_2))^2}{A(\TY)^3}dydx=\\
\Bigl[64\,\bigl( 12-144\,r+924\,r^2-683328\,r^{23}+112976\,r^{24}+757211\,r^8-
10554918\,r^{13}-1513230\,r^9+16242\,r^4-4320\,r^3-51372\,r^5-\\
344988\,r^7-4867848\,r^{21}-18583080\,r^{17}+16493828\,r^{18}-
12883116\,r^{19}+8668124\,r^{20}+2177536\,r^{22}+141366\,r^6+2774371\,r^{10}-\\
4692510\,r^{11}+7331714\,r^{12}+14002613\,r^{14}-16948218\,r^{15}+18708475\,r^{
16} \bigr)  \left( r -1 \right)^2 \left( 2\,r -1 \right)^2\Bigr]\Big/\Bigl[32805\,\left( r^2+1 \right)^5 \left( 2\,r^2+1 \right)^5r^8\Bigr]
\end{multline*}
where
$A(\msP(A,N_1,Q_1,L_3,L_4,L_5,P_2,N_2))=
\Bigl[\sqrt{3} \bigl( 2\,\sqrt{3}r^2y+15-6\,x^2r^2-12\,x-24\,r-12\,r^3x+
9\,r^4 x^2+12\,r^2+12\,r\,x-8\,\sqrt{3} y+6\,x^2+6\,r^4\sqrt{3}y\,x+
14\,y^2-4\,\sqrt{3}y\,x+2\,r^2y^2-4\,r^3\sqrt{3}y+3\,r^4y^2+6\,r^2 x+
4\,\sqrt{3}r\,y \bigr)\Bigr]\Big/\Bigl[12\,r^2\Bigr]
$.

\noindent \textbf{Case 11:}
\begin{multline*}
P(\{X_2,X_3\} \subset \NPE^r(X_1)\cup \G_1^r(X_1),\; X_1 \in T_s)=\\
\left(\int_{s_{12}}^{s_{13}}\int_{r_6(x)}^{r_3(x)} +
\int_{s_{13}}^{1/2}\int_{r_6(x)}^{r_2(x)} \right)
\frac{A(\msP(A,N_1,Q_1,G_3,M_2,N_3,N_2))^2}{A(\TY)^3}dydx=\\
-\Bigl[-253952+1529856\,r^2+601574256\,r^8-385780320\,r^{13}-776518272\,r^9+
7803648\,r^4-70917120\,r^5-396524160\,r^7+\\
209710080\,r^6+869661288\,r^{10}-845940960\,r^{11}+668092108\,r^{12}+
147067614\,r^{14}-32610600\,r^{15}+3173067\,r^{16}\Bigr]\Big/\Bigl[8398080\,r^6\Bigr]
\end{multline*}
where
$A(\msP(A,N_1,Q_1,G_3,M_2,N_3,N_2))=
\Bigl[\sqrt{3} \bigl( 4\,r\,y^2+12\,x+4\,\sqrt{3}r\,y+9\,r-4\,\sqrt{3}y+12\,r\,x-
12+9\,x^2r^3+6\,\sqrt{3}r^3y\,x-12\,r\,x^2-4\,\sqrt{3}r^2y-12\,r^2 x+3\,r^3y^2 \bigr)\Bigr]\Big/\Bigl[24\,r\Bigr]
$.

\noindent \textbf{Case 12:}
\begin{multline*}
P(\{X_2,X_3\} \subset \NPE^r(X_1)\cup \G_1^r(X_1),\; X_1 \in T_s)=\\
\left(\int_{s_{10}}^{s_{12}}\int_{r_6(x)}^{r_2(x)} +
\int_{s_{12}}^{s_{13}}\int_{r_3(x)}^{r_2(x)} \right)
\frac{A(\msP(A,N_1,Q_1,L_3,N_3,N_2))^2}{A(\TY)^3}dydx=\\
\Bigl[\bigl( 64827\,r^{16}-444528\,r^{15}+1223334\,r^{14}-1793232\,r^{13}+
1839416\,r^{12}-2003712\,r^{11}+2286224\,r^{10}-2421504\,r^9+3095088\,r^8-\\
4428288\,r^7+5889152\,r^6-6093312\,r^5+4557056\,r^4-2340864\,r^3+774144\,r^2-147456\,r+12288 \bigr)\\
 \left( -12\,r+7\,r^2+4 \right)^2\Bigr]\Big/\Bigl[8398080\,r^{10}\Bigr]
\end{multline*}
where
$A(\msP(A,N_1,Q_1,L_3,N_3,N_2))=
\Bigl[\sqrt{3} \bigl( -12\,x^2r^2-24\,x-24\,r-12\,r^3x+9\,r^4x^2+4\,y^2-
8\,\sqrt{3}r\,y+6\,r^4\sqrt{3}y\,x+8\,\sqrt{3}y+12\,r^2+24\,r\,x+
12\,x^2-8\,\sqrt{3}y\,x+4\,r^2y^2-4\,r^3\sqrt{3}y+3\,r^4 y^2+
4\,\sqrt{3}r^2 y+12\,r^2 x+12 \bigr)\Bigr]\Big/\Bigl[24\,r^2\Bigr]
$.

\noindent \textbf{Case 13:}
\begin{multline*}
P(\{X_2,X_3\} \subset \NPE^r(X_1)\cup \G_1^r(X_1),\; X_1 \in T_s)=\\
\left(\int_{s_{10}}^{s_{14}}\int_{r_2(x)}^{r_{10}(x)} +
\int_{s_{14}}^{s_{13}}\int_{r_2(x)}^{r_{12}(x)} +
\int_{s_{13}}^{1/2}\int_{r_3(x)}^{r_{12}(x)} \right)
\frac{A(\msP(A,N_1,Q_1,L_3,N_3,N_2))^2}{A(\TY)^3}dydx=\\
\Bigl[196608-3538944\,r+8927944704\,r^7-1883996112\,r^{12}-9492593152\,r^8-
146866176\,r^3+29196288\,r^2+220250112\,r^4-\\
4486594560\,r^6+213597\,r^{20}-259250904\,r^{16}+69124752\,r^{17}-10683306\,r^{18}+864387072\,r^5+5220357120\,r^9-\\
1081136256\,r^{10}+602097408\,r^{11}+2223664128\,r^{13}-1509638512\,r^{14}+
716568768\,r^{15}\Bigr]\Big/\Bigl[16796160\,r^8\Bigr]
\end{multline*}
where
$A(\msP(A,N_1,Q_1,L_3,N_3,N_2))=
\Bigl[\sqrt{3} \bigl( -12\,x^2r^2-24\,x-24\,r-12\,r^3x+9\,r^4x^2+4\,y^2-
8\,\sqrt{3}r\,y+6\,r^4\sqrt{3}y\,x+8\,\sqrt{3}y+12\,r^2+24\,r\,x+
12\,x^2-8\,\sqrt{3}y\,x+4\,r^2 y^2-4\,r^3\sqrt{3}y+3\,r^4 y^2+
4\,\sqrt{3}r^2y+12\,r^2x+12 \bigr)\Bigr]\Big/\Bigl[24\,r^2\Bigr]
$.

\noindent \textbf{Case 14:}
\begin{multline*}
P(\{X_2,X_3\} \subset \NPE^r(X_1)\cup \G_1^r(X_1),\; X_1 \in T_s)=\\
\left(\int_{s_7}^{s_{11}}\int_{r_9(x)}^{\ell_{am}(x)} +
\int_{s_{11}}^{s_{10}}\int_{r_9(x)}^{r_{12}(x)} +
\int_{s_{10}}^{s_{14}}\int_{r_{10}(x)}^{r_{12}(x)} \right)
\frac{A(\msP(A,N_1,Q_1,L_3,L_4,Q_2,N_2))^2}{A(\TY)^3}dydx=\\
-\Bigl[\left( r -1 \right)  \bigl( -16384+278528\,r+215136\,r^{28}+40176\,r^{26}+
215136\,r^{29}-3381264\,r^{27}-2301952\,r^2-99212040\,r^{23}-\\
25050384\,r^{24}-312101312\,r^8-7215869272\,r^{13}-147586784\,r^9-
42770432\,r^4+12591104\,r^3+114049024\,r^5+345810944\,r^7+\\
55914462\,r^{21}-2082969096\,r^{17}+43443459\,r^{18}+826941555\,r^{19}-641846754\,r^{20}+
209930616\,r^{22}-232963072\,r^6+\\
1311322268\,r^{10}-3191747236\,r^{11}+
5434516904\,r^{12}+7756861008\,r^{14}-6865898928\,r^{15}+4727296416\,r^{16}+26115696\,r^{25}
\bigr)\Bigr]\Big/\\
\Bigl[466560\, \left( 2\,r^2+1 \right)^5r^8\Bigr]
\end{multline*}
where
$A(\msP(A,N_1,Q_1,L_3,L_4,Q_2,N_2))=
\Bigl[\sqrt{3} \bigl( -3\,x^2r^2-6\,x-12\,r-6\,r^3 x+3\,r^4 x^2+2\,\sqrt{3}r\,y+
6\,r^2+6\,r\,x+3\,x^2-2\,\sqrt{3}y-2\,\sqrt{3}r^2y\,x+2\,r^4\sqrt{3}y\,x+
2\,\sqrt{3}r^2y-r^2 y^2+5\,y^2-2\,r^3\sqrt{3}y+r^4 y^2-2\,\sqrt{3}y\,x+
6+6\,r^2x \bigr)\Bigr]\Big/\Bigl[6\,r^2\Bigr]
$.

\noindent \textbf{Case 15:}
\begin{multline*}
P(\{X_2,X_3\} \subset \NPE^r(X_1)\cup \G_1^r(X_1),\; X_1 \in T_s)=
\int_{s_{13}}^{1/2}\int_{r_2(x)}^{r_3(x)} \frac{A(\msP(A,N_1,Q_1,G_3,M_2,N_3,N_2))^2}{A(\TY)^3}dydx=\\
\Bigl[\bigl( 63855\,r^{10}-498960\,r^9+1650060\,r^8-3036960\,r^7+3703292\,r^6-
3657696\,r^5+3268368\,r^4-2419200\,r^3+1550448\,r^2-\\
725760\,r+155520 \bigr)  \left( -6+5\,r \right)^2\Bigr]\Big/\Bigl[4199040\,r^2\Bigr]
\end{multline*}
where
$A(\msP(A,N_1,Q_1,G_3,M_2,N_3,N_2))=
\Bigl[\sqrt{3} \bigl( 4\,r\,y^2+12\,x+4\,\sqrt{3}r\,y+9\,r-4\,\sqrt{3}y+12\,r\,x-
12+9\,x^2r^3+6\,\sqrt{3}r^3y\,x-12\,r\,x^2-4\,\sqrt{3}r^2y-12\,r^2 x+
3\,r^3y^2 \bigr)\Bigr]\Big/\Bigl[24\,r\Bigr]
$.

\noindent \textbf{Case 16:}
\begin{multline*}
P(\{X_2,X_3\} \subset \NPE^r(X_1)\cup \G_1^r(X_1),\; X_1 \in T_s)=
\int_{s_{14}}^{1/2}\int_{r_{12}(x)}^{r_{10}(x)} \frac{A(\msP(A,N_1,Q_1,L_3,N_3,N_2))^2}{A(\TY)^3}dydx=\\
-\Bigl[\bigl( 293\,r^{16}+2344\,r^{15}+4662\,r^{14}-9088\,r^{13}-32320\,r^{12}+42976\,r^{11}+
175408\,r^{10}-119680\,r^9-544144\,r^8+372352\,r^7+\\
1216512\,r^6-882688\,r^5-1564672\,r^4+1373184\,r^3+924672\,r^2-1314816\,r+380928 \bigr)\\
\left( -2+r \right)  \left( r^2+2\,r -4 \right)^2\Bigr]\Big/\Bigl[23040\, \left( r+2 \right)^5r^4\Bigr]
\end{multline*}
where
$A(\msP(A,N_1,Q_1,L_3,N_3,N_2))=
\Bigl[\sqrt{3} \bigl( -12\,x^2r^2-24\,x-24\,r-12\,r^3x+9\,r^4x^2+
4\,y^2-8\,\sqrt{3}r\,y+6\,r^4\sqrt{3}y\,x+8\,\sqrt{3}y+12\,r^2+24\,r\,x+
12\,x^2-8\,\sqrt{3}y\,x+4\,r^2y^2-4\,r^3\sqrt{3}y+3\,r^4y^2+4\,\sqrt{3}
r^2 y+12\,r^2x+12 \bigr)\Bigr]\Big/\Bigl[24\,r^2\Bigr]
$.

\noindent \textbf{Case 17:}
\begin{multline*}
P(\{X_2,X_3\} \subset \NPE^r(X_1)\cup \G_1^r(X_1),\; X_1 \in T_s)=\\
\left(\int_{s_{11}}^{s_{14}}\int_{r_{12}(x)}^{\ell_{am}(x)}+
\int_{s_{14}}^{1/2}\int_{r_{10}(x)}^{\ell_{am}(x)} \right)
\frac{A(\msP(A,N_1,Q_1,L_3,L_4,Q_2,N_2))^2}{A(\TY)^3}dydx=\\
\Bigl[\bigl( 6723\,r^{20}+73953\,r^{19}+213678\,r^{18}-433512\,r^{17}-2873232\,r^{16}+
627264\,r^{15}+20218896\,r^{14}+5675184\,r^{13}-97577924\,r^{12}-\\
39916108\,r^{11}+343932568\,r^{10}+108508576\,r^9-906967296\,r^8-96480192\,r^7+1702951296\,r^6-
293251072\,r^5-1994987520\,r^4+\\
981590016\,r^3+1118830592\,r^2-1135919104\,r+287604736 \bigr)  \left( r -1
 \right) \Bigr]\Big/\Bigl[466560\,r^4 \left( r+2 \right)^5\Bigr]
\end{multline*}
where
$A(\msP(A,N_1,Q_1,L_3,N_3,N_2))=
\Bigl[\sqrt{3} \bigl( -3\,x^2r^2-6\,x-12\,r-6\,r^3x+3\,r^4x^2+2\,\sqrt{3}r\,y+
6\,r^2+6\,r\,x+3\,x^2-2\,\sqrt{3}y-2\,\sqrt{3}r^2 y\,x+2\,r^4\sqrt{3}y\,x+
2\,\sqrt{3}r^2y-r^2 y^2+5\,y^2-2\,r^3\sqrt{3}y+r^4 y^2-2\,\sqrt{3}y\,x+
6+6\,r^2x \bigr)\Bigr]\Big/\Bigl[6\,r^2\Bigr]
$.
}

Adding up the $P(\{X_2,X_3\} \subset \NPE^r(X_1)\cup \G_1^r(X_1),\; X_1 \in T_s)$
values in the 17 possible cases above, and multiplying by 6
we get, for $r \in \bigl[6/5,\sqrt{5}-1 \bigr)$,
\begin{multline*}
\nu_\lo(r)=
-\Bigl[-413208\,r+3070468\,r^2-74801558\,r^8+75243552\,r^{13}-4883958\,r^9+
14541630\,r^4+28880-11254002\,r^3-\\
3667716\,r^5+64360782\,r^7+13122\,r^{21}-3300900\,r^{17}+156014\,r^{18}-
175011\,r^{19}+62825\,r^{20}+1458\,r^{22}-19812000\,r^6+\\
99831906\,r^{10}-120628524\,r^{11}+33155180\,r^{12}-67685050\,r^{14}+
5055135\,r^{15}+11053023\,r^{16}\Bigr]\Big/\Bigl[116640\,r^6 \left( r^2+1 \right)\\
 \left( 2\,r^2+1 \right)  \left( r+2 \right)^3 \left( r+1 \right)^3\Bigr].
\end{multline*}
The $\nu_\lo(r)$ values for the other intervals can be calculated similarly.
}

\subsection*{Appendix 4: Proof of Corollary \ref{cor:MT-asy-norm-NYr}:}
Recall that
$S^{\la}_n(r)=\rho^{\la}_{I,n}(r)$ is
the relative edge density of the reflexivity graph for the multiple triangle case.
Then the expectation of $S^{\la}_n(r)$ is
$$ \E\left[S^{\la}_n(r)\right]=\frac{2}{n\,(n-1)}
\sum\hspace*{-0.1 in}\sum_{i < j \hspace*{0.25 in}}
\hspace*{-0.1 in}\,\E\left[h^{\la}_{ij}(r)\right]=\E\left[ h^{\la}_{12}(r) \right]=
P(X_2\in \NPE^r(X_1) \cap \G_1^r(X_1))=\widetilde p_{\la}(r).$$
But, by definition of $\NPE^r(\cdot)$ and $\G_1^r(\cdot)$, if $X_1$
and $X_2$ are in different triangles, then $P(X_2 \in \NPE^r(X_1)\cap
\G_1^r(X_1))=0$. So by the law of total probability
\begin{eqnarray*}
\widetilde p_{\la}(r)&:=&P(X_2 \in \NPE^r(X_1)\cap \G_1^r(X_1))\\
&=& \sum_{i=1}^{J_m}P(X_2 \in \NPE^r(X_1)\cap \G_1^r(X_1)\,|\,\{X_1,X_2\} \subset T_i)\,P(\{X_1,X_2\} \subset T_i)\\
&=& \sum_{i=1}^{J_m}p_{\la}(r)\,P(\{X_1,X_2\} \subset T_i)
\text{ (since $P(X_2 \in \NPE^r(X_1)\cap \G_1^r(X_1)\,|\,\{X_1,X_2\} \subset T_i)=p_{\la}(r)$)}\\
&=& p_{\la}(r) \, \sum_{i=1}^{J_m}\left(\frac{A(T_i)}{\sum_{i=1}^{J_m}A(T_i)}\right)^2
\text{ (since $P(\{X_1,X_2\} \subset T_i)=\left(\frac{A(T_i)}{\sum_{i=1}^{J_m}A(T_i)}\right)^2$)}
= p_{\la}(r) \, \left(\sum_{i=1}^{J_m}w_i^2\right).
\end{eqnarray*}
where $p_{\la}(r)$ is given by Equation \eqref{eqn:Asymean_and}.

Likewise, we get $\widetilde p_{\lo}(r)=p_{\lo}(r)\,\left(\sum_{i=1}^{J_m}w_i^2\right)$
where $p_{\lo}(r)$ is given by Equation \eqref{eqn:Asymean_or}.

Furthermore, the asymptotic variance is
$$\widetilde \nu_{\la}(r)=\E\left[ h^{\la}_{12}(r)h^{\la}_{13}(r) \right]-\E\left[ h^{\la}_{12}(r) \right]\E\left[ h^{\la}_{13}(r) \right]=
P\left(\{X_2,X_3\} \subset \NPE^r(X_1)\cap \G_1^r(X_1)\right)-\left(\widetilde p_{\la}(r)\right)^2.$$
Then for $J_m>1$, we have
\begin{multline*}
P(\{X_2,X_3\} \subset \NPE^r(X_1)\cap \G_1^r(X_1))=
\sum_{i=1}^{J_m}P(\{X_2,X_3\} \subset \NPE^r(X_1)\cap \G_1^r(X_1)\,|\, \{X_1,X_2,X_3\} \subset T_i)\, P(\{X_1,X_2,X_3\} \subset T_i)\\
 = P(\{X_2,X_3\} \subset \NPE^r(X_1)\cap \G_1^r(X_1)\,|\,\{X_1,X_2,X_3\} \subset T_e)\, \left(\sum_{i=1}^{J_m}w_i^3 \right).
\end{multline*}

Hence,
\begin{eqnarray*}
\widetilde \nu_{\la}(r)&=&P(\{X_2,X_3\} \subset \NPE^r(X_1)\cap \G_1^r(X_1)\,|\, \{X_1,X_2,X_3\} \subset T_e)\,
\left(\sum_{i=1}^{J_m}w_i^3 \right)-\left(\widetilde p_{\la}(r)\right)^2\\
&=&\nu_{\la}(r)\,\left(\sum_{i=1}^{J_m}w_i^3\right)+
p_{\la}(r)^2\,\left(\sum_{i=1}^{J_m}w_i^3-\left(\sum_{i=1}^{J_m}w_i^2\right)^2\right).
\end{eqnarray*}

Likewise, we get
$\widetilde \nu_{\lo}(r)=\nu_{\lo}(r)\,\left(\sum_{i=1}^{J_m}w_i^3\right)+
p_{\lo}(r)^2\,\left(\sum_{i=1}^{J_m}w_i^3-\left(\sum_{i=1}^{J_m}w_i^2\right)^2\right).$

So conditional on $\Y_m$, if $\widetilde \nu_{\la}(r)>0$ then
$\sqrt{n}\,\left(S^{\la}_n(r)-\widetilde p_{\la}(r)\right)
\stackrel {\mathcal L}{\longrightarrow} \mathcal N\left(0,\widetilde \nu_{\la}(r)\right)$.
A similar result holds for the underlying graph version.
$\blacksquare$

\subsection*{Appendix 5: Proof of Theorem \ref{thm:MT-asy-norm-II}:}
Recall that
$\rho^{\la}_{II,n}(r)$ is the version II of the
relative edge density of the reflexivity graph for the multiple triangle case.
Then the expectation of $\rho^{\la}_{II,n}(r)$ is
$$ \E\left[\rho^{\la}_{II,n}(r)\right]=
\sum_{i=1}^{J_m}\frac{n_i\,(n_i-1)}{2\,n_t}\,\E\left[\rho^{\la}_{{}_{[i]}}(r)\right]=p_{\la}(r)$$
since by \eqref{eqn:E[rho-and-D]} we have
$$\E[\rho^{\la}_{{}_{[i]}}(r)]=\frac{2}{n_i(n_i-1)}\sum\hspace*{-0.1 in}\sum_{k < l \hspace*{0.25 in}}
\hspace*{-0.1 in}\,\E\left[h^{\la}_{kl}(r)\right]=\E\left[ h^{\la}_{12}(r) \right]=p_{\la}(r)$$
where
$p_{\la}(r)$ is given by Equation \eqref{eqn:Asymean_and}.
Likewise, we get $\widetilde p_{\lo}(r)=p_{\lo}(r)$
where $p_{\lo}(r)$ is given by Equation \eqref{eqn:Asymean_or}.

Next,
$$ \Var\left[\rho^{\la}_{II,n}(r)\right]=
\sum_{i=1}^{J_m}\frac{n_i^2\,(n_i-1)^2}{4\,n_t^2}\,
\Var\left[\rho^{\la}_{{}_{[i]}}(r)\right]$$
since $\rho^{\la}_{{}_{[k]}}(r)$ and $\rho^{\la}_{{}_{[l]}}(r)$ are independent for $k\not=l$.
Then by \eqref{eqn:Var[rho-and-D]} we have
$$\Var\left[\rho^{\la}_{{}_{[i]}}(r)\right]=
\frac{2}{n_i\,(n_i-1)}\Var\left[ h^{\la}_{12}(r) \right]+
\frac{4\,(n_i-2)}{n_i\,(n_i-1)} \, \Cov\left[ h^{\la}_{12}(r),h^{\la}_{13}(r) \right].$$
So,
$$ \Var\left[\rho^{\la}_{II,n}(r)\right]=
\sum_{i=1}^{J_m}\frac{n_i\,(n_i-1)}{2\,n_t^2}\, \Var\left[ h^{\la}_{12}(r) \right]+
\sum_{i=1}^{J_m}\frac{n_i\,(n_i-1)\,(n_i-2)}{n_t^2}\, \Cov\left[ h^{\la}_{12}(r),h^{\la}_{13}(r) \right].$$
Here
$\sum_{i=1}^{J_m}\frac{n_i\,(n_i-1)}{2\,n_t^2}\, \Var\left[ h^{\la}_{12}(r) \right]=
\frac{1}{n_t}\, \Var\left[ h^{\la}_{12}(r) \right]$.
Then
for large $n_i$ and $n$,
$$\frac{1}{n_t}\, \Var\left[ h^{\la}_{12}(r) \right] \approx
\frac{2}{n^2 \sum_{i=1}^{J_m} w_i^2}\, \Var\left[ h^{\la}_{12}(r) \right]$$
since $\frac{n_t^2}{n^2}=\sum_{i=1}^{J_m}\frac{n_i\,(n_i-1)}{2\,n^2}$
and $n_i/n \rightarrow w_i$ as $n_i,n \rightarrow \infty$.
Similarly,
for large $n_i$ and $n$,
$$\sum_{i=1}^{J_m}\frac{n_i\,(n_i-1)\,(n_i-2)}{n_t^2}\, \Cov\left[ h^{\la}_{12}(r),h^{\la}_{13}(r) \right]
\approx
\left[\frac{4}{n} \left(\sum_{j=1}^{J_m} w_i^3\right)\Bigg/\left(\sum_{i=1}^{J_m} w_i^2\right)^2 \right] \Cov\left[ h^{\la}_{12}(r),h^{\la}_{13}(r) \right].$$

Hence, conditional on $\Y_m$,
$\sqrt{n}\,\left(\rho^{\la}_{II,n}(r)-\widetilde p_{\la}(r)\right)
\stackrel {\mathcal L}{\longrightarrow} \mathcal N\left(0,4\,\breve \nu_{\la}(r)\right)$
provided that $\breve \nu_{\la}(r)>0$
where $\breve p_{\la}(r)=p_{\la}(r)$ and
$\breve \nu_{\la}(r)=\nu_{\la}(r) \left(\sum_{i=1}^{J_m} w_i^3\right)\Big/\left(\sum_{i=1}^{J_m} w_i^2\right)^2 $.
A similar result holds for the underlying graph version.
$\blacksquare$

\clearpage

\newpage

\clearpage

\end{document}

%% file: Nofnu2.pstex_t
\begin{picture}(0,0)%
\includegraphics{Nofnu2.pstex}%
\end{picture}%
\setlength{\unitlength}{3947sp}%
\begingroup\makeatletter\ifx\SetFigFont\undefined%
\gdef\SetFigFont#1#2#3#4#5{%
  \reset@font\fontsize{#1}{#2pt}%
  \fontfamily{#3}\fontseries{#4}\fontshape{#5}%
  \selectfont}%
\fi\endgroup%
\begin{picture}(10684,8079)(589,-7348)
\put(3376,-4936){\makebox(0,0)[lb]{\smash{{\SetFigFont{14}{16.8}{\rmdefault}{\mddefault}{\updefault}{\color[rgb]{0,0,0}$x$}%
}}}}
\put(2992,143){\rotatebox{315.0}{\makebox(0,0)[lb]{\smash{{\SetFigFont{14}{16.8}{\rmdefault}{\mddefault}{\updefault}{\color[rgb]{0,0,0}$\ell_2(v(x),x)$}%
}}}}}
\put(1405,-2508){\rotatebox{320.0}{\makebox(0,0)[lb]{\smash{{\SetFigFont{14}{16.8}{\rmdefault}{\mddefault}{\updefault}{\color[rgb]{0,0,0}$\ell(v(x),x)$}%
}}}}}
\put(6975,-1389){\rotatebox{310.0}{\makebox(0,0)[lb]{\smash{{\SetFigFont{14}{16.8}{\rmdefault}{\mddefault}{\updefault}{\color[rgb]{0,0,0}$e(x)$}%
}}}}}
\put(1051,-6511){\makebox(0,0)[lb]{\smash{{\SetFigFont{14}{16.8}{\rmdefault}{\mddefault}{\updefault}{\color[rgb]{0,0,0}$\y_1=v(x)$}%
}}}}
\put(11251,-6136){\makebox(0,0)[lb]{\smash{{\SetFigFont{14}{16.8}{\rmdefault}{\mddefault}{\updefault}{\color[rgb]{0,0,0}$\y_2$}%
}}}}
\put(4876,539){\makebox(0,0)[lb]{\smash{{\SetFigFont{14}{16.8}{\rmdefault}{\mddefault}{\updefault}{\color[rgb]{0,0,0}$\y_3$}%
}}}}
\put(3376,-6286){\rotatebox{45.0}{\makebox(0,0)[lb]{\smash{{\SetFigFont{14}{16.8}{\rmdefault}{\mddefault}{\updefault}{\color[rgb]{0,0,0}$d(v(x),\ell_2(v(x),x))=2\,d(v(x),\ell(v(x),x))$}%
}}}}}
\put(960,-5299){\rotatebox{45.0}{\makebox(0,0)[lb]{\smash{{\SetFigFont{12}{14.4}{\rmdefault}{\mddefault}{\updefault}{\color[rgb]{0,0,0} $d(v(x),\ell(v(x),x))$}%
}}}}}
\end{picture}%

%% file: Gammaofnu2.pstex_t
\begin{picture}(0,0)%
\includegraphics{Gammaofnu2.pstex}%
\end{picture}%
\setlength{\unitlength}{3947sp}%
\begingroup\makeatletter\ifx\SetFigFont\undefined%
\gdef\SetFigFont#1#2#3#4#5{%
  \reset@font\fontsize{#1}{#2pt}%
  \fontfamily{#3}\fontseries{#4}\fontshape{#5}%
  \selectfont}%
\fi\endgroup%
\begin{picture}(10684,8043)(589,-7348)
\put(3376,-4936){\makebox(0,0)[lb]{\smash{{\SetFigFont{12}{14.4}{\rmdefault}{\bfdefault}{\updefault}{\color[rgb]{0,0,0}$x$}%
}}}}
\put(2626,-2011){\makebox(0,0)[lb]{\smash{{\SetFigFont{12}{14.4}{\rmdefault}{\mddefault}{\updefault}{\color[rgb]{0,0,0}$\xi_3(x)$}%
}}}}
\put(3867,-6839){\rotatebox{320.0}{\makebox(0,0)[lb]{\smash{{\SetFigFont{12}{14.4}{\rmdefault}{\mddefault}{\updefault}{\color[rgb]{0,0,0}$\xi_1(x)$}%
}}}}}
\put(7533,-5107){\rotatebox{65.0}{\makebox(0,0)[lb]{\smash{{\SetFigFont{12}{14.4}{\rmdefault}{\mddefault}{\updefault}{\color[rgb]{0,0,0}$\xi_2(x)$}%
}}}}}
\put(2101,-7036){\rotatebox{45.0}{\makebox(0,0)[lb]{\smash{{\SetFigFont{12}{14.4}{\rmdefault}{\mddefault}{\updefault}{\color[rgb]{0,0,0}$d(y_1,\xi_1(x))$}%
}}}}}
\put(1405,-2508){\rotatebox{320.0}{\makebox(0,0)[lb]{\smash{{\SetFigFont{12}{14.4}{\rmdefault}{\mddefault}{\updefault}{\color[rgb]{0,0,0}$\ell(y_1,x)$}%
}}}}}
\put(960,-5299){\rotatebox{45.0}{\makebox(0,0)[lb]{\smash{{\SetFigFont{12}{14.4}{\rmdefault}{\mddefault}{\updefault}{\color[rgb]{0,0,0}$d(y_1,\ell(y_1,x))=r\,d(y_1,\xi_1(x))$}%
}}}}}
\put(1051,-6511){\makebox(0,0)[lb]{\smash{{\SetFigFont{12}{14.4}{\rmdefault}{\mddefault}{\updefault}{\color[rgb]{0,0,0}$\y_1$}%
}}}}
\put(4876,539){\makebox(0,0)[lb]{\smash{{\SetFigFont{12}{14.4}{\rmdefault}{\mddefault}{\updefault}{\color[rgb]{0,0,0}$\y_3$}%
}}}}
\put(11251,-6136){\makebox(0,0)[lb]{\smash{{\SetFigFont{12}{14.4}{\rmdefault}{\mddefault}{\updefault}{\color[rgb]{0,0,0}$\y_2$}%
}}}}
\end{picture}%

%% file: means3.pstex_t
\begin{picture}(0,0)%
\includegraphics{means3.pstex}%
\end{picture}%
\setlength{\unitlength}{3947sp}%
\begingroup\makeatletter\ifx\SetFigFont\undefined%
\gdef\SetFigFont#1#2#3#4#5{%
  \reset@font\fontsize{#1}{#2pt}%
  \fontfamily{#3}\fontseries{#4}\fontshape{#5}%
  \selectfont}%
\fi\endgroup%
\begin{picture}(10211,8344)(-5,-9933)
\put(5655,-4426){\rotatebox{20.0}{\makebox(0,0)[lb]{\smash{{\SetFigFont{25}{14.4}{\rmdefault}{\mddefault}{\updefault}{\color[rgb]{0,0,0}$p_{\la}(r)$}%
}}}}}
\put(5490,-3367){\rotatebox{10.0}{\makebox(0,0)[lb]{\smash{{\SetFigFont{25}{14.4}{\rmdefault}{\mddefault}{\updefault}{\color[rgb]{0,0,0}$p(r)$}%
}}}}}
\put(5483,-2686){\rotatebox{5.0}{\makebox(0,0)[lb]{\smash{{\SetFigFont{25}{14.4}{\rmdefault}{\mddefault}{\updefault}{\color[rgb]{0,0,0}$p_{\lo}(r)$}%
}}}}}
\end{picture}%

%% file: var3.pstex_t
\begin{picture}(0,0)%
\includegraphics{var3.pstex}%
\end{picture}%
\setlength{\unitlength}{3947sp}%
\begingroup\makeatletter\ifx\SetFigFont\undefined%
\gdef\SetFigFont#1#2#3#4#5{%
  \reset@font\fontsize{#1}{#2pt}%
  \fontfamily{#3}\fontseries{#4}\fontshape{#5}%
  \selectfont}%
\fi\endgroup%
\begin{picture}(10211,8344)(-5,-9933)
\put(8551,-4786){\makebox(0,0)[lb]{\smash{{\SetFigFont{25}{14.4}{\rmdefault}{\mddefault}{\updefault}{\color[rgb]{0,0,0}$4\,\nu_{\la}(r)$}%
}}}}
\put(6451,-6886){\makebox(0,0)[lb]{\smash{{\SetFigFont{25}{14.4}{\rmdefault}{\mddefault}{\updefault}{\color[rgb]{0,0,0}$4\,\nu(r)$}%
}}}}
\put(5551,-8311){\makebox(0,0)[lb]{\smash{{\SetFigFont{25}{14.4}{\rmdefault}{\mddefault}{\updefault}{\color[rgb]{0,0,0}$4\,\nu_{\lo}(r)$}%
}}}}
\end{picture}%

%% file: ls_lam_cases.pstex_t
\begin{picture}(0,0)%
\includegraphics{ls_lam_cases.pstex}%
\end{picture}%
\setlength{\unitlength}{3947sp}%
\begingroup\makeatletter\ifx\SetFigFont\undefined%
\gdef\SetFigFont#1#2#3#4#5{%
  \reset@font\fontsize{#1}{#2pt}%
  \fontfamily{#3}\fontseries{#4}\fontshape{#5}%
  \selectfont}%
\fi\endgroup%
\begin{picture}(11349,8280)(289,-7498)
\put(10426,-6586){\makebox(0,0)[lb]{\smash{{\SetFigFont{22}{16.8}{\rmdefault}{\mddefault}{\updefault}{\color[rgb]{0,0,0}$\y_2=(1,0)$}%
}}}}
\put(7201,-6586){\makebox(0,0)[lb]{\smash{{\SetFigFont{22}{16.8}{\rmdefault}{\mddefault}{\updefault}{\color[rgb]{0,0,0}$e_3$}%
}}}}
\put(5101,-6586){\makebox(0,0)[lb]{\smash{{\SetFigFont{22}{16.8}{\rmdefault}{\mddefault}{\updefault}{\color[rgb]{0,0,0}$M_3$}%
}}}}
\put(3301,-6586){\makebox(0,0)[lb]{\smash{{\SetFigFont{22}{14.4}{\rmdefault}{\mddefault}{\updefault}{\color[rgb]{0,0,0}$s_1$}%
}}}}
\put(1201,-2611){\makebox(0,0)[lb]{\smash{{\SetFigFont{22}{16.8}{\rmdefault}{\mddefault}{\updefault}{\color[rgb]{0,0,0}$\ell_s(r=4,x)$}%
}}}}
\put(1801,-736){\makebox(0,0)[lb]{\smash{{\SetFigFont{22}{16.8}{\rmdefault}{\mddefault}{\updefault}{\color[rgb]{0,0,0}$\ell_s(r=1.75,x)$}%
}}}}
\put(3226,314){\makebox(0,0)[lb]{\smash{{\SetFigFont{22}{16.8}{\rmdefault}{\mddefault}{\updefault}{\color[rgb]{0,0,0}$\ell_s\bigl(r=\sqrt{2},x\bigr)$}%
}}}}
\put(5926,614){\makebox(0,0)[lb]{\smash{{\SetFigFont{22}{16.8}{\rmdefault}{\mddefault}{\updefault}{\color[rgb]{0,0,0}$\y_3=\bigl(1/2,\sqrt{3}/2\bigr)$}%
}}}}
\put(7876,-2236){\makebox(0,0)[lb]{\smash{{\SetFigFont{22}{16.8}{\rmdefault}{\mddefault}{\updefault}{\color[rgb]{0,0,0}$e_1$}%
}}}}
\put(3976,-6586){\makebox(0,0)[lb]{\smash{{\SetFigFont{22}{14.4}{\rmdefault}{\mddefault}{\updefault}{\color[rgb]{0,0,0}$s_2$}%
}}}}
\put(826,-6736){\makebox(0,0)[lb]{\smash{{\SetFigFont{22}{16.8}{\rmdefault}{\mddefault}{\updefault}{\color[rgb]{0,0,0}$\y_1=(0,0)$}%
}}}}
\put(5851,-4111){\makebox(0,0)[lb]{\smash{{\SetFigFont{22}{16.8}{\rmdefault}{\mddefault}{\updefault}{\color[rgb]{0,0,0}$M_{C}$}%
}}}}
\put(3001,-2761){\makebox(0,0)[lb]{\smash{{\SetFigFont{22}{16.8}{\rmdefault}{\mddefault}{\updefault}{\color[rgb]{0,0,0}$e_2$}%
}}}}
\end{picture}%

%% file: G1ofxCase1.pstex_t
\begin{picture}(0,0)%
\includegraphics{G1ofxCase1.pstex}%
\end{picture}%
\setlength{\unitlength}{3947sp}%
\begingroup\makeatletter\ifx\SetFigFont\undefined%
\gdef\SetFigFont#1#2#3#4#5{%
  \reset@font\fontsize{#1}{#2pt}%
  \fontfamily{#3}\fontseries{#4}\fontshape{#5}%
  \selectfont}%
\fi\endgroup%
\begin{picture}(11349,8130)(289,-7348)
\put(10426,-6586){\makebox(0,0)[lb]{\smash{{\SetFigFont{25}{16.8}{\rmdefault}{\mddefault}{\updefault}{\color[rgb]{0,0,0}$\y_2=(1,0)$}%
}}}}
\put(826,-6586){\makebox(0,0)[lb]{\smash{{\SetFigFont{25}{16.8}{\rmdefault}{\mddefault}{\updefault}{\color[rgb]{0,0,0}$\y_1=(0,0)$}%
}}}}
\put(5326,614){\makebox(0,0)[lb]{\smash{{\SetFigFont{25}{16.8}{\rmdefault}{\mddefault}{\updefault}{\color[rgb]{0,0,0}$\y_3=(1/2,\sqrt{3}/2)$}%
}}}}
\put(7726,-2236){\makebox(0,0)[lb]{\smash{{\SetFigFont{25}{16.8}{\rmdefault}{\mddefault}{\updefault}{\color[rgb]{0,0,0}$e_1$}%
}}}}
\put(3001,-2686){\makebox(0,0)[lb]{\smash{{\SetFigFont{25}{16.8}{\rmdefault}{\mddefault}{\updefault}{\color[rgb]{0,0,0}$e_2$}%
}}}}
\put(7201,-6586){\makebox(0,0)[lb]{\smash{{\SetFigFont{25}{16.8}{\rmdefault}{\mddefault}{\updefault}{\color[rgb]{0,0,0}$e_3$}%
}}}}
\put(5626,-6586){\makebox(0,0)[lb]{\smash{{\SetFigFont{25}{16.8}{\rmdefault}{\mddefault}{\updefault}{\color[rgb]{0,0,0}$M_3$}%
}}}}
\put(2026,-5986){\makebox(0,0)[lb]{\smash{{\SetFigFont{25}{16.8}{\rmdefault}{\mddefault}{\updefault}{\color[rgb]{0,0,0}$\xi_1(r,x)$}%
}}}}
\put(5701,-4111){\makebox(0,0)[lb]{\smash{{\SetFigFont{25}{16.8}{\rmdefault}{\mddefault}{\updefault}{\color[rgb]{0,0,0}$M_{C}$}%
}}}}
\put(1426,-6211){\makebox(0,0)[lb]{\smash{{\SetFigFont{25}{16.8}{\rmdefault}{\mddefault}{\updefault}{\color[rgb]{0,0,0}$x_1$}%
}}}}
\end{picture}%

%% file: G1ofxCase2.pstex_t
\begin{picture}(0,0)%
\includegraphics{G1ofxCase2.pstex}%
\end{picture}%
\setlength{\unitlength}{3947sp}%
\begingroup\makeatletter\ifx\SetFigFont\undefined%
\gdef\SetFigFont#1#2#3#4#5{%
  \reset@font\fontsize{#1}{#2pt}%
  \fontfamily{#3}\fontseries{#4}\fontshape{#5}%
  \selectfont}%
\fi\endgroup%
\begin{picture}(11349,8130)(289,-7348)
\put(10426,-6586){\makebox(0,0)[lb]{\smash{{\SetFigFont{25}{16.8}{\rmdefault}{\mddefault}{\updefault}{\color[rgb]{0,0,0}$\y_2=(1,0)$}%
}}}}
\put(826,-6586){\makebox(0,0)[lb]{\smash{{\SetFigFont{25}{16.8}{\rmdefault}{\mddefault}{\updefault}{\color[rgb]{0,0,0}$\y_1=(0,0)$}%
}}}}
\put(5326,614){\makebox(0,0)[lb]{\smash{{\SetFigFont{25}{16.8}{\rmdefault}{\mddefault}{\updefault}{\color[rgb]{0,0,0}$\y_3=(1/2,\sqrt{3}/2)$}%
}}}}
\put(7726,-2236){\makebox(0,0)[lb]{\smash{{\SetFigFont{25}{16.8}{\rmdefault}{\mddefault}{\updefault}{\color[rgb]{0,0,0}$e_1$}%
}}}}
\put(3001,-2686){\makebox(0,0)[lb]{\smash{{\SetFigFont{25}{16.8}{\rmdefault}{\mddefault}{\updefault}{\color[rgb]{0,0,0}$e_2$}%
}}}}
\put(7201,-6586){\makebox(0,0)[lb]{\smash{{\SetFigFont{25}{16.8}{\rmdefault}{\mddefault}{\updefault}{\color[rgb]{0,0,0}$e_3$}%
}}}}
\put(5626,-6586){\makebox(0,0)[lb]{\smash{{\SetFigFont{25}{16.8}{\rmdefault}{\mddefault}{\updefault}{\color[rgb]{0,0,0}$M_3$}%
}}}}
\put(2026,-5986){\makebox(0,0)[lb]{\smash{{\SetFigFont{25}{16.8}{\rmdefault}{\mddefault}{\updefault}{\color[rgb]{0,0,0}$\xi_1(r,x)$}%
}}}}
\put(5701,-4111){\makebox(0,0)[lb]{\smash{{\SetFigFont{25}{16.8}{\rmdefault}{\mddefault}{\updefault}{\color[rgb]{0,0,0}$M_{C}$}%
}}}}
\put(3601,-6136){\makebox(0,0)[lb]{\smash{{\SetFigFont{25}{16.8}{\rmdefault}{\mddefault}{\updefault}{\color[rgb]{0,0,0}$x_1$}%
}}}}
\put(6226,-4636){\makebox(0,0)[lb]{\smash{{\SetFigFont{25}{16.8}{\rmdefault}{\mddefault}{\updefault}{\color[rgb]{0,0,0}$\xi_2(r,x)$}%
}}}}
\end{picture}%

%% file: G1ofxCase3.pstex_t
\begin{picture}(0,0)%
\includegraphics{G1ofxCase3.pstex}%
\end{picture}%
\setlength{\unitlength}{3947sp}%
\begingroup\makeatletter\ifx\SetFigFont\undefined%
\gdef\SetFigFont#1#2#3#4#5{%
  \reset@font\fontsize{#1}{#2pt}%
  \fontfamily{#3}\fontseries{#4}\fontshape{#5}%
  \selectfont}%
\fi\endgroup%
\begin{picture}(11349,8130)(289,-7348)
\put(10426,-6586){\makebox(0,0)[lb]{\smash{{\SetFigFont{25}{16.8}{\rmdefault}{\mddefault}{\updefault}{\color[rgb]{0,0,0}$\y_2=(1,0)$}%
}}}}
\put(826,-6586){\makebox(0,0)[lb]{\smash{{\SetFigFont{25}{16.8}{\rmdefault}{\mddefault}{\updefault}{\color[rgb]{0,0,0}$\y_1=(0,0)$}%
}}}}
\put(5326,614){\makebox(0,0)[lb]{\smash{{\SetFigFont{25}{16.8}{\rmdefault}{\mddefault}{\updefault}{\color[rgb]{0,0,0}$\y_3=(1/2,\sqrt{3}/2)$}%
}}}}
\put(7726,-2236){\makebox(0,0)[lb]{\smash{{\SetFigFont{25}{16.8}{\rmdefault}{\mddefault}{\updefault}{\color[rgb]{0,0,0}$e_1$}%
}}}}
\put(3001,-2686){\makebox(0,0)[lb]{\smash{{\SetFigFont{25}{16.8}{\rmdefault}{\mddefault}{\updefault}{\color[rgb]{0,0,0}$e_2$}%
}}}}
\put(7201,-6586){\makebox(0,0)[lb]{\smash{{\SetFigFont{25}{16.8}{\rmdefault}{\mddefault}{\updefault}{\color[rgb]{0,0,0}$e_3$}%
}}}}
\put(5626,-6586){\makebox(0,0)[lb]{\smash{{\SetFigFont{25}{16.8}{\rmdefault}{\mddefault}{\updefault}{\color[rgb]{0,0,0}$M_3$}%
}}}}
\put(5701,-4111){\makebox(0,0)[lb]{\smash{{\SetFigFont{25}{16.8}{\rmdefault}{\mddefault}{\updefault}{\color[rgb]{0,0,0}$M_{C}$}%
}}}}
\put(7726,-5236){\makebox(0,0)[lb]{\smash{{\SetFigFont{25}{16.8}{\rmdefault}{\mddefault}{\updefault}{\color[rgb]{0,0,0}$\xi_2(r,x)$}%
}}}}
\put(4651,-6136){\makebox(0,0)[lb]{\smash{{\SetFigFont{25}{16.8}{\rmdefault}{\mddefault}{\updefault}{\color[rgb]{0,0,0}$x_1$}%
}}}}
\put(2851,-5536){\makebox(0,0)[lb]{\smash{{\SetFigFont{25}{16.8}{\rmdefault}{\mddefault}{\updefault}{\color[rgb]{0,0,0}$\xi_1(r,x)$}%
}}}}
\end{picture}%

%% file: G1ofxCase4.pstex_t
\begin{picture}(0,0)%
\includegraphics{G1ofxCase4.pstex}%
\end{picture}%
\setlength{\unitlength}{3947sp}%
\begingroup\makeatletter\ifx\SetFigFont\undefined%
\gdef\SetFigFont#1#2#3#4#5{%
  \reset@font\fontsize{#1}{#2pt}%
  \fontfamily{#3}\fontseries{#4}\fontshape{#5}%
  \selectfont}%
\fi\endgroup%
\begin{picture}(11349,8130)(289,-7348)
\put(10426,-6586){\makebox(0,0)[lb]{\smash{{\SetFigFont{25}{16.8}{\rmdefault}{\mddefault}{\updefault}{\color[rgb]{0,0,0}$\y_2=(1,0)$}%
}}}}
\put(826,-6586){\makebox(0,0)[lb]{\smash{{\SetFigFont{25}{16.8}{\rmdefault}{\mddefault}{\updefault}{\color[rgb]{0,0,0}$\y_1=(0,0)$}%
}}}}
\put(5326,614){\makebox(0,0)[lb]{\smash{{\SetFigFont{25}{16.8}{\rmdefault}{\mddefault}{\updefault}{\color[rgb]{0,0,0}$\y_3=(1/2,\sqrt{3}/2)$}%
}}}}
\put(7201,-6586){\makebox(0,0)[lb]{\smash{{\SetFigFont{25}{16.8}{\rmdefault}{\mddefault}{\updefault}{\color[rgb]{0,0,0}$e_3$}%
}}}}
\put(2026,-5986){\makebox(0,0)[lb]{\smash{{\SetFigFont{25}{16.8}{\rmdefault}{\mddefault}{\updefault}{\color[rgb]{0,0,0}$\xi_1(r,x)$}%
}}}}
\put(5701,-4111){\makebox(0,0)[lb]{\smash{{\SetFigFont{25}{16.8}{\rmdefault}{\mddefault}{\updefault}{\color[rgb]{0,0,0}$M_{C}$}%
}}}}
\put(3826,-5536){\makebox(0,0)[lb]{\smash{{\SetFigFont{25}{16.8}{\rmdefault}{\mddefault}{\updefault}{\color[rgb]{0,0,0}$x_1$}%
}}}}
\put(3901,-1111){\makebox(0,0)[lb]{\smash{{\SetFigFont{25}{16.8}{\rmdefault}{\mddefault}{\updefault}{\color[rgb]{0,0,0}$e_2$}%
}}}}
\put(6901,-1111){\makebox(0,0)[lb]{\smash{{\SetFigFont{25}{16.8}{\rmdefault}{\mddefault}{\updefault}{\color[rgb]{0,0,0}$e_1$}%
}}}}
\put(5326,-6586){\makebox(0,0)[lb]{\smash{{\SetFigFont{25}{16.8}{\rmdefault}{\mddefault}{\updefault}{\color[rgb]{0,0,0}$M_3$}%
}}}}
\put(3151,-6661){\makebox(0,0)[lb]{\smash{{\SetFigFont{25}{16.8}{\rmdefault}{\mddefault}{\updefault}{\color[rgb]{0,0,0}$G_1$}%
}}}}
\put(1501,-4936){\makebox(0,0)[lb]{\smash{{\SetFigFont{25}{16.8}{\rmdefault}{\mddefault}{\updefault}{\color[rgb]{0,0,0}$G_6$}%
}}}}
\put(2851,-2686){\makebox(0,0)[lb]{\smash{{\SetFigFont{25}{16.8}{\rmdefault}{\mddefault}{\updefault}{\color[rgb]{0,0,0}$M_2$}%
}}}}
\put(3976,-2986){\makebox(0,0)[lb]{\smash{{\SetFigFont{25}{16.8}{\rmdefault}{\mddefault}{\updefault}{\color[rgb]{0,0,0}$L_5$}%
}}}}
\put(5251,-2986){\makebox(0,0)[lb]{\smash{{\SetFigFont{25}{16.8}{\rmdefault}{\mddefault}{\updefault}{\color[rgb]{0,0,0}$\xi_3(r,x)$}%
}}}}
\put(6901,-2986){\makebox(0,0)[lb]{\smash{{\SetFigFont{25}{16.8}{\rmdefault}{\mddefault}{\updefault}{\color[rgb]{0,0,0}$L_4$}%
}}}}
\put(7651,-3061){\makebox(0,0)[lb]{\smash{{\SetFigFont{25}{16.8}{\rmdefault}{\mddefault}{\updefault}{\color[rgb]{0,0,0}$L_3$}%
}}}}
\put(5626,-6211){\makebox(0,0)[lb]{\smash{{\SetFigFont{25}{16.8}{\rmdefault}{\mddefault}{\updefault}{\color[rgb]{0,0,0}$L_2$}%
}}}}
\put(6526,-4786){\makebox(0,0)[lb]{\smash{{\SetFigFont{25}{16.8}{\rmdefault}{\mddefault}{\updefault}{\color[rgb]{0,0,0}$\xi_2(r,x)$}%
}}}}
\put(7951,-2686){\makebox(0,0)[lb]{\smash{{\SetFigFont{25}{16.8}{\rmdefault}{\mddefault}{\updefault}{\color[rgb]{0,0,0}$M_1$}%
}}}}
\end{picture}%

%% file: G1ofxCase5.pstex_t
\begin{picture}(0,0)%
\includegraphics{G1ofxCase5.pstex}%
\end{picture}%
\setlength{\unitlength}{3947sp}%
\begingroup\makeatletter\ifx\SetFigFont\undefined%
\gdef\SetFigFont#1#2#3#4#5{%
  \reset@font\fontsize{#1}{#2pt}%
  \fontfamily{#3}\fontseries{#4}\fontshape{#5}%
  \selectfont}%
\fi\endgroup%
\begin{picture}(11349,8130)(289,-7348)
\put(10426,-6586){\makebox(0,0)[lb]{\smash{{\SetFigFont{25}{16.8}{\rmdefault}{\mddefault}{\updefault}{\color[rgb]{0,0,0}$\y_2=(1,0)$}%
}}}}
\put(826,-6586){\makebox(0,0)[lb]{\smash{{\SetFigFont{25}{16.8}{\rmdefault}{\mddefault}{\updefault}{\color[rgb]{0,0,0}$\y_1=(0,0)$}%
}}}}
\put(5326,614){\makebox(0,0)[lb]{\smash{{\SetFigFont{25}{16.8}{\rmdefault}{\mddefault}{\updefault}{\color[rgb]{0,0,0}$\y_3=(1/2,\sqrt{3}/2)$}%
}}}}
\put(7726,-2236){\makebox(0,0)[lb]{\smash{{\SetFigFont{25}{16.8}{\rmdefault}{\mddefault}{\updefault}{\color[rgb]{0,0,0}$e_1$}%
}}}}
\put(3001,-2686){\makebox(0,0)[lb]{\smash{{\SetFigFont{25}{16.8}{\rmdefault}{\mddefault}{\updefault}{\color[rgb]{0,0,0}$e_2$}%
}}}}
\put(7201,-6586){\makebox(0,0)[lb]{\smash{{\SetFigFont{25}{16.8}{\rmdefault}{\mddefault}{\updefault}{\color[rgb]{0,0,0}$e_3$}%
}}}}
\put(5626,-6586){\makebox(0,0)[lb]{\smash{{\SetFigFont{25}{16.8}{\rmdefault}{\mddefault}{\updefault}{\color[rgb]{0,0,0}$M_3$}%
}}}}
\put(5701,-4111){\makebox(0,0)[lb]{\smash{{\SetFigFont{25}{16.8}{\rmdefault}{\mddefault}{\updefault}{\color[rgb]{0,0,0}$M_{C}$}%
}}}}
\put(4801,-2986){\makebox(0,0)[lb]{\smash{{\SetFigFont{25}{16.8}{\rmdefault}{\mddefault}{\updefault}{\color[rgb]{0,0,0}$\xi_3(r,x)$}%
}}}}
\put(4951,-5386){\makebox(0,0)[lb]{\smash{{\SetFigFont{25}{16.8}{\rmdefault}{\mddefault}{\updefault}{\color[rgb]{0,0,0}$x_1$}%
}}}}
\put(7576,-5236){\makebox(0,0)[lb]{\smash{{\SetFigFont{25}{16.8}{\rmdefault}{\mddefault}{\updefault}{\color[rgb]{0,0,0}$\xi_2(r,x)$}%
}}}}
\put(3076,-5536){\makebox(0,0)[lb]{\smash{{\SetFigFont{25}{16.8}{\rmdefault}{\mddefault}{\updefault}{\color[rgb]{0,0,0}$\xi_1(r,x)$}%
}}}}
\end{picture}%

%% file: G1ofxCase6.pstex_t
\begin{picture}(0,0)%
\includegraphics{G1ofxCase6.pstex}%
\end{picture}%
\setlength{\unitlength}{3947sp}%
\begingroup\makeatletter\ifx\SetFigFont\undefined%
\gdef\SetFigFont#1#2#3#4#5{%
  \reset@font\fontsize{#1}{#2pt}%
  \fontfamily{#3}\fontseries{#4}\fontshape{#5}%
  \selectfont}%
\fi\endgroup%
\begin{picture}(11349,8130)(289,-7348)
\put(10426,-6586){\makebox(0,0)[lb]{\smash{{\SetFigFont{25}{16.8}{\rmdefault}{\mddefault}{\updefault}{$\y_2=(1,0)$}%
}}}}
\put(826,-6586){\makebox(0,0)[lb]{\smash{{\SetFigFont{25}{16.8}{\rmdefault}{\mddefault}{\updefault}{$\y_1=(0,0)$}%
}}}}
\put(5326,614){\makebox(0,0)[lb]{\smash{{\SetFigFont{25}{16.8}{\rmdefault}{\mddefault}{\updefault}{$\y_3=(1/2,\sqrt{3}/2)$}%
}}}}
\put(7726,-2236){\makebox(0,0)[lb]{\smash{{\SetFigFont{25}{16.8}{\rmdefault}{\mddefault}{\updefault}{$e_1$}%
}}}}
\put(3001,-2686){\makebox(0,0)[lb]{\smash{{\SetFigFont{25}{16.8}{\rmdefault}{\mddefault}{\updefault}{$e_2$}%
}}}}
\put(7201,-6586){\makebox(0,0)[lb]{\smash{{\SetFigFont{25}{16.8}{\rmdefault}{\mddefault}{\updefault}{$e_3$}%
}}}}
\put(5626,-6586){\makebox(0,0)[lb]{\smash{{\SetFigFont{25}{16.8}{\rmdefault}{\mddefault}{\updefault}{$M_3$}%
}}}}
\put(5701,-4111){\makebox(0,0)[lb]{\smash{{\SetFigFont{25}{16.8}{\rmdefault}{\mddefault}{\updefault}{$M_{C}$}%
}}}}
\put(5251,-4411){\makebox(0,0)[lb]{\smash{{\SetFigFont{25}{16.8}{\rmdefault}{\mddefault}{\updefault}{$x_1$}%
}}}}
\put(5851,-4636){\makebox(0,0)[lb]{\smash{{\SetFigFont{25}{16.8}{\rmdefault}{\mddefault}{\updefault}{$\xi_2(r,x)$}%
}}}}
\put(5176,-3436){\makebox(0,0)[lb]{\smash{{\SetFigFont{25}{16.8}{\rmdefault}{\mddefault}{\updefault}{$\xi_3(r,x)$}%
}}}}
\put(4201,-4711){\makebox(0,0)[lb]{\smash{{\SetFigFont{25}{16.8}{\rmdefault}{\mddefault}{\updefault}{$\xi_1(r,x)$}%
}}}}
\put(5476,-6961){\makebox(0,0)[lb]{\smash{{\SetFigFont{25}{14.4}{\rmdefault}{\mddefault}{\updefault}{case-6}%
}}}}
\end{picture}%